\documentstyle{amsart}
\font\tenmsb=msbm10
\font\sevenmsb=msbm7
\font\fivemsb=msbm5
\font\largemsb=msbm10 at11pt
\font\largesevenmsb=msbm7 at 7.7pt
\font\largefivemsb=msbm5 at 5.5pt
\font\smallmsb=msbm7 at7.7pt
\newfam\msbfam
\textfont\msbfam=\tenmsb
\scriptfont\msbfam=\sevenmsb
\scriptscriptfont\msbfam=\fivemsb
\newfam\largemsbfam
\textfont\largemsbfam=\largemsb
\scriptfont\largemsbfam=\largesevenmsb
\scriptscriptfont\largemsbfam=\largefivemsb
\newfam\smallmsbfam
\textfont\smallmsbfam=\smallmsb
\def\Bbb#1{\ifdim1em>10.1pt{\fam\largemsbfam\relax#1\kern.9pt}
   \else\ifdim1em<9.9pt{\fam\smallmsbfam\relax#1\kern.5pt}
   \else\fam\msbfam\relax#1\kern.8pt\fi\fi}
\def\qed{\hfill\vbox{\hrule height.09ex
   \hbox{\vrule width.09ex height1.7ex depth.3ex \kern2ex
   \vrule width.09ex height1.7ex depth.3ex}\hrule height.09ex}\bigskip}
\def\B{{\Bbb B}}
\def\Z{{\Bbb Z}}
\def\R{{\Bbb R}}
\def\C{{\Bbb C}}
\def\S{{\Bbb S}}

\def\M{{\Bbb M}}

\def\interior{\hspace{0.05in}
\begin{picture}(6,6)
\put (0,0){\line(1,0){6}}
\put (6,0){\line(0,1){6}}
\end{picture} \hspace{0.05in}}
\newcommand{\dis}{\displaystyle}

\newtheorem{theorem}{Theorem}
\newtheorem{proposition}{Proposition}[section]
\newtheorem{lemma}[proposition]{Lemma}

\begin{document}

\large

\title[Global $\bar\partial_{\bold M}$-homotopy with $C^k$ estimates]
{Global $\bar\partial_{\bold M}$-homotopy with $C^k$ estimates
for a family\\ of compact, regular $q$-pseudoconcave CR manifolds}
\author[Peter L. Polyakov]{Peter L. Polyakov}
\address{Department of Mathematics, University of Wyoming, Laramie, WY 82071}
\email{polyakov@@uwyo.edu}
\today
\maketitle

\begin{abstract}
Let ${\bold M}_0$ be a compact, regular q-pseudoconcave compact CR submanifold
of a complex manifold ${\bold G}$ and ${\cal B}$ - a holomorphic vector bundle
on ${\bold G}$ such that $\dim H^r\left({\bold M}_0,
{\cal B}\big|_{\bold M}\right)=0$ for some fixed $r<q$. We prove a
global homotopy formula with $C^k$ estimates for $r$-cohomology of
${\cal B}$  on arbitrary CR submanifold ${\bold M}$ close enough to
${\bold M}_0$.
\end{abstract}

\section{Introduction.}\label{Introduction}

\indent
Let ${\bold M}$ be a compact generic CR submanifold in a complex
$n$ - dimensional manifold ${\bold G}$, i.e. for any $z \in {\bold M}$
there exist a neighborhood $V \ni z$ in ${\bold G}$
and smooth real valued functions
$$\{ \rho_{k},  \ k = 1, \dots , m  \ ( 1 < m < n-1)\}$$
on $V$ such that
$$\begin{array}{ll}
{\bold M} \cap V = \{ z \in {\bold G} \cap V: \rho_{1}(z) =
\dots = \rho_{m}(z) = 0\},\vspace{0.1in}\\
\partial\rho_1\wedge\cdots\wedge\partial\rho_m\neq 0 
\ \ \mbox{on} \ {\bold M} \cap V.
\end{array}
\eqno(\arabic{equation})
\newcounter{manifold}
\setcounter{manifold}{\value{equation}}
\addtocounter{equation}{1}$$
\indent
In \cite{P2}, where our motivation was to obtain sharp
estimates for the operator $\bar\partial_{\bold M}$, we introduced
a special nonhomogeneous Lipschitz scale $\Pi^{\alpha}({\bold M})$,
based on the Stein's scale of spaces $\Gamma^{\alpha}({\bold M})$
(cf. \cite{S}, \cite{FS}), and proved an ``almost homotopy'' in this
scale. The sharpness of estimates in the scale $\Pi^{\alpha}({\bold M})$
is expressed by the fact that $\bar\partial_{\bold M}$ maps
$\Pi^{\alpha}({\bold M})$ into $\Pi^{\alpha-1}({\bold M})$ and the
solution operators map $\Pi^{\alpha}({\bold M})$ into
$\Pi^{\alpha+1}({\bold M})$.\\
\indent
In this paper our motivation is not the sharpness of estimates but rather
stability of the global solvability of the $\bar\partial_{\bold M}$-equation
on a compact CR submanifold of a complex manifold. Namely, we prove a global
homotopy formula with $C^k$- estimates for operator $\bar\partial_{\bold M}$
on a family of compact CR submanifolds, close to a fixed submanifold
${\bold M}_0$ and satisfying special concavity condition. 
The proof basically consists of two steps. On the first step,
using a simplified version of techniques from \cite{P2} we construct
a local ``almost homotopy'' formula with estimates and on the second step
we globalize the constructed formula with control of estimates. Estimates
in Propositions~\ref{REstimates} and ~\ref{Rk+sEstimates} represent an important
byproduct. They provide local ``tame estimates'' in terminology of R. Hamilton
(cf. \cite{Ha}) for operator $\bar\partial_{\bold M}$ and are used in \cite{P3} to
study deformations of a fixed embedded CR structure on ${\bold M}_0$.\\
\indent
Before formulating the main result we will introduce necessary
notations and definitions.\\
\indent
The CR structure on ${\bold M}$ is induced from ${\bold G}$ and
is defined by the subbundles
$$T^{\prime\prime}({\bold M})
= T^{\prime\prime}({\bold G})|_{\bold M} \cap {\bold C}T({\bold M})
\hspace{0.2in} \mbox{and} \hspace{0.2in}
T^{\prime}({\bold M}) = T^{\prime}({\bold G})|_{\bold M}
\cap {\bold C}T({\bold M}),$$
where ${\bold C}T({\bold M})$ is the complexified tangent bundle
of ${\bold M}$ and the subbundles $T^{\prime\prime}({\bold G})$ and
$T^{\prime}({\bold G})= \overline{T^{\prime\prime}}({\bold G})$ of
the complexified tangent bundle ${\bold C}T({\bold G})$ define the complex
structure on ${\bold G}$.\\
\indent
We will denote by $T^c({\bold M})$ the subbundle $T({\bold M}) \cap
\left[ T^{\prime}({\bold M}) \oplus T^{\prime\prime}({\bold M})\right].$
If we fix a hermitian scalar product on ${\bold G}$ then we can choose
a subbundle $N \in T({\bold M})$ of real dimension $m$ such that
$T^c({\bold M}) \perp N$ and for a complex subbundle 
${\bold N} = {\bold C}N$ of ${\bold C}T({\bold M})$ we have
$${\bold C}T({\bold M}) =
T^{\prime}({\bold M}) \oplus T^{\prime\prime}({\bold M})
\oplus {\bold N},
\hspace{0.1in}T^{\prime}({\bold M})
\perp {\bold N} \hspace{0.1in}\mbox{and}
\hspace{0.1in}T^{\prime\prime}({\bold M}) \perp {\bold N}.$$
\indent
The Levi form of ${\bold M}$ is defined as the hermitian form on
$T^{\prime}({\bold M})$ with values in ${\bold N}$
$${\cal L}_z(L(z)) = \sqrt{-1} \cdot \pi
\left( \left[\overline L, L\right]\right)(z)
\hspace{0.2in}\left( L(z) \in T^{\prime}_z({\bold M})\right),$$
where $\left[ \overline L, L \right] =
\overline L L - L \overline L$ and $\pi$ is the projection
of ${\bold C}T({\bold M})$ along
$T^{\prime}({\bold M}) \oplus T^{\prime\prime}({\bold M})$
onto $N$.\\
\indent
If functions $\{ \rho_{k} \}$ are chosen so that the vectors
$\{ \mbox{grad} \rho_{k} \}$ are orthonormal then 
the Levi form of ${\bold M}$ may be defined as
$$L_z({\bold M}) = - \sum_{k=1}^m\left(L_{z}\rho_k(\zeta)\right) \cdot 
\mbox{grad $\rho_k(z)$},$$
where $L_{z}\rho(\zeta)$ is the Levi form
of the real valued function $\rho \in C^4({\bold D})$ at the point $z$:
$$L_{z}\rho(\zeta) = \sum_{i,j}\frac{{\partial}^2\rho}
{\partial \zeta_{i} \partial\bar \zeta_{j}}(z) \ \zeta_{i} 
\cdot \bar{\zeta_{j}}.$$
\indent
For a pair of vectors $\mu = (\mu_1, \dots, \mu_n )$ and
$\nu = (\nu_1, \dots, \nu_n )$ in $\C^n$ we will denote
$\langle \mu, \nu \rangle = \sum_{i=1}^n \mu_i \cdot \nu_i$.\\
\indent
For a unit vector
$\theta = (\theta_1, \dots \theta_m) \in \mbox{Re}{\bold N}_z$ we
define the Levi form of 
${\bold M}$ at the point $z \in {\bold M}$ in the direction 
$\theta$ as the scalar hermitian form on
${\bold C}T^c_{z}({\bold M})$
$$\langle \theta,\ L_z({\bold M}) \rangle = - L_{z}\rho_\theta(\zeta),$$
where $\rho_\theta(\zeta) = \sum_{k=1}^m \theta_k \rho_k(\zeta)$.\\
\indent
Following \cite{H2} we call ${\bold M}$ q-pseudoconcave
at $z \in {\bold M}$ in the direction $\theta$ if the Levi form
of ${\bold M}$ at $z$ in this direction 
$\langle \theta,\ L_z({\bold M}) \rangle$ has at least $q$ negative 
eigenvalues on ${\bold C}T^c_{z}({\bold M})$. Correspondingly ${\bold M}$
is called q-pseudoconcave at $z \in {\bold M}$ if it
is q-pseudoconcave in all directions.\\
\indent
We call a q-pseudoconcave CR manifold ${\bold M}$ a
regular q-pseudoconcave CR manifold (cf. \cite{P1})
if for any $z \in {\bold M}$ there exist an open neighborhood ${\cal U} \ni z$
in  ${\bold M}$ and a family $E_{q}(\theta,z)$ of
$q$-dimensional complex linear subspaces in ${\bold C}T^c_{z}({\bold M})$
smoothly depending on $(\theta, z) \in \S^{m-1} \times {\cal U}$ and such that
the Levi form $\langle \theta,\ L_z({\bold M}) \rangle$
is strictly negative on $E_{q}(\theta,z).$\\
\indent
For a function $f$ on ${\bold M}$ ($h$ on ${\bold G}$)
we denote by $\left| f \right|_k$ (respectively $\left| h \right|_k$)
the $C^k\left({\bold M}\right)$-norm of $f$ (respectively the
$C^k\left({\bold G}\right)$-norm of $h$).\\
\indent
Let ${\cal E}_0:\M \to {\bold M}_0 \subset {\bold G}$ be a fixed
embedding of a compact $C^p$ manifold $\M$ into a complex manifold
${\bold G}$ such that ${\bold M}_0 = {\cal E}_0\left(\M\right)$ is
a regular q-pseudoconcave CR submanifold of ${\bold G}$ and let
$\left\{{\cal U}^{\iota}\right\}_1^N$ be a finite cover of some
neighborhood of ${\bold M}_0$ in ${\bold G}$ such that in each
${\cal U}^{\iota}$ the manifold ${\bold M}_0 \cap {\cal U}^{\iota}$ has
the form (\arabic{manifold}) with defining functions
$\left\{\rho^{(0)}_{\iota,l}\right\}_{1\leq l \leq m}$.
If ${\cal E}$ is an embedding of $\M$ into ${\bold G}$ close to
${\cal E}_0$ with ${\bold M} = {\cal E}\left(\M\right)$ 
then the map ${\cal F}={\cal E}\circ{\cal E}_0^{-1}:{\bold M}_0\to {\bold M}$
may be defined in some small enough neighborhood
$U^{\iota}={\bold M}_0\cap {\cal U}^{\iota}$ as ${\cal F}(z) = z+f_{\iota}(z)$,
with $f_{\iota} \in \left[C^p(U^{\iota})\right]^n$.
If $\left|f\right|_p = \max_{\iota}\{\left|f_{\iota}\right|_p\}$
is small enough then ${\cal G} = {\cal F}^{-1}:{\bold M}\to {\bold M}_0$
is also well defined and has the form ${\cal G}(z) = z+g_{\iota}(z)$
in some neighborhood $V^{\iota} \subset {\cal F}\left(U^{\iota}\right)$ with
$g_{\iota} \in \left[C^p\left(V^{\iota}\right)\right]^n$. We denote
$$|{\cal E}|_p = \max_{1\leq \iota \leq N}\left\{|g_{\iota}|_p\right\}.$$
\indent
For a compact, regular q-pseudoconcave $C^p$ submanifold
${\bold M} \subset {\bold G}$ and a holomorphic
vector bundle ${\cal B}$ on ${\bold G}$ we say that
$\dim H^r\left({\bold M},{\cal B}\big|_{\bold M}\right)=0$
for $r \in \Z^{+}$ if for any $\bar\partial_{\bold M}$-closed form
$f \in C^{k}_{(0,r)}\left({\bold M},{\cal B}\big|_{\bold M}\right)$,
$(1\leq k \leq p)$, there exists a form $h \in C^{k}_{(0,r-1)}
\left({\bold M},{\cal B}\big|_{\bold M}\right)$ such that
$\bar\partial_{\bold M}h = f$.\\
\indent
The following theorem represents the main result of the paper.

\begin{theorem}\label{HomotopyTheorem}
Let ${\bold M}_0 \subset {\bold G}$ be a compact, regular
q-pseudoconcave $C^p$ submanifold and let ${\cal B}$ be a holomorphic
vector bundle on ${\bold G}$. Let for some fixed
$1 \leq r \leq q-1$ the condition $\dim H^r\left({\bold M}_0,
{\cal B}\big|_{{\bold M}_0}\right)=0$ be satisfied.
Then there exist $C,\delta > 0$ such that for any $C^p$ embedding
$${\cal E}:\M\to {\bold M} \subset {\bold G}$$
with $|{\cal E}|_p < \delta$ and $k \leq p-7$ there exist linear
bounded operators
$${\bold Q}^i_{\bold M}: C^{k}_{(0,i)}
\left({\bold M},{\cal B}\big|_{\bold M}\right)
\rightarrow C^{k-1}_{(0,i-1)}
\left({\bold M},{\cal B}\big|_{\bold M}\right)
\hspace{0.05in}\mbox{for}\hspace{0.05in}i=r,r+1$$
such that for any differential form $h \in C^{k}_{(0,r)}({\bold M})$
equality:
$$h = \bar\partial_{\bold M} {\bold Q}^r_{\bold M}(h)
+{\bold Q}^{r+1}_{\bold M}(\bar\partial_{\bold M} h)
\eqno(\arabic{equation})
\newcounter{HomotopyFormula}
\setcounter{HomotopyFormula}{\value{equation}}
\addtocounter{equation}{1}$$
and estimates
$$\left|{\bold Q}^i_{\bold M}(h)\right|_{k-1}\leq C|h|_k
\eqno(\arabic{equation})
\newcounter{HomotopyEstimate}
\setcounter{HomotopyEstimate}{\value{equation}}
\addtocounter{equation}{1}$$
hold.
\end{theorem}

\indent
Author thanks G. Henkin, S. Krantz and A. Tumanov for helpful discussions and
the referee for pointing out a mistake in the original version of this article.\\

\section{Construction of "almost-homotopy" formula.}
\label{Construction}

\indent
As an intermediate step in the proof of Theorem~\ref{HomotopyTheorem} we
prove the following proposition.

\begin{proposition}\label{AlmostHomotopy}
Let ${\bold M} \subset {\bold G}$ be a compact, regular
q-pseudoconcave $C^p$ submanifold and let
$\left\{{\cal U}^{\iota}\right\}_1^N$ be a finite cover of some
neighborhood of ${\bold M}$ in ${\bold G}$ such that in each
${\cal U}^{\iota}$ the manifold ${\bold M} \cap {\cal U}^{\iota}$ has
the form (\arabic{manifold}) with defining functions
$\left\{\rho_{\iota,l}\right\}_1^m$.
Then for any $r = 1,\dots, q-1$ and $k \leq p-4$ there exist
linear bounded operators
$${\bold R}^r_{\bold M}: C^{k}_{(0,r)}\left({\bold M},
{\cal B}\big|_{\bold M}\right) \rightarrow
C^{k+1/2}_{(0,r-1)}\left({\bold M},{\cal B}\big|_{\bold M}
\right)\hspace{0.1in}\mbox{and}\hspace{0.1in}
{\bold H}^r_{\bold M}: C^{k}_{(0,r)}\left({\bold M},
{\cal B}\big|_{\bold M}\right) \rightarrow
C^{k+1/2}_{(0,r)}\left({\bold M},{\cal B}\big|_{\bold M}\right)$$
such that for any differential form $h \in
C^{k}_{(0,r)}\left({\bold M},{\cal B}\big|_{\bold M}\right)$
we have equality:
$$h = \bar\partial_{\bold M} {\bold R}^r_{\bold M}(h) +
{\bold R}^{r+1}_{\bold M}(\bar\partial_{\bold M} h)
+ {\bold H}^r_{\bold M}(h)
\eqno(\arabic{equation})
\newcounter{AlmostFormula}
\setcounter{AlmostFormula}{\value{equation}}
\addtocounter{equation}{1}$$
and estimates
$$\left|{\bold R}^r_{\bold M}\right|_k
\leq C(k)\left(1 + |\rho|_{k+3}\right)^{P(k)},\hspace{0.1in}
\left|{\bold H}^r_{\bold M}\right|_k
\leq C(k)\left(1 + |\rho|_{k+3}\right)^{P(k)},
\eqno(\arabic{equation})
\newcounter{AlmostEstimate1}
\setcounter{AlmostEstimate1}{\value{equation}}
\addtocounter{equation}{1}$$
$$\left|{\bold R}^r_{\bold M}\right|_{k+1/2}
\leq C(k)\left(1 + |\rho|_{k+4}\right)^{P(k)},\hspace{0.1in}
\left|{\bold H}^r_{\bold M}\right|_{k+1/2}
\leq C(k)\left(1 + |\rho|_{k+4}\right)^{P(k)},
\eqno(\arabic{equation})
\newcounter{AlmostEstimate2}
\setcounter{AlmostEstimate2}{\value{equation}}
\addtocounter{equation}{1}$$
where $|\rho|_s=\max_{1\leq \iota \leq N,1\leq l \leq m}\left\{|\rho_{\iota,l}|_s\right\}$.
\end{proposition}

\indent
Boundedness of different "solution operators"
${\bold R}_r :C^{k}_{(0,r)}({\bold M}) \rightarrow
C^{k+1/2 -\epsilon}_{(0,r)}({\bold M})$ for any $\epsilon >0$ was
first proved in \cite{AiH}.

\indent
Below we introduce notations and definitions necessary for the
construction of a formula.\\
\indent 
For a vector-valued $C^1$-function $\eta =(\eta_1,\dots,\eta_n)$ we will
use the notation:
$$\omega^{\prime}(\eta) = \sum_{k=1}^{n} (-1)^{k-1} \eta_k
\wedge_{j \neq k} d\eta_j, \hspace{0.2in} \omega(\eta) =
\wedge_{j = 1}^n d\eta_j.$$
\indent
If $\eta = \eta(\zeta, z,t)$ is a $C^1$-function of $\zeta \in \C^n,
z \in \C^n$ and a real parameter $t \in \R^l$ satisfying
the condition 
$$\sum_{k=1}^{n} \eta_k(\zeta, z,t) \cdot (\zeta_k - z_k) = 1
\eqno(\arabic{equation})
\newcounter{leray}
\setcounter{leray}{\value{equation}}
\addtocounter{equation}{1}$$
then
$$d \omega^{\prime}(\eta) \wedge \omega(\zeta) \wedge \omega(z) = 0$$\\
or, separating differentials,
$$d_{t} \omega^{\prime}(\eta) + \bar\partial_{\zeta} 
\omega^{\prime}(\eta) + \bar\partial_z \omega^{\prime}(\eta) = 0.
\eqno(\arabic{equation})
\newcounter{dofomega}
\setcounter{dofomega}{\value{equation}}
\addtocounter{equation}{1}$$
\indent
Also, if $\eta(\zeta,z,t)$ satisfies (\arabic{leray}) then the
differential form $\omega^{\prime}(\eta) \wedge \omega(\zeta) 
\wedge \omega(z)$ can be represented as:
$$\sum_{r=0}^{n-1}\omega^{\prime}_r(\eta)\wedge\omega(\zeta)
\wedge\omega(z),
\eqno(\arabic{equation})
\newcounter{sumomega}
\setcounter{sumomega}{\value{equation}}
\addtocounter{equation}{1}$$
where $\omega^{\prime}_r(\eta)$ is a differential form of the order $r$
in $d{\bar z}$ and respectively of the order $n-r-1$ in $d{\bar \zeta}$
and $dt$. From (\arabic{dofomega}) and (\arabic{sumomega}) follow
equalities:
$$d_t \omega^{\prime}_{r}(\eta)
+\bar\partial_{\zeta}\omega^{\prime}_{r}(\eta)
+ \bar\partial_z \omega^{\prime}_{r-1}(\eta) = 0
\hspace{0.3in} (r=1,\dots,n),
\eqno(\arabic{equation})
\newcounter{domega}
\setcounter{domega}{\value{equation}}
\addtocounter{equation}{1}$$
and
$$\omega^{\prime}_r(\eta) = \frac{1}{(n-r-1)!r!}
\mbox{Det} \Big[\eta, \hspace{0.05in}
\overbrace{\bar\partial_z \eta}^{r}, \hspace{0.05in}
\overbrace{\bar\partial_{\zeta,t} \eta}^{n-r-1}\Big],
\eqno(\arabic{equation})
\newcounter{determinantomega}
\setcounter{determinantomega}{\value{equation}}
\addtocounter{equation}{1}$$
where the determinant is calculated by the usual rules but with external
products of elements and the position of the element in the external product
is defined by the number of its column.\\
\indent
Let ${\cal U}$ be an open neighborhood in ${\bold G}$ and
$U = {\cal U} \cap {\bold M}$. We call a vector function
$$P(\zeta,z) = (P_1(\zeta,z),\dots, P_n(\zeta,z)) \hspace{0.2in}
\mbox{for} \hspace{0.2in} (\zeta,z) \in
\left( {\cal U} \setminus U \right)\times U$$
by strong ${\bold M}$-barrier for ${\cal U}$ if there exists
$C >0$ such that the inequality:
$$\left| \Phi(\zeta,z) \right|
> C \cdot \left( \rho(\zeta) + | \zeta - z|^2 \right)
\eqno(\arabic{equation})
\newcounter{barrier}
\setcounter{barrier}{\value{equation}}
\addtocounter{equation}{1}$$
holds for $(\zeta,z) \in \left({\cal U} \setminus U \right)\times U$,
where
$$\Phi(\zeta,z) = \langle P(\zeta,z), \zeta - z \rangle = 
\sum_{i=1}^n P_i(\zeta,z) \cdot (\zeta_i - z_i).$$
\indent
According to (\arabic{manifold}) we may assume that
$U = {\cal U} \cap {\bold M}$ is a set of common zeros of smooth functions
$\{ \rho_k, \hspace{0.1in}k=1, \dots, m \}$.
The Levi form of the function $\rho^2(z)=\sum_{j=1}^m \rho_j^2(z)$ is positive
definite on the complex subspaces ${\bold N}_z$ for any $z\in U$.
Therefore, scaling functions $\{ \rho_k, \hspace{0.1in}k=1, \dots, m \}$
if necessary, and using q-pseudoconcavity of ${\bold M}$
for any $z \in {\bold M}$ we can find an open neighborhood
$U \ni z$, and a family $E_{q+m}(\theta,z)$ of
$q+m$ dimensional complex linear subspaces in $\C^n$
smoothly depending on $(\theta, z) \in \S^{m-1} \times U$ and such that
$$-{\cal L}_{z}\rho_{\theta}-{\cal L}_{z}\rho^2$$
is strictly negative on $E_{q+m}(\theta,z)$ with all negative eigenvalues not exceeding
some $c < 0$.\\
\indent
For a set of functions $\rho_1, \dots , \rho_m \in C^p({\cal U})$ we
consider a family $E_{n-q-m}^{\bot}(\theta,z)$ of (n-q-m)-dimensional
subspaces in $T({\bold G})$, orthogonal to $E_{q+m}(\theta,z)$ and
a set of $C^p$-smooth vector functions
$$a_j(\theta,z) = \left( a_{j1}(\theta,z), \dots, a_{jn}(\theta,z)
\right) \hspace{0.05in} \mbox{for} \hspace{0.05in}j = 1, \dots, n-q-m$$
representing an orthonormal basis in $E_{n-q-m}^{\bot}(\theta,z).$\\
\indent
Defining for $(\theta, z, w) \in \S^{m-1}\times{\cal U}\times \C^n$
$$A_j(\theta,z,w) = \sum_{i=1}^n a_{ji}(\theta,z) \cdot w_i,
\hspace{0.1in}(j=1, \dots, n-q-m)$$
we construct the form
$${\cal A}(\theta, z, w) = \sum_{j=1}^{n-q-m} A_j(\theta,z,w) \cdot
{\bar A_j(\theta,z,w)}$$
such that the hermitian form
$${\dis\frac{1}{2}\left({\cal L}_{z}\rho_{\theta}(w)+{\cal L}_{z}\rho^2(w)\right)
+{\cal A}(\theta, z, w)}$$
is strictly positive definite in $w$ for $(\theta,z) \in \S^{m-1}
\times {\cal U}$.\\
\indent
Then we define for $\zeta, z \in\left({\cal U}\setminus U\right) \times U$:
$$\begin{array}{lllllll}
{\dis\theta_k(\zeta) = - \frac{\rho_k(\zeta)}{\rho(\zeta)} 
\hspace{0.2in}\mbox{for $k=1,\dots, m$},}\vspace{0.2in}\\
{\dis Q^{(k)}_i(\zeta, z)
= - \frac{\partial \rho_k}{\partial \zeta_i}(z)
-\frac{1}{2}\sum_{j=1}^{n} \frac{{\partial}^2 \rho_k}{\partial
\zeta_i \partial \zeta_j}(z)(\zeta_j - z_j)
-\frac{1}{2}\theta_k(\zeta)\sum_{j=1}^{n}\frac{{\partial}^2\rho^2}{\partial
\zeta_i \partial \zeta_j}(z)(\zeta_j - z_j),}\vspace{0.2in}\\
{\dis F^{(k)} (\zeta, z) = \langle Q^{(k)}(\zeta,z),
\zeta - z \rangle,}\vspace{0.2in}\\
{\dis P_i(\zeta, z) = \sum_{k=1}^m \theta_k(\zeta) \cdot Q^{(k)}_i(\zeta,z)
+ \sum_{j=1}^{n-q-m} a_{ji}(\theta(\zeta), z) \cdot
{\bar A_j}(\theta(\zeta),z,\zeta-z),}\vspace{0.2in}\\
{\dis\Phi (\zeta, z) = \langle P(\zeta, z), \zeta - z \rangle
= \sum_{k=1}^m \theta_k(\zeta) \cdot F^{(k)} (\zeta, z)
+ {\cal A}(\theta(\zeta),z,\zeta-z).}
\end{array}
\eqno(\arabic{equation})
\newcounter{Barrier}
\setcounter{Barrier}{\value{equation}}
\addtocounter{equation}{1}$$
\indent
To prove that $P_i(\zeta, z)$ is a strong
${\bold M}$-barrier for some $U \ni z$
we consider the Taylor expansions of $\rho_k$ for $k=1,\dots,m$,
and of $\rho^2$, and using that ${\dis \frac{\partial\rho^2_k}{\partial\zeta_i}(z)=0}$
for $z\in{\bf M}$ obtain
$$2\mbox{Re}F^{(k)}(\zeta, z)=\rho_k(z)-\rho_k(\zeta)
+{\cal L}_{z}\rho_k(\zeta - z) + O(|\zeta - z |^3)$$
$$+\theta_k(\zeta)\left[\rho^2(z)-\rho^2(\zeta)
+{\cal L}_{z}\rho^2(\zeta - z) + O(|\zeta - z |^3)\right].$$
\indent
Then we obtain for some $U$ and $(\zeta,z) \in
\left(U \setminus (U \cap {\bold M}) \right)
\times \left( U \cap {\bold M} \right)$:
$$\mbox{Re}\Phi(\zeta, z) =
\sum_{k=1}^m \theta_k(\zeta) \cdot \mbox{Re}F^{(k)}(\zeta, z)
+ \sum_{j=1}^{n-q-m} A_j(\theta(\zeta),\zeta,z)
\cdot {\bar A}_j(\theta(\zeta),\zeta,z)
\eqno(\arabic{equation})
\newcounter{RePhiestimate}
\setcounter{RePhiestimate}{\value{equation}}
\addtocounter{equation}{1}$$
$$=\frac{1}{2}\left(\rho(\zeta)-\rho^2(\zeta)+{\cal L}_{z}\rho_{\theta}(\zeta - z)
+{\cal L}_{z}\rho^2(\zeta - z)\right)+{\cal A}(\theta(\zeta),z,\zeta-z)+O(|\zeta - z |^3),$$
which implies the existence of an open
neighborhood $U \ni z$ in $\C^n,$ satisfying (\arabic{barrier}).\\
\indent
For $A_j(\zeta,z):= A_j(\theta(\zeta),z,\zeta-z)$ and $Q^{(k)}_i(\zeta, z)$
we have the following equalities that will be used in the estimates below
$$\begin{array}{ll}
\bar\partial_{\zeta}{\bar A_j}(\zeta,z)=
\mu_{\tau}^{(j)}(\zeta,z) + \mu_{\nu}^{(j)}(\zeta,z),\vspace{0.2in}\\
\bar\partial_{\zeta}Q^{(k)}_i(\zeta, z)=\chi^{(k)}_i(\zeta, z),
\end{array}
\eqno(\arabic{equation})
\newcounter{muandQ}
\setcounter{muandQ}{\value{equation}}
\addtocounter{equation}{1}$$
where
$$\mu_{\tau}^{(j)}(\zeta,z)
= \sum_{i=1}^n {\bar a_{ji}}(\theta(\zeta),z) d{\bar \zeta_i},
\hspace{0.2in}\mu_{\nu}^{(j)}(\zeta,z)= \sum_{i=1}^n
(\bar\zeta_i - \bar z_i) \bar\partial_{\zeta} {\bar a_{ji}}
(\theta(\zeta),z),$$
and
$$\chi^{(k)}_i(\zeta, z)=-\frac{1}{2}\left(\sum_{j=1}^{n}
\frac{{\partial}^2\rho^2}{\partial
\zeta_i\partial\zeta_j}(z)(\zeta_j - z_j)\right)
\bar\partial_{\zeta}\theta_k(\zeta).$$
\indent
In our description of local integral formulas on ${\bold M}$ and in the
future estimates we will also need the following notations.\\
\indent
We define the tubular neighborhood ${\cal U}(\epsilon)$ of ${\bold M}$ in
${\cal U}$ as follows:
$${\cal U}(\epsilon)= \{z \in  {\cal U}: \rho(z) < \epsilon \},$$
where $\rho(z) = {\left(\sum_{k=1}^m \rho_k^2(z)\right)}^{\frac{1}{2}}$.
The boundary of ${\cal U}(\epsilon)$ - $U(\epsilon)$ is defined
by the condition
$$U(\epsilon) = \{ z \in {\cal U}: \rho(z) = \epsilon \}.$$
We consider the fibration of ${\cal U}(\epsilon)$ by the manifolds
$${\bold M}\left(\delta_1,\dots,\delta_m\right)
= \left\{\zeta \in {\cal U}(\epsilon): \rho_1(\zeta)=\delta_1,\dots,
\rho_m(\zeta)=\delta_m\right\}$$
and denote by ${\cal T} \in {\bold C}T\left({\bold G}\right)$ the
subbundle of vectors tangent to fibers
${\bold M}\left(\delta_1,\dots,\delta_m\right)$.\\
\indent
For a sufficiently small neighborhood
${\cal U} \in {\bold G}$ we may assume that functions
$$\rho_k(\zeta),\hspace{0.05in}\mbox{Im}F^{(k)}(\zeta,z)
\hspace{0.05in}\{ k=1,\dots, m \}$$
have a nonzero jacobian with respect to $\mbox{Re}\zeta_{i_1}, \dots,
\mbox{Re}\zeta_{i_m},\hspace{0.05in}\mbox{Im}\zeta_{i_1}, \dots,
\mbox{Im}\zeta_{i_m}$ for $z, \zeta \in {\cal U}$.
Therefore, for any fixed $z \in {\cal U}$ these functions
may be chosen as local $C^p$ coordinates in $\zeta$.
We may also complement the functions above by holomorphic functions
$w_j(\zeta)=u_j(\zeta) + iv_j(\zeta)$ with
$j=1,\dots, n-m$ so that the functions
$$\rho_k(\zeta),\hspace{0.05in}\mbox{Im}F^{(k)}(\zeta,z)
\hspace{0.05in} \{ k=1,\dots, m \},$$
$$u_j(\zeta), v_j(\zeta)\hspace{0.05in} \{ j=1,\dots, n-m \},$$
represent a complete system of local coordinates in
$\zeta \in {\cal U}$ for any fixed $z \in {\cal U}$.\\
\indent
We consider complex valued vector fields on
${\cal U}$ for any fixed $z \in {\cal U}$:
$$Y_{i,\zeta}(z) = \frac{\partial}{\partial\mbox{Im}F^{(i)}(\zeta,z)}
\hspace{0.05in}\mbox{for}\hspace{0.05in}i=1,\dots, m,$$
$$W_{i,\zeta} = \frac{\partial}{\partial w_i},\hspace{0.05in}
{\overline W}_{i,\zeta} = \frac{\partial}{\partial {\bar w}_i}
\hspace{0.05in}\mbox{for}\hspace{0.05in} i=1,\dots, n-m,$$
and denote
$$Y_{\zeta}(z) = \sum_{k=1}^m \theta_k(\zeta)Y_{k,\zeta}(z).$$
\indent
We also introduce a local extension operator of functions and forms
from $U = {\cal U} \cap {\bold M}$ to ${\cal U}$
$$E:g \to E(g).$$
Assuming that locally manifold $\bold M$ in $\cal U$ with coordinates
$z_j=x_j+iy_j, j=1,\dots, n$ is defined as
$${\cal U} \cap {\bold M}=\left\{z\in {\cal U}:
\rho_j(z)\equiv x_j-h_j(y_1,\dots, y_m,z_{m+1},\dots, z_n)=0,
j=1,\dots, m\right\},
\eqno(\arabic{equation})
\newcounter{coordinates}
\setcounter{coordinates}{\value{equation}}
\addtocounter{equation}{1}$$
we define for a function $g(y_1,\dots y_m,z_{m+1},\dots, z_n)$ on $U$
$$E(g)(z_1,\dots,z_n)=g(y_1,\dots, y_m,z_{m+1},\dots, z_n),$$
extending a function identically with respect to
$x_1,\dots,x_m$. For a differential form
$$g = \sum_{I,J,K} g_{I,J,K}dy_I\wedge dz_J\wedge d{\overline z}_K$$
with multiindices $I \in (1,\dots,m),J,K \in (m+1,\dots,n)$
we define extension operator by extending coefficients as in the
formula above.\\
\indent
In our proof of Proposition~\ref{AlmostHomotopy} we will also use a special norm
for functions and forms on $\left({\cal U}\setminus U\right)_{\zeta}\times U_z$.
Namely, using functions $\rho,\theta_1,\dots,\theta_m$ and coordinates
$y_1,\dots, y_m,z_{m+1},\dots,z_n$ from (\arabic{coordinates})
we define a $C^p$-diffeomorphism
$$\Psi:{\cal U}(\epsilon)\setminus U\to\left(0,\epsilon\right)
\times \S^{m-1}\times U$$
by the formula
$$\Psi(\zeta)=\left(\rho(\zeta),\theta_1(\zeta),\dots,\theta_m(\zeta),
\eta_1,\dots,\eta_m,\zeta_{m+1},\dots,\zeta_n\right)$$
for $\zeta=\xi+i\eta$. Then for a form $h(\zeta,z)$ on
$\left({\cal U}(\epsilon)\setminus U\right)_{\zeta}\times U_z$
and $l\leq p$ we denote
$$\left\|g\right\|_l
=\left|g\left(\Psi^{-1}(\rho,\theta,\eta,\zeta^{\prime\prime}),
z\right)\right|_{C^l\left(\left[0,\epsilon\right]\times\S^{m-1}
\times U_{\zeta}\times U_z\right)},$$
with $\zeta^{\prime\prime}(\zeta)=(\zeta_{m+1},\dots,\zeta_n)$.\\
\indent
In what follows we will assume that the defining functions
$\rho_1,\dots,\rho_m$ satisfy condition
$$\min_{\cal U}\left|\mbox{Det}\left(\frac{\partial\rho_k}
{\partial\zeta_J}\right)^{k=1,\dots,m}_{|J|=m, J\in 1,\dots,n}\right| > c
\eqno(\arabic{equation})
\newcounter{NonDegeneracy}
\setcounter{NonDegeneracy}{\value{equation}}
\addtocounter{equation}{1}$$
for some fixed $c > 0$. Then from the construction of $A_j(\zeta,z)$
and $\Phi(\zeta,z)$ we conclude that for $l \leq p-2$ the estimates
$$\begin{array}{lll}
\left\|A_j\right\|_l\leq C,\vspace{0.1in}\\
\left\|Q^{(k)}_i\right\|_l,\ \left\|F^{(k)}\right\|_l,\
\left\|\Phi\right\|_l
\leq C \cdot \left( 1 + |\rho|_{l+2} \right),\vspace{0.1in}\\
\left\|Y_{\zeta}(z)\right\|_l \leq C
\left(1+|\rho|_{l+2}\right)^m.
\end{array}
\eqno(\arabic{equation})
\newcounter{lnorms}
\setcounter{lnorms}{\value{equation}}
\addtocounter{equation}{1}$$
hold for some $ C > 0$.\\
\indent
The following proposition provides local integral formula for
$\bar\partial_{\bold M}$.

\begin{proposition}\label{LocalHomotopy}
Let ${\bold M} \subset {\bold G}$ be a generic, regular q-pseudoconcave
CR submanifold of the class $C^p$ and let
${\cal U}$ be an open neighborhood in ${\bold G}$ with analytic
coordinates $z_1, \dots, z_n$.\\
\indent
Then for $r = 1,...,q-1$, $k \leq p$, and any differential form
$g \in C^k_{(0,r)}({\bold M})$ with compact support in
$U = {\cal U}\cap{\bold M}$ the following equality
$$g = \bar\partial_{\bold M} R_r(g) +
R_{r+1}(\bar\partial_{\bold M} g) + H_r(g),
\eqno(\arabic{equation})
\newcounter{LocalFormula}
\setcounter{LocalFormula}{\value{equation}}
\addtocounter{equation}{1}$$
holds, where
$$R_r(g)(z)$$
$$= (-1)^r \frac{(n-1)!}{(2\pi i)^n} \cdot \mbox{pr}_{\bold M} \circ
\lim_{\epsilon \rightarrow 0}
\int_{U(\epsilon)\times [0,1]} E(g)(\zeta)
\wedge\omega^{\prime}_{r-1}\left((1-t)\frac{\bar\zeta - \bar z}
{{\mid \zeta - z \mid}^2} + t\frac{P(\zeta,z)}
{\Phi(\zeta,z)}\right) \wedge\omega(\zeta),$$
$$H_r(g)(z) = (-1)^{r} \frac{(n-1)!}{(2\pi i)^n}
\cdot \mbox{pr}_{\bold M} \circ \lim_{\epsilon \rightarrow 0}
\int_{U(\epsilon)}E(g)(\zeta)
\wedge \omega^{\prime}_{r} \left( \frac{P(\zeta,z)}
{\Phi(\zeta,z)}\right) \wedge\omega(\zeta),$$
$E(g)$ is an extension of $g$ to ${\cal U}$,
$\Phi(\zeta,z)$ is a local barrier for ${\cal U}$ constructed
in (\arabic{Barrier}) and $\mbox{pr}_{\bold M}$ denotes the operator
of projection to the space of tangential differential forms on ${\bold M}$.
\end{proposition}
\vspace{0.2in}

\indent
We omit the proof of Proposition~\ref{LocalHomotopy} because it is
completely analogous to the proof of formula (\arabic{LocalFormula})
for another barrier function in \cite{P1}.\\

\indent
To construct now global formula on ${\bold M}$ we consider two finite
coverings
$\{ {\cal U}_{\iota} \subset
{\cal U}^{\prime}_{\iota}\}$ of
${\bold G}$ and two partitions of unity $\{\vartheta_{\iota} \}$ and
$\{\vartheta^{\prime}_{\iota} \}$ subordinate to these coverings and
such that $\vartheta^{\prime}_{\iota}(z) = 1$ for
$z \in \mbox{supp}(\vartheta_{\iota})$.\\
\indent
Applying Corollary~\ref{LocalHomotopy} to the form $\vartheta_{\iota}g$ in
${\cal U}^{\prime}_{\iota}$ we obtain
$$\vartheta_{\iota}(z)g(z) =
\bar\partial_{\bold M} R_r^{\iota}(\vartheta_{\iota}g)(z)
+ R_{r+1}^{\iota}(\bar\partial_{\bold M} \vartheta_{\iota}g)(z)
+ H_r^{\iota}(\vartheta_{\iota}g)(z).$$
\indent
Multiplying the equality above by $\vartheta^{\prime}_{\iota}(z)$ and
using equalities
$$\vartheta^{\prime}_{\iota}(z) \cdot \bar\partial_{\bold M}
R_r^{\iota}(\vartheta_{\iota}g)(z) = \bar\partial_{\bold M} \left[
\vartheta^{\prime}_{\iota}(z) \cdot
R_r^{\iota}(\vartheta_{\iota}g)(z) \right] - \bar\partial_{\bold M}
\vartheta^{\prime}_{\iota}(z) \wedge R_r^{\iota}(\vartheta_{\iota}g)(z)$$
and
$$R_{r+1}^{\iota}(\bar\partial_{\bold M} \vartheta_{\iota}g)(z) =
R_{r+1}^{\iota}(\bar\partial_{\bold M} \vartheta_{\iota} \wedge g)(z)
+ R_{r+1}^{\iota}(\vartheta_{\iota} \bar\partial_{\bold M} g)(z)$$
we obtain
$$\vartheta_{\iota}(z)g(z) = \bar\partial_{\bold M}
{\bold R}_r^{\iota}(g)(z)+{\bold R}_{r+1}^{\iota}(\bar\partial_{\bold M}g)(z)
+{\bold H}_r^{\iota}(g)(z)
\eqno(\arabic{equation})
\newcounter{varthetag}
\setcounter{varthetag}{\value{equation}}
\addtocounter{equation}{1}$$
with
$${\bold R}_r^{\iota}(g)(z) = \vartheta^{\prime}_{\iota}(z) \cdot
R_r^{\iota}(\vartheta_{\iota}g)(z)$$
and
$${\bold H}_r^{\iota}(g)(z) = - \bar\partial_{\bold M}
\vartheta^{\prime}_{\iota}(z) \wedge R_r^{\iota}(\vartheta_{\iota}g)(z)
+ \vartheta^{\prime}_{\iota}(z) \cdot
R_{r+1}^{\iota}(\bar\partial_{\bold M} \vartheta_{\iota} \wedge g)(z)
+ \vartheta^{\prime}_{\iota}(z) \cdot H_r^{\iota}(\vartheta_{\iota}g)(z).$$
\indent
Adding equalities (\arabic{varthetag}) for all $\iota$ we obtain

\begin{proposition}\label{GlobalHomotopy}
Let ${\bold M} \subset {\bold G}$ be a generic, regular q-pseudoconcave
compact CR submanifold of the class $C^p$.\\
\indent
Then for $r = 1,...,q-1$, $k \leq p$, and any differential form
$g \in C^k_{(0,r)}({\bold M})$ the following equality
$$g = \bar\partial_{\bold M}{\bold R}^r_{\bold M}(g) +
{\bold R}^{r+1}_{\bold M}(\bar\partial_{\bold M} g)
+ {\bold H}^r_{\bold M}(g),
\eqno(\arabic{equation})
\newcounter{GlobalFormula}
\setcounter{GlobalFormula}{\value{equation}}
\addtocounter{equation}{1}$$
holds, where
$${\bold R}^r_{\bold M}(g)(z) = \sum_{\iota} \vartheta^{\prime}_{\iota}(z)
\cdot R_r^{\iota}(\vartheta_{\iota}g)(z)
\eqno(\arabic{equation})
\newcounter{boldR}
\setcounter{boldR}{\value{equation}}
\addtocounter{equation}{1}$$
and
$${\bold H}^r_{\bold M}(g)(z) = \sum_{\iota} \left[ - \bar\partial_{\bold M}
\vartheta^{\prime}_{\iota}(z) \wedge R_r^{\iota}(\vartheta_{\iota}g)(z)
+ \vartheta^{\prime}_{\iota}(z)
R_{r+1}^{\iota}(\bar\partial_{\bold M} \vartheta_{\iota} \wedge g)(z)\right.
\eqno(\arabic{equation})
\newcounter{boldH}
\setcounter{boldH}{\value{equation}}
\addtocounter{equation}{1}$$
$$\left.+ \vartheta^{\prime}_{\iota}(z)
\cdot H_r^{\iota}(\vartheta_{\iota}g)(z)\right].$$
\end{proposition}
\vspace{0.2in}

\section{Estimates for ${\bold R}^r_{\bold M}$
and ${\bold H}^r_{\bold M}$.}\label{RHEstimates}

\indent
From the construction of operator ${\bold R}^r_{\bold M}$ we conclude that
in order to prove necessary estimates for operator ${\bold R}^r_{\bold M}$
it suffices to prove these estimates for operator $R_r$. In the proposition
below we state necessary estimates for operator $R_r$.\\

\begin{proposition}\label{REstimates}
Let ${\bold M} \subset {\bold G}$ be a generic, regular q-pseudoconcave
CR submanifold of the class $C^p$ in ${\cal U}$ satisfying condition
(\arabic{NonDegeneracy}) and let
$g \in C^k_{(0,r)}({\bold M})$ $(1 \leq k \leq p-4)$ be a form with
compact support in $U = {\cal U} \cap {\bold M}$.\\
\indent
Then $R_r(g)$ defined in (\arabic{LocalFormula}) satisfies the
following estimates
$$\left| R_r(g) \right|_k
\leq C(k)\left(1 + |\rho|_{k+3}\right)^{P(k)}\left| g\right|_k,
\eqno(\arabic{equation})
\newcounter{REstimate1}
\setcounter{REstimate1}{\value{equation}}
\addtocounter{equation}{1}$$
$$\left| R_r(g) \right|_{k+1/2}
\leq C(k)\left(1 + |\rho|_{k+4}\right)^{P(k)}\left| g\right|_k,
\eqno(\arabic{equation})
\newcounter{REstimate2}
\setcounter{REstimate2}{\value{equation}}
\addtocounter{equation}{1}$$
with $P(k)$ a polynomial in $k$ and a constant $C(k)$ independent of $g$.
\end{proposition}

\indent
{\bf Proof.}\hspace{0.05in}
In our proof of Proposition~\ref{REstimates} we will
use the approximation of $R_r$ by the operators
$$R_r(\epsilon)(f)(z) = (-1)^{r} \frac{(n-1)!}{(2\pi i)^n}
\eqno(\arabic{equation})
\newcounter{OperatorREpsilon}
\setcounter{OperatorREpsilon}{\value{equation}}
\addtocounter{equation}{1}$$
$$\times\int_{U(\epsilon)\times[0,1]}
\vartheta(\zeta)E(f)(\zeta)\omega^{\prime}_{r-1}
\left((1-t)\frac{\bar\zeta - \bar z}{{\mid \zeta - z \mid}^2}
+ t\frac{P(\zeta,z)}{\Phi(\zeta,z)}\right) \wedge\omega(\zeta)$$
when $\epsilon$ goes to $0$.\\
\indent
Using equalities (\arabic{muandQ}) we obtain the following representation
of kernels of these integrals on ${\cal U}\times[0,1]\times{\bold M}$:
$$\left.\vartheta(\zeta)\cdot \omega^{\prime}_{r-1}
\left((1-t)\frac{\bar\zeta - \bar z}{{\mid \zeta - z \mid}^2}
+ t\frac{P(\zeta,z)}{\Phi(\zeta,z)}\right) \wedge\omega(\zeta)
\right|_{{\cal U} \times [0,1] \times {\bold M} }
\eqno(\arabic{equation})
\newcounter{cauchymartinellikernel}
\setcounter{cauchymartinellikernel}{\value{equation}}
\addtocounter{equation}{1}$$
$$= \sum_{i,J} a_{(i,J)}(t,\zeta,z)dt \wedge
\lambda^{i,J}_{r-1}(\zeta, z) +
\sum_{i,J} b_{(i,J)}(t,\zeta,z) dt \wedge
\gamma^{i,J}_{r-1}(\zeta, z),$$
where $i$ is an index, $J= \cup_{i=1}^{8} J_i$ is a
multiindex such that $i \not \in J,$
$a_{(i,J)}(t,\zeta,z)$ and $b_{(i,J)}(t,\zeta,z)$ are polynomials in $t$
with coefficients that are $C^p$-functions of  $z$, $\zeta$ and
$\theta(\zeta)$, and $\lambda^{i,J}_{r-1}(\zeta, z)$
and $\gamma^{i,J}_{r-1}(\zeta, z)$ are defined as follows:
$$\lambda^{i,J}_{r-1}(\zeta, z) =\frac{1}
{|\zeta - z|^{2(|J_1|+|J_5|+1)} \cdot
{\Phi(\zeta,z)}^{n-|J_1|-|J_5|-1}}
\eqno(\arabic{equation})
\newcounter{lambdaForm}
\setcounter{lambdaForm}{\value{equation}}
\addtocounter{equation}{1}$$
$$\times \sum \mbox{Det} \left[\bar{\zeta} - \bar{z},\hspace{0.03in}
Q^{(i)},\hspace{0.03in}
\overbrace{d\bar{\zeta}}^{j \in J_1},\hspace{0.03in}
\overbrace{{\bar A} \cdot \bar\partial_{\zeta}a}^{j \in J_2},\hspace{0.03in}
\overbrace{a \cdot \mu_{\nu},\hspace{0.03in}\chi}^{j \in J_3},\hspace{0.03in}
\overbrace{a \cdot \mu_{\tau}}^{j \in J_4}, \hspace{0.03in}
\overbrace{d\bar{z}}^{j \in J_5},\hspace{0.03in}
\overbrace{{\bar A} \cdot \bar \partial_z a}^{j \in J_6},\hspace{0.03in}
\overbrace{a \cdot \bar\partial_z {\bar A}}^{j \in J_7},\hspace{0.03in}
\overbrace{\bar\partial_z Q}^{j \in J_8}
\right] \wedge\omega(\zeta),$$
and
$$\gamma^{i,J}_{r-1}(\zeta, z) =\frac{1}
{|\zeta - z|^{2(|J_1|+|J_5|+1)} \cdot
{\Phi(\zeta,z)}^{n-|J_1|-|J_5|-1}}
\eqno(\arabic{equation})
\newcounter{gammaForm}
\setcounter{gammaForm}{\value{equation}}
\addtocounter{equation}{1}$$
$$\times \sum \mbox{Det} \left[\bar{\zeta} - \bar{z},\hspace{0.03in}
a_i {\bar A_i},\hspace{0.03in}
\overbrace{d\bar{\zeta}}^{j \in J_1},\hspace{0.03in}
\overbrace{{\bar A} \cdot \bar\partial_{\zeta}a}^{j \in J_2},\hspace{0.03in}
\overbrace{a \cdot \mu_{\nu},\hspace{0.03in}\chi}^{j \in J_3},\hspace{0.03in}
\overbrace{a \cdot \mu_{\tau}}^{j \in J_4}, \hspace{0.03in}
\overbrace{d\bar{z}}^{j \in J_5},\hspace{0.03in}
\overbrace{{\bar A} \cdot \bar \partial_z a}^{j \in J_6},\hspace{0.03in}
\overbrace{a \cdot \bar\partial_z {\bar A}}^{j \in J_7},\hspace{0.03in}
\overbrace{\bar\partial_z Q}^{j \in J_8}
\right] \wedge\omega(\zeta).$$
\indent
We joined together terms of the form $a \cdot \mu_{\nu}$ and $\chi$ since according
to (\arabic{muandQ}) they are similar from the point of view of estimates.\\
\indent
In the proof of boundedness of operators
$R_r:C^{k}_{(0,r)}({\bold M}) \rightarrow C^{k+1/2}_{(0,r-1)}({\bold M})$
we will use smoothness of integrals with kernels
$\lambda^{i,J}_{r-1}(\zeta, z)$ and $\gamma^{i,J}_{r-1}(\zeta, z)$.
Smoothness of these integrals with respect to "CR tangent" and " CR normal"
vector fields was investigated in \cite{P2}. Below we describe some of
the constructions from there.\\
\indent
We consider kernels:
$${\cal K}^{I}_{d,h}(\zeta, z) =
\frac{ \{ \rho(\zeta) \}^{I_1}
(\zeta - z)^{I_2}(\bar{\zeta} - \bar{z})^{I_3} }
{|\zeta -z|^d \cdot \Phi(\zeta,z)^{h}}
\overbrace{\wedge d\rho_i}^{i \in I_4}
\overbrace{\wedge d\theta_i(\zeta)}^{i \in I_5}
\wedge d\sigma_{2n-m}(\zeta,z),$$
where $I= \cup_{j=1}^5 I_j$ and $I_j$ for $j=1, \dots, 5$ are
multiindices such that $I_1$ contains $m$ indices, $I_2$, $I_3$
contain $n$ indices, $I_4 \cup I_5$ contains $m-1$ indices,
$|I_4|+|I_5|=m-1$,
$\{ \rho(\zeta) \}^{I_1}= \prod_{i_s \in I_1} \rho_s(\zeta)^{i_s}$,
$(\zeta - z)^{I_2} = \prod_{i_s \in I_2} (\zeta_s - z_s)^{i_s}$,
$(\bar{\zeta} - \bar{z})^{I_3} = \prod_{i_s \in I_3}
(\bar{\zeta}_s - \bar{z}_s)^{i_s}$, and
$$d\sigma_{2n-m}(\zeta,z)=\bigwedge_{k=1}^m
d_{\zeta}\mbox{Im}F^{(k)}(\zeta,z)\bigwedge_{i=1}^{n-m}
\left(dw_i\wedge d{\bar w}_i\right).$$
\indent
For kernels ${\cal K}^{I}_{d,h}$ we use the following notation
$$\begin{array}{lll}
k\left({\cal K}^{I}_{d,h}\right)=d-|I_2|-|I_3|,\vspace{0.1in}\\
h\left({\cal K}^{I}_{d,h}\right)=h,\vspace{0.1in}\\
l\left({\cal K}^{I}_{d,h}\right)=|I_1|+|I_4|.
\end{array}$$
\indent
The lemma below is a refinement of Lemma 3.2 from \cite{P2} for
$C^k$ forms and $C^p$ vector fields.\\

\begin{lemma}\label{Smoothness}
Let ${\bold M} \subset {\bold G}$ be a compact, generic, regular
q-pseudoconcave $C^p$ submanifold,
${\cal U}$ and $U = {\cal U} \cap {\bold M}$ be neighborhoods such that
(\arabic{NonDegeneracy}) is satisfied for some fixed $c > 0$ in ${\cal U}$
and let $\Phi(\zeta, z)$ be as in (\arabic{Barrier}).
Let $g(\zeta,z,\theta,t)$ be a $C^l$ form $(1 \leq l \leq p-3)$
with compact support in
${\cal U}_{\zeta} \times {\cal U}_{z}
\times \S^{m-1} \times[0,1]$.\\
\indent
Then for $g(\zeta, z, t)=g(\zeta,z,\theta(\zeta),t)$ and a $C^p$
vector field
$$D_z = \sum_{j=1}^n a_j(z) \frac{\partial}{\partial z_j} +
\sum_{j=1}^n b_j(z) \frac{\partial}{\partial\bar z_j}
\in C^p\left({\cal U},{\cal T}\right)$$
such that $\left\| D_z \right\|_p \leq 1$ the following equality holds
$$D_z \left( \int_{U(\epsilon) \times [0,1]} g(\zeta, z, t) \cdot
{\cal K}^{I}_{d,h}(\zeta, z) dt \right)
\eqno(\arabic{equation})
\newcounter{DKRepresentation}
\setcounter{DKRepresentation}{\value{equation}}
\addtocounter{equation}{1}$$
$$= \int_{U(\epsilon) \times [0,1]}
\left[ D_z g(\zeta, z, t) \right]
\cdot {\cal K}^{I}_{d,h}(\zeta, z)dt
+ \int_{U(\epsilon) \times [0,1]}
\left[ D_{\zeta} g(\zeta, z, t) \right] \cdot
{\cal K}^{I}_{d,h}(\zeta, z) dt$$
$$+ \sum_{S,a,b} \int_{U(\epsilon) \times [0,1]}
c_{ \{ S,a,b \}}(\zeta, z, t) \cdot
g(\zeta, z, t) \cdot {\cal K}^{S}_{a,b}(\zeta, z) dt$$
$$+ \sum_{L,i} \int_{U(\epsilon) \times [0,1]}
c_{ \{L,i \}}(\zeta, z, t) \cdot
\left[ Y_{\zeta}(z) g(\zeta, z, t) \right]
\cdot {\cal K}^{L}_{i,h}(\zeta, z) dt,$$
where vector field $D_{\zeta}$ is defined as
$$D_{\zeta} = \sum_{j=1}^n a_j(\zeta) \frac{\partial}{\partial {\zeta}_j} +
\sum_{j=1}^n b_j(\zeta) \frac{\partial}{\partial\bar{\zeta}_j},$$
$$\left\|c_{ \{ S,a,b \}}\right\|_l,
\left\|c_{ \{ L,i,e,k \}}\right\|_l \leq C\left(1+|\rho|_{l+3}\right)^{3m},
\eqno(\arabic{equation})
\newcounter{cEstimates}
\setcounter{cEstimates}{\value{equation}}
\addtocounter{equation}{1}$$
with some $C > 0$ and kernels ${\cal K}^{S}_{a,b}$ and
${\cal K}^{L}_{i,h}$ satisfy the following conditions
$$\begin{array}{ll}
k\left({\cal K}^{S}_{a,b}\right)+b
-l\left({\cal K}^{S}_{a,b}\right)
\leq k\left({\cal K}^{I}_{d,h}\right)+h
-l\left({\cal K}^{I}_{d,h}\right),
\vspace{0.1in}\\
k\left({\cal K}^{S}_{a,b}\right)+2b
-2l\left({\cal K}^{S}_{a,b}\right)
\leq k\left({\cal K}^{I}_{d,h}\right)+2h
-2l\left({\cal K}^{I}_{d,h}\right).
\end{array}
\eqno(\arabic{equation})
\newcounter{SmoothnessIndices}
\setcounter{SmoothnessIndices}{\value{equation}}
\addtocounter{equation}{1}$$
\end{lemma}
\vspace{0.2in}

\indent
{\bf Proof.}\hspace{0.05in}
To prove the lemma we represent the integral from the left hand side of
(\arabic{DKRepresentation}) as
$$D_z \left( \int_{U(\epsilon) \times [0,1]} g(\zeta, z, t) \cdot
{\cal K}^{I}_{d,h}(\zeta, z) dt \right)
\eqno(\arabic{equation})
\newcounter{DzFormula}
\setcounter{DzFormula}{\value{equation}}
\addtocounter{equation}{1}$$
$$= \int_{U(\epsilon) \times [0,1]} \left[ D_z g(\zeta, z, t) \right]
\cdot {\cal K}^{I}_{d,h}(\zeta, z) dt
- \int_{U(\epsilon) \times [0,1]} g(\zeta, z, t) 
\left[ D_{\zeta} {\cal K}^{I}_{d,h}(\zeta, z) \right]dt$$
$$+ \int_{U(\epsilon) \times [0,1]} g(\zeta, z, t) 
\left[ \left( D_z + D_{\zeta} \right)
{\cal K}^{I}_{d,h}(\zeta, z) \right]dt.$$
\indent
To transform the second term of the right hand side of
(\arabic{DzFormula}) we apply integration by parts and obtain
$$\int_{U(\epsilon) \times [0,1]} g(\zeta, z, t) 
\left[ D_{\zeta} {\cal K}^{I}_{d,h}(\zeta, z) \right]dt
= - \int_{U(\epsilon) \times [0,1]}
\left[ D_{\zeta} g(\zeta, z, t) \right]
{\cal K}^{I}_{d,h}(\zeta, z) dt$$
$$+ \sum_{S,a,b} \int_{U(\epsilon) \times [0,1]}
c_{ \{ S,a,b \}}(\zeta, z, t) \cdot 
g(\zeta, z, t) \cdot {\cal K}^{S}_{a,b}(\zeta, z)dt,$$
with $\left\|c_{ \{ S,a,b \}}\right\|_l \leq 1$
and kernels ${\cal K}^{S}_{a,b}$ satisfying (\arabic{SmoothnessIndices}).\\
\indent
To transform the third term of the right hand side of
(\arabic{DzFormula}) we will use the formulas below that follow from
the definitions of $F^{(k)}(\zeta,z)$ and ${\cal A}(\zeta,z)$, estimates
(\arabic{lnorms}) and
from the fact that $D_z \in C^p\left({\cal U},{\cal T}\right)$:
$$\begin{array}{lllll}
\left( D_z + D_{\zeta} \right)
{\cal A}(\zeta, z) = \sum_{l_i+l_j=2} \alpha^{l_il_j}_{ij}(\zeta,z,t)
(\zeta_i - z_i)^{l_i}({\bar\zeta}_j - {\bar z}_j)^{l_j},
\vspace{0.2in}\\
\left( D_z + D_{\zeta} \right)
\mbox{Re} F^{(k)}(\zeta,z) = \sum_{l_i+l_j=2} \beta^{l_il_j}_{ij}(\zeta,z,t)
(\zeta_i - z_i)^{l_i}({\bar\zeta}_j - {\bar z}_j)^{l_j},
\vspace{0.2in}\\
\left( D_z + D_{\zeta} \right)
\mbox{Im} F^{(k)}(\zeta,z) = \sum_{i=1}^n \beta_{i}(\zeta,z,t)
(\zeta_i - z_i) + \sum_{i=1}^n {\bar\beta}_{i}(\zeta,z,t)
({\bar\zeta}_i - {\bar z}_i),
\vspace{0.2in}\\
\left( D_z + D_{\zeta} \right)(\zeta_j - z_j)
= \sum_{i=1}^n \alpha^1_i(\zeta,z,t)
(\zeta_i - z_i) + \sum_{i=1}^n {\bar\alpha}^1_{i}(\zeta,z,t)
({\bar\zeta}_i - {\bar z}_i),
\vspace{0.2in}\\
\left( D_z + D_{\zeta} \right)({\bar\zeta}_j - {\bar z}_j)
= \sum_{i=1}^n \alpha^2_i(\zeta,z,t)
(\zeta_i - z_i) + \sum_{i=1}^n {\bar\alpha}^2_i(\zeta,z,t)
({\bar\zeta}_i - {\bar z}_i),
\end{array}
\eqno(\arabic{equation})
\newcounter{DzDzetaFA}
\setcounter{DzDzetaFA}{\value{equation}}
\addtocounter{equation}{1}$$
with
$$\max\left\{\left\| \alpha^{k_ik_j}_{ij} \right\|_l,\
\left\| \alpha^k_i \right\|_l,\
\left\| {\bar\alpha}^k_i \right\|_l\right\}\leq C,\ \mbox{and}\
\max\left\{\left\| \beta^{k_ik_j}_{ij} \right\|_l,\
\left\| \beta_i \right\|_l,\
\left\| {\bar\beta}_i \right\|_l\right\}
\leq C \cdot \left(1+ |\rho|_{l+3}\right)$$
for some $C > 0$.\\
\indent
Applying operators $D_z$ and $D_{\zeta}$ to
${\cal K}^{I}_{d,h}(\zeta, z)$ and using formulas (\arabic{DzDzetaFA})
we obtain
$$\left( D_z + D_{\zeta} \right) {\cal K}^{I}_{d,h}(\zeta, z)
\eqno(\arabic{equation})
\newcounter{DzDzetaK2}
\setcounter{DzDzetaK2}{\value{equation}}
\addtocounter{equation}{1}$$
$$= (-h) \frac{ \{ \rho(\zeta) \}^{I_1}
(\zeta - z)^{I_2} (\bar{\zeta} - \bar{z})^{I_3} }
{|\zeta -z|^d \cdot \Phi(\zeta,z)^{h+1}}
\left( \left[ \left( D_z + D_{\zeta} \right) {\cal A} \right]
\overbrace{\wedge d\rho_i}^{i \in I_4}
\wedge \overbrace{d\theta_i(\zeta)}^{i \in I_5}
\wedge d\sigma_{2n-m}(\zeta,z) \right.$$
$$+ \sum_{k=1}^m \theta_k(\zeta) \cdot
\left[ \left( D_z + D_{\zeta} \right) \mbox{Re}F^{(k)} \right]
\overbrace{\wedge d\rho_i}^{i \in I_4}
\wedge \overbrace{d\theta_i(\zeta)}^{i \in I_5}
\wedge d\sigma_{2n-m}(\zeta,z)$$
$$\left. + \sum_{k=1}^m \theta_k(\zeta) \cdot
\left[ \left( D_z + D_{\zeta} \right) \mbox{Im}F^{(k)} \right]
\overbrace{\wedge d\rho_i}^{i \in I_4}
\wedge \overbrace{d\theta_i(\zeta)}^{i \in I_5}
\wedge d\sigma_{2n-m}(\zeta,z) \right)$$
$$+ \sum_{S,a,b} c_{ \{ S,a,b \}}(\zeta, z, t) \cdot 
{\cal K}^{S}_{a,b}(\zeta, z)$$
with $\left\|c_{ \{ S,a,b \}}\right\|_l \leq C$
and kernels ${\cal K}^{S}_{a,b}$ satisfying (\arabic{SmoothnessIndices}).\\
\indent
From estimates (\arabic{DzDzetaFA}) we conclude that the first two
terms of the right hand side of (\arabic{DzDzetaK2})
can be represented as linear combinations with coefficients satisfying
$\left\|c_{ \{ S,a,b \}}\right\|_l \leq C\left(1+|\rho|_{l+3}\right)$
of kernels ${\cal K}^{S}_{d,h+1}$ with
$$|S_2| + |S_3| = |I_2| + |I_3| + 2,$$
or
$$k\left({\cal K}^{S}_{d,h+1}\right) \leq
k\left({\cal K}^{I}_{d,h+1}\right)-2,$$
and, therefore, satisfying (\arabic{SmoothnessIndices}).\\
\indent
Considering equality
$$d\sigma_{2n-m}(\zeta,z)=c(\zeta,z)\bigwedge_{j=1}^m
d_{\zeta}\mbox{Im}F^{(j)}(\zeta,z)\bigwedge_{i=1}^{n-m}
\left(dw_i\wedge d{\bar w}_i\right)$$
with
$$\begin{array}{ll}
\left\|c\right\|_0 >C_1,\vspace{0.1in}\\
\left\|c\right\|_l \leq C_2\left(1+|\rho|_{l+2}\right)^m,
\end{array}
\eqno(\arabic{equation})
\newcounter{cEstimate}
\setcounter{cEstimate}{\value{equation}}
\addtocounter{equation}{1}$$
and defining
$$d_{\zeta}\mbox{Im}\Phi(\zeta,z)\interior d\sigma_{2n-m}(\zeta,z)
=c(\zeta,z)\left[\sum_{k=1}^m(-1)^{k-1}\theta_k(\zeta)\bigwedge_{j\neq k}
d_{\zeta}\mbox{Im}F^{(j)}(\zeta,z)\bigwedge_{i=1}^{n-m}
\left(dw_i\wedge d{\bar w}_i\right)\right],$$
we represent
the third term of the right hand side of (\arabic{DzDzetaK2}) as
$$(\pm 1)\sum_{k=1}^m\frac{\theta_k(\zeta)\cdot
\left[ \left( D_z + D_{\zeta} \right)
\mbox{Im}F^{(k)}(\zeta,z) \right] \{ \rho(\zeta) \}^{I_1} 
(\zeta - z)^{I_2} (\bar{\zeta} - \bar{z})^{I_3} }{|\zeta -z|^d }
\cdot\frac{1}{c(\zeta,z)}$$
$$\times d_{\zeta}\left[\frac{\left( d_{\zeta}\mbox{Im}\Phi(\zeta,z)\interior
d\sigma_{2n-m}(\zeta,z)\right)\wedge\overbrace{d\rho_i}^{i \in I_4}
\wedge \overbrace{d\theta_i(\zeta)}^{i \in I_5}}
{\Phi(\zeta,z)^{h}} \right]$$
$$+\sum_{S,a,b} c_{ \{ S,a,b \}}(\zeta, z, t) \cdot 
{\cal K}^{S}_{a,b}(\zeta, z)$$
with $\left\|c_{ \{ S,a,b \}}\right\|_l \leq C\left(1+|\rho|_{l+3}\right)^{2m}$
and kernels ${\cal K}^{S}_{a,b}$ satisfying 
(\arabic{SmoothnessIndices}).\\
\indent
Applying integration by parts to corresponding integrals,
using third formula from (\arabic{DzDzetaFA}), estimate
(\arabic{lnorms}) for $\left\|Y_{\zeta}(z)\right\|$, and estimate
(\arabic{cEstimate}) we obtain
$$\sum_{k=1}^m \int_{U(\epsilon) \times [0,1]}g(\zeta, z, t)
\cdot \frac{\theta_k(\zeta) \cdot \left[ \left( D_z + D_{\zeta} \right)
\mbox{Im}F^{(k)}(\zeta,z) \right] \{ \rho(\zeta) \}^{I_1}
(\zeta - z)^{I_2} (\bar{\zeta} - \bar{z})^{I_3} }{|\zeta -z|^d }
\cdot\frac{1}{c(\zeta,z)}
\eqno(\arabic{equation})
\newcounter{DzDzetaK3}
\setcounter{DzDzetaK3}{\value{equation}}
\addtocounter{equation}{1}$$
$$\times d_{\zeta} \left[\frac{\left( d_{\zeta}\mbox{Im}\Phi(\zeta,z)\interior
d\sigma_{2n-m}(\zeta,z)\right)\overbrace{\wedge d\rho_i}^{i \in I_4}
\wedge \overbrace{d\theta_i(\zeta)}^{i \in I_5}}{\Phi(\zeta,z)^{h}} \right]$$
$$=-\sum_{k=1}^m \int_{U(\epsilon) \times [0,1]}
\theta_k(\zeta) g(\zeta, z, t) \cdot Y_{\zeta}(z)
\left[\frac{\left[ \left( D_z + D_{\zeta} \right)
\mbox{Im}F^{(k)}(\zeta,z) \right] \{ \rho(\zeta) \}^{I_1}
(\zeta - z)^{I_2} (\bar{\zeta} - \bar{z})^{I_3} }{|\zeta -z|^d c(\zeta,z)}
\right]$$
$$\times \frac{\left( d_{\zeta}\mbox{Im}\Phi(\zeta,z)\interior
d\sigma_{2n-m}(\zeta,z)\right)\wedge\overbrace{d\rho_i}^{i \in I_4}
\wedge \overbrace{d\theta_i(\zeta)}^{i \in I_5}}{\Phi(\zeta,z)^{h}}$$
$$-\sum_{k=1}^m \int_{U(\epsilon) \times [0,1]}
\theta_k(\zeta) \left[ Y_{\zeta}(z) g(\zeta, z, t) \right] \cdot
\frac{\left[ \left( D_z + D_{\zeta} \right)
\mbox{Im}F^{(k)}(\zeta,z) \right] \{ \rho(\zeta) \}^{I_1} 
(\zeta - z)^{I_2} (\bar{\zeta} - \bar{z})^{I_3} }{|\zeta -z|^d c(\zeta,z)}$$
$$\times \frac{\left( d_{\zeta}\mbox{Im}\Phi(\zeta,z)\interior
d\sigma_{2n-m}(\zeta,z)\right)\overbrace{\wedge d\rho_i}^{i \in I_4}
\wedge \overbrace{d\theta_i(\zeta)}^{i \in I_5}}{\Phi(\zeta,z)^{h}}$$
$$= \sum_{S,a} \int_{U(\epsilon) \times [0,1]}
c_{ \{ S,a \}}(\zeta, z, t) \cdot 
g(\zeta, z, t) \cdot {\cal K}^{S}_{a,h}(\zeta, z) dt$$
$$+ \sum_{L,i} \int_{U(\epsilon) \times [0,1]}
c_{ \{L,i \}}(\zeta, z, t) \cdot
\left[ Y_{\zeta}(z) g(\zeta, z, t) \right]
\cdot {\cal K}^{L}_{i,h}(\zeta, z) dt$$
with $\left\|c_{ \{ S,a \}}\right\|_l< C\left(1+|\rho|_{l+3}\right)^{3m}$,
$\left\|c_{ \{L,i \}}\right\|_l < C\left(1+|\rho|_{l+3}\right)^m$
and kernels satisfying (\arabic{SmoothnessIndices}).\\
\indent
Combining now formulas (\arabic{DzFormula}), (\arabic{DzDzetaK2})
and (\arabic{DzDzetaK3}) we obtain statement of the Lemma.\qed

\indent
Using now Lemma~\ref{Smoothness} we reduce the statement of
Proposition~\ref{REstimates} to the case $k=0$. In this reduction
we will need the following simple lemma.
                                                                                                                                                                                                                                                                                                                                                                                                                                                                                                                           
\begin{lemma}\label{TangentDifferentials}
Let ${\bold M}$ be a generic CR submanifold of class $C^p$ in the unit ball
${\bold B}^n$ in $\C^n$ of the form:
$${\bold M} = \{\zeta \in {\bold B}^n: \rho_{1}(\zeta)=\dots
=\rho_{m}(\zeta)=0\},$$
where $\{ \rho_{k} \}, \  k = 1, \dots , m  \ ( m < n)$ are real
valued functions of the class $C^p$ satisfying conditions
(\arabic{NonDegeneracy}) on ${\bold B}^n$ for some $c > 0$.\\
\indent
Then for any point $\zeta_0 \in {\bold M}$ there exists a neighborhood
${\bold V}_{\epsilon}(\zeta_0) = \{ \zeta : | \zeta - \zeta_0 | <
\epsilon \}$ such that for any $l>n-m$ the
following representation holds in ${\bold V}_{\epsilon}:$
$$d{\bar\zeta}_{i_1}\wedge\dots\wedge d{\bar\zeta}_{i_l}\wedge
d\zeta_1\wedge\dots\wedge d\zeta_n = \sum d\rho_{j_1}
\wedge\dots\wedge d\rho_{j_{l-n+m}}\wedge
g_{j_1 \dots j_{l-n+m}}^{i_1 \dots i_l}(\zeta)
\eqno(\arabic{equation})
\newcounter{tangentestimate}
\setcounter{tangentestimate}{\value{equation}}
\addtocounter{equation}{1}$$
with $\left|g_{j_1 \dots j_{l-n+m}}^{i_1 \dots i_l}\right|_p \leq C$.
\end{lemma}
\qed\\
\indent
According to (\arabic{OperatorREpsilon}) and (\arabic{cauchymartinellikernel})
in order to prove the statement of the Proposition~\ref{REstimates} it suffices
to prove the estimates
$$\begin{array}{ll}
\left| \int_{U(\epsilon) \times [0,1]}
a_{(i,J)}(t,\zeta,z)dt \wedge E(g)(\zeta) \wedge
\lambda^{i,J}_{r-1}(\zeta, z) \right|_k
\leq C(k)\left(1 + |\rho|_{k+3}\right)^{P(k)}
\left| g \right|_{k},
\vspace{0.2in}\\
\left| \int_{U(\epsilon) \times [0,1]}
b_{(i,J)}(t,\zeta,z)dt \wedge E(g)(\zeta) \wedge
\gamma^{i,J}_{r-1}(\zeta, z) \right|_k
\leq C(k)\left(1 + |\rho|_{k+3}\right)^{P(k)}
\left| g \right|_{k}
\end{array}
\eqno(\arabic{equation})
\newcounter{lambdagammaEstimates1}
\setcounter{lambdagammaEstimates1}{\value{equation}}
\addtocounter{equation}{1}$$
and
$$\begin{array}{ll}
\left| \int_{U(\epsilon) \times [0,1]}
a_{(i,J)}(t,\zeta,z)dt \wedge E(g)(\zeta) \wedge
\lambda^{i,J}_{r-1}(\zeta, z) \right|_{k+1/2}
\leq C(k)\left(1 + |\rho|_{k+4}\right)^{P(k)}
\left| g \right|_{k},
\vspace{0.2in}\\
\left| \int_{U(\epsilon) \times [0,1]}
b_{(i,J)}(t,\zeta,z)dt \wedge E(g)(\zeta) \wedge
\gamma^{i,J}_{r-1}(\zeta, z) \right|_{k+1/2}
\leq C(k)\left(1 + |\rho|_{k+4}\right)^{P(k)}
\left| g \right|_{k}$$
\end{array}
\eqno(\arabic{equation})
\newcounter{lambdagammaEstimates2}
\setcounter{lambdagammaEstimates2}{\value{equation}}
\addtocounter{equation}{1}$$
with constants $C(k)$ and $P(k)$ independent of $g$ and $\epsilon$.\\
\indent
Using formulas (\arabic{muandQ}) and estimates (\arabic{lnorms}) for
the terms of determinants in (\arabic{lambdaForm}) and (\arabic{gammaForm})
and applying Lemma~\ref{TangentDifferentials} to the differential form
$$\overbrace{d{\bar\zeta}}^{|J_1|+r} \wedge
\overbrace{\mu_{\tau}}^{|J_4|} \wedge
\omega(\zeta)$$
we obtain representations
$$a_{(i,J)}(t,\zeta,z) dt \wedge E(g)(\zeta) \wedge
\lambda^{i,J}_{r-1}(\zeta, z)
= \sum_{\{I,d,j \}}c_{ \{I,d,j \} }(\zeta, z, t)
E(g)(\zeta){\cal K}^{I}_{d,h}(\zeta, z),
\eqno(\arabic{equation})
\newcounter{lambdaRepresentation}
\setcounter{lambdaRepresentation}{\value{equation}}
\addtocounter{equation}{1}$$
and
$$b_{(i,J)}(t,\zeta,z) dt \wedge E(g)(\zeta) \wedge
\gamma^{i,J}_{r-1}(\zeta, z)
= \sum_{\{I,d,j \}}c_{ \{I,d,j \} }(\zeta, z, t)
E(g)(\zeta){\cal K}^{I}_{d,h}(\zeta, z)
\eqno(\arabic{equation})
\newcounter{gammaRepresentation}
\setcounter{gammaRepresentation}{\value{equation}}
\addtocounter{equation}{1}$$
with coefficients $c_{ \{I,d,j \} }$ satisfying
$$\left\|c_{ \{ I,d,j \}}\right\|_k
\leq C(k)\left(1 + |\rho|_{k+3}\right)^n
\eqno(\arabic{equation})
\newcounter{cDetEstimates}
\setcounter{cDetEstimates}{\value{equation}}
\addtocounter{equation}{1}$$
for $k \leq p-3$.\\
\indent
Multiindices $I_i$ for $i=1, \dots, 5$ and indices $d,h$ in
(\arabic{lambdaRepresentation}) satisfy conditions
$$\begin{array}{lllll}
\hspace{0.3in}d = 2(|J_1|+|J_5|+1),\vspace{0.1in}\\
\hspace{0.3in}h = n-|J_1|-|J_5|-1,\vspace{0.1in}\\
\hspace{0.3in}|I_2| + |I_3| = 1+|J_2|+|J_3|+|J_6|,\vspace{0.1in}\\
\hspace{0.3in}|I_4| = |J_1|+|J_4|+r+m-n.
\end{array}
\eqno(\arabic{equation})
\newcounter{KlambdaIndices}
\setcounter{KlambdaIndices}{\value{equation}}
\addtocounter{equation}{1}$$
\indent
Multiindices $I_i$ for $i=1, \dots, 5$ and indices $d,h$ in
(\arabic{gammaRepresentation}) satisfy conditions
$$\begin{array}{lllll}
\hspace{0.3in}d = 2(|J_1|+|J_5|+1),\vspace{0.1in}\\
\hspace{0.3in}h = n-|J_1|-|J_5|-1,\vspace{0.1in}\\
\hspace{0.3in}|I_2| + |I_3| = 2+|J_2|+|J_3|+|J_6|,\vspace{0.1in}\\
\hspace{0.3in}|I_4| = |J_1|+|J_4|+r+m-n.
\end{array}
\eqno(\arabic{equation})
\newcounter{KgammaIndices}
\setcounter{KgammaIndices}{\value{equation}}
\addtocounter{equation}{1}$$
\indent
Using representations (\arabic{lambdaRepresentation}) and
(\arabic{gammaRepresentation}) we reduce the problem
of proving (\arabic{lambdagammaEstimates2}) to each term
$$E(g)(\zeta)\cdot c_{ \{I,d,j \} }(\zeta, z, t)
{\cal K}^{I}_{d,h}(\zeta, z)$$
of the right hand sides of these representations.\\
\indent
In the lemma below we reduce the proof of (\arabic{lambdagammaEstimates1})
and (\arabic{lambdagammaEstimates2}) to the case $k=0$.

\begin{lemma}\label{ReductionTok0}
Statement of Proposition~\ref{REstimates} follows from the corresponding
statement for $k=0$. Namely, for the proof of estimates
(\arabic{lambdagammaEstimates1}) and (\arabic{lambdagammaEstimates2})
it suffices to prove that for a function $h(\zeta,z,t)$ on
${\cal U}_{\zeta}\times U_z\times [0,1]$ and
kernels ${\cal K}^{I}_{d,h}$, obtained from
$\lambda^{i,J}_{r-1}$ and $\gamma^{i,J}_{r-1}$ after application
of Lemma~\ref{Smoothness} the following estimates hold
$$\begin{array}{ll}
\left| \int_{U(\epsilon) \times [0,1]} h(\zeta, z, t)
{\cal K}^{I}_{d,h}(\zeta, z) dt \right|_0 \leq C |h |_0,\vspace{0.1in}\\
\left| \int_{U(\epsilon) \times [0,1]} h(\zeta, z, t)
{\cal K}^{I}_{d,h}(\zeta, z) dt \right|_{1/2} \leq C |h |_{1,z},
\end{array}
\eqno(\arabic{equation})
\newcounter{FinalEstimate}
\setcounter{FinalEstimate}{\value{equation}}
\addtocounter{equation}{1}$$
where $|h |_{1,z}= \sup_{\zeta,t}|h(\zeta,\cdot,t)|_1$.
\end{lemma}

\indent
{\bf Proof.}\hspace{0.05in}
In order to prove Proposition~\ref{REstimates}
we have to prove that for any set of
vector fields $D_1, \dots, D_k \in C^p\left({\cal U},{\cal T}\right)$
such that $\left| D_i \right|_p \leq 1$
$$\begin{array}{ll}
\left| D_1 \circ \dots \circ D_k
\int_{U(\epsilon) \times [0,1]}
E(g)(\zeta)\cdot c_{ \{I,d,j \} }(\zeta, z, t)
{\cal K}^{I}_{d,h}(\zeta, z)dt \right|_0
\leq C(k)\left(1 + |\rho|_{k+3}\right)^{P(k)}
\left| g \right|_{k},\vspace{0.1in}\\
\left| D_1 \circ \dots \circ D_k
\int_{U(\epsilon) \times [0,1]}
E(g)(\zeta)\cdot c_{ \{I,d,j \} }(\zeta, z, t)
{\cal K}^{I}_{d,h}(\zeta, z)dt \right|_{1/2}
\leq C(k)\left(1 + |\rho|_{k+4}\right)^{P(k)}
\left| g \right|_{k}
\end{array}
\eqno(\arabic{equation})
\newcounter{DEstimate}
\setcounter{DEstimate}{\value{equation}}
\addtocounter{equation}{1}$$
with constants $C(k)$ and $P(k)$ independent of $g$ and $\epsilon$.\\
\indent
We apply operator
$${\cal D} = D_1 \circ \dots \circ D_k$$
to integral
$$\int_{U(\epsilon) \times [0,1]}
E(g)(\zeta)\cdot c_{ \{I,d,j \} }(\zeta, z, t)
{\cal K}^{I}_{d,h}(\zeta, z)dt$$
using Lemma~\ref{Smoothness}.\\
\indent
Then we obtain representation
$${\cal D}\int_{U(\epsilon) \times [0,1]}
E(g)(\zeta)\cdot c_{ \{I,d,j \} }(\zeta, z, t)
{\cal K}^{I}_{d,h}(\zeta, z)dt
\eqno(\arabic{equation})
\newcounter{DOperator}
\setcounter{DOperator}{\value{equation}}
\addtocounter{equation}{1}$$
$$= \sum_{\|R\| \leq k} \sum_{\{ R,S,a,b \}}
\int_{U(\epsilon) \times [0,1]}
c_{ \{ R,S,a,b \}}(\zeta, z, t) \cdot
\left[ \{ Y_{\zeta}(z), D_{\zeta} \}^R
E(g)(\zeta) \right]
\cdot {\cal K}^{S}_{a,b}(\zeta, z)dt,$$
where $R=(r_1,r_2)$, $\{ Y_{\zeta}(z), D_{\zeta} \}^R$ denotes
a composition of $r_1$ differentiations $Y_{\zeta}(z)$
and $r_2$ differentiations $D_{i,\zeta}$ with
$\|R\| = r_1 + r_2$,
$$\left\|c_{ \{ R,S,a,b \}}\right\|_0
\leq C(k)\left(1 + |\rho|_{k+3}\right)^{P(k)},
\left\|c_{ \{ R,S,a,b \}}\right\|_1
\leq C(k)\left(1 + |\rho|_{k+4}\right)^{P(k)},
\eqno(\arabic{equation})
\newcounter{cRSEstimates}
\setcounter{cRSEstimates}{\value{equation}}
\addtocounter{equation}{1}$$
and kernels ${\cal K}^{S}_{a,b}$ satisfy (\arabic{SmoothnessIndices}).\\
\indent
Applying then (\arabic{FinalEstimate}) to each term of the right
hand side of (\arabic{DOperator}) and using estimates
$$\begin{array}{ll}
\left|\{ Y_{\zeta}(z), D_{\zeta} \}^R E(g) \right|_0
\leq C(k) \left(1 + |\rho|_{k+3}\right)^{P(k)}\left| g \right|_k,
\vspace{0.1in}\\
\left|\{ Y_{\zeta}(z), D_{\zeta} \}^R E(g) \right|_{1,z}
\leq C(k) \left(1 + |\rho|_{k+4}\right)^{P(k)}\left| g \right|_k,
\end{array}
\eqno(\arabic{equation})
\newcounter{YDEstimates}
\setcounter{YDEstimates}{\value{equation}}
\addtocounter{equation}{1}$$
we obtain (\arabic{DEstimate}).\qed

\indent
The following lemma which is a part of Lemma 3.5 from \cite{P2} will
be used in the proof of (\arabic{FinalEstimate}).\\

\begin{lemma}\label{Integral}
Let
$${\bold B}(\delta) = \{ (\eta, w) \in \R^m \times
\C^{n-m}: \sum_{i=1}^m \eta_i^2+ \sum_{i=1}^{n-m} |w|^2 < \delta \},$$
$${\bold V}(\delta) = \{ (\eta, w) \in \R^m \times
\C^{n-m}: |\eta_1|+\sum_{i=2}^m |\eta_i|^2
+\sum_{i=1}^{n-m} |w|^2 < \delta^2 \},$$
$$K\left\{ k,h\right\}(\eta,w,\epsilon)$$
$$= \frac{\wedge_{i=1}^{m} d\eta_i
\wedge_{i=1}^{n-m} (dw_i \wedge d\bar{w}_i) }
{ (\epsilon + \sum_{i=1}^{m}|\eta_i|+ \sum_{i=1}^{n-m} |w_i|)^{k}
(\sqrt{\epsilon} + \sqrt{|\eta_1|} + \sum_{i=2}^{m}|\eta_i|
+ \sum_{i=1}^{n-m}|w_i|)^{2h}},$$
with $k, 2h \in \Z.$\\
\indent
Let
$${\cal I}_1\left\{k,h\right\}(\epsilon,\delta)
= \int_{{\bold V}(\delta)} K\left\{ k,h\right\}(\eta,w,\epsilon),$$
and
$${\cal I}_2\left\{k,h\right\}(\epsilon,\delta)
= \int_{{\bold B}(1) \setminus {\bold V}(\delta)}
K\left\{k,h\right\}(\eta, w, \epsilon).$$

\indent
Then
$${\cal I}_1\left\{k,h\right\}(\epsilon,\delta)$$
$$= \left\{ \begin{array}{llll}
{\cal O}\left(\epsilon^{2n-m-k-h} \cdot (\log{\epsilon})^2 \right)
& \mbox{if} \hspace{0.05in} k \geq 2n-m-1 \hspace{0.05in}
\mbox{and} \hspace{0.05in} k+h \geq 2n-m,\\
{\cal O}\left(\delta \right)
& \mbox{if} \hspace{0.05in}   k \geq 2n-m-1 \hspace{0.05in}
\mbox{and} \hspace{0.05in} k+h \leq 2n-m-1,\\
{\cal O}\left( \epsilon^{(2n-m-k-2h+1)/2} \cdot \log{\epsilon}\right) &
\mbox{if} \hspace{0.05in}
k \leq 2n-m-2 \hspace{0.05in} \mbox{and} \hspace{0.05in} k+2h \geq 2n-m+1,\\
{\cal O}\left( \delta \right) & \mbox{if} \hspace{0.05in}
k \leq 2n-m-2 \hspace{0.05in} \mbox{and} \hspace{0.05in} k+2h \leq 2n-m,
\end{array} \right.$$
$${\cal I}_2\left\{k,h\right\}(\epsilon,\delta)
= {\cal O}\left(\delta^{-1} \right)\hspace{0.1in}\mbox{if}
\hspace{0.1in}\left\{
\begin{array}{ll}
k \geq 2n-m-1 \hspace{0.05in}\mbox{and}\hspace{0.05in} k+h \leq 2n-m,\\
k \leq 2n-m-2 \hspace{0.05in}\mbox{and}\hspace{0.05in} k+2h \leq 2n-m+2,
\end{array} \right.$$
\end{lemma}

\indent
We will prove estimate (\arabic{FinalEstimate}) as a corollary of the
lemma below.

\begin{lemma}\label{Chalf}
Let $g$ be a function with compact support in
${\cal U}_{\zeta}\times U_z\times[0,1]$ such that
$\left|g\right|_{1,z} < \infty$
and let ${\cal K}^{S}_{a,b}$ satisfy conditions
$$\begin{array}{ll}
k\left({\cal K}^{S}_{a,b}\right)+b
-l\left({\cal K}^{S}_{a,b}\right)\leq 2n-m-2, \vspace{0.1in}\\
k\left({\cal K}^{S}_{a,b}\right)+2b
-2l\left({\cal K}^{S}_{a,b}\right)\leq 2n-m.
\end{array}
\eqno(\arabic{equation})
\newcounter{FinalIndices}
\setcounter{FinalIndices}{\value{equation}}
\addtocounter{equation}{1}$$
\indent
Then
$$f_{\epsilon}(z):= \left( \int_{U(\epsilon)
\times [0,1]} g(\zeta, z, t){\cal K}^{S}_{a,b}(\zeta, z) dt \right)
\in C^{1/2}(U)$$
and
$$\left| f_{\epsilon} \right|_0 \leq C\left(1+|\rho|_3\right)^P
\left| g \right|_0,\hspace{0.1in}
\left| f_{\epsilon} \right|_{1/2} \leq C\left(1+|\rho|_3\right)^P
\left| g \right|_{1,z}
\eqno(\arabic{equation})
\newcounter{fEstimate}
\setcounter{fEstimate}{\value{equation}}
\addtocounter{equation}{1}$$
with $C$ independent of $g$ and $\epsilon$.
\end{lemma}

\indent
{\bf Proof.}\hspace{0.05in}
To prove inclusion
$$f_{\epsilon} \in C^{1/2}(U)$$
we consider for $w \in U$ and arbitrary $\delta > 0$ neighborhoods
$$W(w,\epsilon,\sqrt{\delta})
=\left\{ \zeta \in U: \rho(\zeta) = \epsilon,
\hspace{0.1in}\left| \Phi(\zeta,w)\right| \leq C \cdot \delta \right\},$$
such that for $|z-w| \leq \delta$
$$W(w,\epsilon,c\sqrt{\delta}) \subset {\cal W}(z,\epsilon,\sqrt{\delta})$$
with constants $c,C>0$ independent of $w$, $z$ and $\delta$.\\
\indent
Then we represent $f_{\epsilon}(z)$ as
$$f_{\epsilon}(z) = \int_{{\cal W}(z,\epsilon,\sqrt{\delta}) \times [0,1]}
g(\zeta,z,t){\cal K}^{S}_{a,b}(\zeta, z) dt
\eqno(\arabic{equation})
\newcounter{fRepresentation1}
\setcounter{fRepresentation1}{\value{equation}}
\addtocounter{equation}{1}$$
$$+ \int_{\left( U(\epsilon)
\setminus W(w,\epsilon,\sqrt{\delta}) \right) \times [0,1]}
g(\zeta,z,t){\cal K}^{S}_{a,b}(\zeta, z) dt.$$
\indent
Applying then formula
$$\left. d\rho_i \right|_{U(\epsilon)} = \epsilon d\theta_i$$
and Lemma~\ref{Integral} we obtain for the first term of the right
hand side of (\arabic{fRepresentation1})
$$\left| \int_{{\cal W}(z,\epsilon,\sqrt{\delta}) \times [0,1]}
g(\zeta, z, t) {\cal K}^{S}_{a,b}(\zeta, z) dt \right|$$
$$\leq C | g |_0 \cdot \left( \epsilon^{l({\cal K})}
\cdot \int_{{\cal W}(z,\epsilon,\sqrt{\delta}) \times [0,1]}
\frac{ \overbrace{\wedge_{i} d\theta_i(\zeta)}^{m-1}
\wedge d\sigma_{2n-m}(\zeta,z)} {|\zeta - z|^{k({\cal K})} \cdot
|\Phi(\zeta,z)|^{h({\cal K})} } \right)
\eqno(\arabic{equation})
\newcounter{WEstimate}
\setcounter{WEstimate}{\value{equation}}
\addtocounter{equation}{1}$$
$$\leq C\left(1+|\rho|_3\right)^P| g |_0\left(\epsilon^{l({\cal K})} \cdot
{\cal I}_1\left\{k({\cal K}), h({\cal K})\right\}
\left( \epsilon, \sqrt{\delta} \right) \right)
\leq C\left(1+|\rho|_3\right)^P| g |_0 \cdot \delta^{1/2}.$$
\indent
The same argument proves the first estimate from (\arabic{fEstimate}).\\
\indent
For the second term of the right hand side of (\arabic{fRepresentation1})
using the estimate
$$\left| F^{(k)}(\zeta, z) - F^{(k)}(\zeta, w) \right|
\leq C \left(1+ |\rho|_3\right)\delta
\eqno(\arabic{equation})
\newcounter{Fz-Fw}
\setcounter{Fz-Fw}{\value{equation}}
\addtocounter{equation}{1}$$
for $z,w$ such that $|z-w| < \delta$ and Lemma~\ref{Integral} we obtain
$$\left| \int_{\left( U(\epsilon)
\setminus {\cal W}(z,\epsilon,\sqrt{\delta}) \right) \times [0,1]}
g(\zeta, z, t) {\cal K}^{S}_{a,b}(\zeta, z) dt
- \int_{\left( {\cal U}(\epsilon)
\setminus {\cal W}(z,\epsilon,\sqrt{\delta}) \right) \times [0,1]}
g(\zeta, w, t){\cal K}^{S}_{a,b}(\zeta, w) dt \right|$$
$$\leq \left| \int_{\left( {\cal U}(\epsilon)
\setminus {\cal W}(z,\epsilon,\sqrt{\delta}) \right) \times [0,1]}
g(\zeta,z,t) \left[ {\cal K}^{S}_{a,b}(\zeta,z)
- {\cal K}^{S}_{a,b}(\zeta, w) \right] dt\right|$$
$$+ \left| \int_{\left( {\cal U}(\epsilon)
\setminus {\cal W}(z,\epsilon,\sqrt{\delta}) \right) \times [0,1]}
\left[ g(\zeta, z, t) - g(\zeta, w, t) \right]
{\cal K}^{S}_{a,b}(\zeta, w) dt\right|$$
$$\leq C\left(1+|\rho|_3\right)^P| g |_0\delta
\cdot\left[ {\cal I}_2\left\{ k({\cal K})+1,b-l({\cal K}) \right\}
\left( \epsilon, \sqrt{\delta} \right)
+ {\cal I}_2\left\{ k({\cal K}),
b-l({\cal K})+1 \right\}
\left( \epsilon,\sqrt{\delta} \right)\right]$$
$$+ C\left(1+|\rho|_3\right)^P| g |_{1,z}\delta \cdot
\left[{\cal I}_2\left\{ k({\cal K}),b-l({\cal K}) \right\}
\left( \epsilon,\sqrt{\delta} \right) \right]$$
$$\leq C\left(1+|\rho|_3\right)^P| g |_{1,z}\delta^{1/2}.$$
\indent
Representation (\arabic{fRepresentation1}) together with the estimates
above show that 
$$\left| f_{\epsilon} \right|_{1/2}
\leq C\left(1+|\rho|_3\right)^P| g |_{1,z}$$
uniformly with respect to $\epsilon$.\\
\indent
To prove the first estimate from (\arabic{fEstimate}) we use estimate
(\arabic{WEstimate}) for ${\cal W}(z,\epsilon,1)$.\qed

\indent
In order to complete the proof of Proposition~\ref{REstimates} we
have to prove applicability of Lemma~\ref{Chalf} to the kernels
obtained from $\lambda^{i,J}_{r-1}$ and $\gamma^{i,J}_{r-1}$ after
applications of Lemma~\ref{Smoothness}. We will achieve this goal by
proving relations (\arabic{FinalIndices}) for these kernels.
According to Lemma~\ref{Smoothness} expressions in the left hand
sides of these relations don't increase under transformations from
this lemma. Therefore it suffices to prove relations (\arabic{FinalIndices})
for the original kernels ${\cal K}^{I}_{d,h}(\zeta, z)$ satisfying
conditions (\arabic{KlambdaIndices}) and (\arabic{KgammaIndices}).\\
\indent
Second condition from (\arabic{FinalIndices}) is always satisfied for
the indices satisfying (\arabic{KlambdaIndices})
as can be seen from the inequality
$$k({\cal K})+2h({\cal K})-2l({\cal K})\leq 2n-m-|J_6|\leq 2n-m,
\eqno(\arabic{equation})
\newcounter{SecondIndices}
\setcounter{SecondIndices}{\value{equation}}
\addtocounter{equation}{1}$$
where we used relations
$$\sum_{i=1}^4 |J_{i}|= n-r-1,$$
$$|J_2| + |J_3|\leq m-1,$$
for the multiindices of $\lambda^{i,J}_{r-1}$.\\
\indent
The same arguments show that
condition (\arabic{SecondIndices}) is also satisfied for
the indices defined by (\arabic{KgammaIndices}).\\
\indent
First condition from (\arabic{FinalIndices}) is not satisfied
for all kernels ${\cal K}^{I}_{d,h}(\zeta, z)$. But in the lemma
below we show that if this condition
is not satisfied then the corresponding term
of the integral formula for $R_r(\epsilon)$ does not survive under the
limit when $\epsilon \rightarrow 0$.\\

\begin{lemma}\label{RDoomed}
If $k({\cal K}), h({\cal K}), l({\cal K}) \in \Z$ and
$$k({\cal K})+h({\cal K})-l({\cal K}) \geq 2n-m-1$$
then
$$\left| \int_{U(\epsilon) \times [0,1]}
E(g)(\zeta) c(\zeta,z,t) {\cal K}^{I}_{d,h}(\zeta, z) dt
\right|_k \leq C(k) (\sqrt{\epsilon} \cdot \log{\epsilon})
\left(1 + |\rho|_{k+3}\right)^{P(k)}\left| g \right|_k.
\eqno(\arabic{equation})
\newcounter{DoomedEstimate}
\setcounter{DoomedEstimate}{\value{equation}}
\addtocounter{equation}{1}$$
\end{lemma}
\vspace{0.2in}

\indent
{\bf Proof.}\hspace{0.05in}
Under every application of Lemma~\ref{Smoothness} the quantity
$$k({\cal K})+h({\cal K})-l({\cal K})$$
doesn't increase, therefore using representation (\arabic{DOperator})
and estimates (\arabic{cRSEstimates}) and (\arabic{YDEstimates})
from Lemma~\ref{ReductionTok0} we reduce the
statement of the lemma to the case $k=0$:
$$\left| \int_{U(\epsilon) \times [0,1]}
E(g)(\zeta) c(\zeta,z,t) {\cal K}^{I}_{d,h}(\zeta, z) dt
\right|_0 \leq C(\sqrt{\epsilon} \cdot \log{\epsilon})
\left(1+|\rho|_3\right)^P\left| g \right|_0.
\eqno(\arabic{equation})
\newcounter{Doomed0}
\setcounter{Doomed0}{\value{equation}}
\addtocounter{equation}{1}$$
\indent
To prove estimate (\arabic{Doomed0}) we use inequality
$$2n-m+ l({\cal K}) - k({\cal K}) - h({\cal K})
= n-|J_1|-|J_5|+|J_6|-1 \geq 1,
\eqno(\arabic{equation})
\newcounter{FirstIndices}
\setcounter{FirstIndices}{\value{equation}}
\addtocounter{equation}{1}$$
which is a corollary of definitions of $k({\cal K})$, $h({\cal K})$
and $l({\cal K})$, equality
$$\sum_{i=1}^5 |J_{i}| = n-r-1$$
and inequality
$$n-1-|J_1|-|J_5|\geq 1.$$
\indent
From the condition of the lemma and inequality (\arabic{FirstIndices})
we obtain
$$k({\cal K})+h({\cal K})-l({\cal K}) = 2n-m-1$$
and
$$n-|J_1|-|J_5|-1 = 1,$$
which leads to
$$|J_1|=n-r-1,\hspace{0.1in}|J_3|=0,\hspace{0.1in}|J_4|=0,
\hspace{0.1in}|J_5|=r-1,$$
and
$$l({\cal K}) \geq |J_1|+|J_4|+r+m-n = m-1 \geq 1.$$
\indent
Using Lemma~\ref{Integral} to estimate the integral in the left hand side
of (\arabic{Doomed0}) we obtain
$$\left| \int_{U(\epsilon) \times [0,1]}
E(g)(\zeta) c(\zeta,z,t)
{\cal K}^{I}_{d,h}(\zeta, z) dt \right|_0$$
$$\leq \left| g \right|_0 \cdot
\epsilon^{l({\cal K})} \cdot {\cal O}
\left( {\cal I}_1\left\{ k({\cal K}), h({\cal K}), 0\right\}
\left( \epsilon, 1 \right) \right)$$
$$\leq \left| g \right|_0 \cdot
\left\{ \begin{array}{ll}
\epsilon^{l({\cal K})} \cdot
{\cal O}\left(\epsilon^{2n-m-k({\cal K})-h({\cal K})}
\cdot (\log{\epsilon})^2 \right) & \mbox{if} \hspace{0.05in}
k({\cal K}) \geq 2n-m-1,
\vspace{0.1in}\\
\epsilon^{l({\cal K})} \cdot
{\cal O}\left( \epsilon^{(2n-m-k({\cal K})-2h({\cal K})+1)/2}
\cdot \log{\epsilon}\right) & \mbox{if} \hspace{0.05in}
k({\cal K}) \leq 2n-m-2.
\end{array} \right.$$
\indent
Using then inequality (\arabic{FirstIndices}) in the first subcase of the above
and inequality (\arabic{SecondIndices}) in the second subcase we obtain
estimate (\arabic{Doomed0}).\qed\\
\indent
This completes the proof of Proposition~\ref{REstimates}.\\

\indent
In the proposition below we refine estimate (\arabic{REstimate1}) for a
special case $k \ll p$.\\

\begin{proposition}\label{Rk+sEstimates}
Let ${\bold M} \subset {\bold G}$ be a generic, regular q-pseudoconcave
CR submanifold of the class $C^p$ in ${\cal U}$ satisfying condition
(\arabic{NonDegeneracy}). Let $s,k \in \Z$ be such that $s\leq k$ and
$k+s\leq p-3$.\\
\indent
Then $R_r(g)$ defined in (\arabic{LocalFormula}) satisfies the
estimate
$$\left| R_r(g) \right|_{k+s} \leq
C(k)\left(1 + |\rho|_{k+3}\right)^{P(k)}
\Bigl[|g|_{k+s}+\left(1+|\rho|_{k+s+3}\right)|g|_k\Bigr],
\eqno(\arabic{equation})
\newcounter{Rk+sEstimate}
\setcounter{Rk+sEstimate}{\value{equation}}
\addtocounter{equation}{1}$$
with $P(k)$ a polynomial in $k$ and a constant $C(k)$ independent of $g$.
\end{proposition}

\indent
{\bf Proof.}\ Proof of estimate (\arabic{Rk+sEstimate}) is analogous to the
proof of (\arabic{REstimate1}). Namely, we inductively use Lemma~\ref{Smoothness}
and reduce the statement of the Proposition to the estimate
(\arabic{FinalEstimate}). The only difference is that we consider separately
two groups of terms: with derivatives of $g$ of order higher that $k$ and the rest.
For the terms with derivatives of $g$ of higher order, derivatives of functions
$c_{ \{ S,a,b \}}$ and $c_{ \{ L,i,e,k \}}$, appearing in the
Lemma~\ref{Smoothness} will be of the order lower than $k$. Therefore,
using estimate (\arabic{cEstimates}) and estimate (\arabic{fEstimate}) from
Lemma~\ref{Chalf}, we obtain that these terms are dominated
by the first term of the right hand side of (\arabic{Rk+sEstimate}).\\
\indent
For the terms from the second group we have to estimate the derivatives
of functions $c_{ \{ S,a,b \}}$ and $c_{ \{ L,i,e,k \}}$ of the higher
order but derivatives of $g$ wiil be of the order, less or equal to $k$.
Using estimates (\arabic{cEstimates}) and (\arabic{fEstimate}) we
obtain the second term of the right hand side of (\arabic{Rk+sEstimate}).\\
\indent
Lemma~\ref{RDoomed} assures that only the terms with "good" indices have
to be estimated.
\qed

\indent
To complete the proof of Proposition~\ref{AlmostHomotopy} we have to prove
estimates (\arabic{AlmostEstimate1}) and (\arabic{AlmostEstimate2}) for
operator ${\bold H}^r_{\bold M}$.
From the definition of operator ${\bold H}^r_{\bold M}$ we conclude that
it suffices to prove these estimates for each of the terms below
$$\bar\partial_{\bold M} \vartheta^{\prime}_{\iota}(z)
\wedge R_r^{\iota}(\vartheta_{\iota}g)(z),\hspace{0.1in}
\vartheta^{\prime}_{\iota}(z) \cdot
R_{r+1}^{\iota}(\bar\partial_{\bold M} \vartheta_{\iota} \wedge g)(z)
\hspace{0.1in}\mbox{and}\hspace{0.1in}
\vartheta^{\prime}_{\iota}(z) \cdot H_r^{\iota}(\vartheta_{\iota}g)(z).$$
\indent
Estimates of the first two of these terms follow from the corresponding
estimates of operators $R_r$ proved in Proposition~\ref{REstimates}.
The proposition below takes care of the third term
of ${\bold H}^r_{\bold M}$.

\begin{lemma}\label{HZero}
Let $r < q$. Then
$$H_r(g)(z) = 0.
\eqno(\arabic{equation})
\newcounter{Hzero}
\setcounter{Hzero}{\value{equation}}
\addtocounter{equation}{1}$$
\end{lemma}
\vspace{0.2in}

\indent
{\bf Proof.}\hspace{0.05in}
Using approximation of $H_r$ by the operators
$$H_r(\epsilon)(g)(z) = (-1)^{r} \frac{(n-1)!}{(2\pi i)^n}
\cdot \mbox{pr}_{\bold M} \circ
\int_{U(\epsilon)} \vartheta(\zeta)E(g)(\zeta)
\wedge \omega^{\prime}_{r} \left( \frac{P(\zeta,z)}
{\Phi(\zeta,z)}\right) \wedge\omega(\zeta)$$
we conclude that it suffices to prove equality
$$\omega^{\prime}_{r} \left( \frac{P(\zeta,z)}
{\Phi(\zeta,z)}\right) \wedge\omega(\zeta) =0
\eqno(\arabic{equation})
\newcounter{Omegazero}
\setcounter{Omegazero}{\value{equation}}
\addtocounter{equation}{1}$$
for $r < q$.\\
\indent
This kernel with the use of (\arabic{muandQ})
may be represented on ${\cal U} \times {\cal U}$ as
$$\left. \omega^{\prime}_{r}\left(\frac{P(\zeta,z)}
{\Phi(\zeta,z)}\right) \wedge\omega(\zeta)
\right|_{{\cal U} \times {\cal U} }
\eqno(\arabic{equation})
\newcounter{cauchykernel}
\setcounter{cauchykernel}{\value{equation}}
\addtocounter{equation}{1}$$
$$= \sum_{i,J} a_{(i,J)}(\zeta,z) \wedge
{\widetilde \phi}^{i,J}_{r}(\zeta, z) +
\sum_{i,J} b_{(i,J)}(\zeta,z) \wedge
{\widetilde \psi}^{i,J}_{r}(\zeta, z),$$
where $i$ is an index, $J= \cup_{i=1}^{6} J_i$ is a
multiindex such that $i \not \in J,$
$a_{(i,J)}(\zeta,z)$ and $b_{(i,J)}(\zeta,z)$ are $C^p$-functions of
$z$, $\zeta$ and $\theta(\zeta)$, and ${\widetilde \phi}^{i,J}_{r}(\zeta, z)$
and ${\widetilde \psi}^{i,J}_{r}(\zeta, z)$ are defined as follows:
$${\widetilde \phi}^{i,J}_{r}(\zeta, z) =\frac{1}
{{\Phi(\zeta,z)}^{n}} \times \mbox{Det} \left[
Q^{(i)},\hspace{0.03in}
\overbrace{{\bar A} \cdot \bar\partial_{\zeta}a}^{j \in J_1},\hspace{0.03in}
\overbrace{a \cdot \mu_{\nu},\hspace{0.03in}\chi}^{j \in J_2},\hspace{0.03in}
\overbrace{a \cdot \mu_{\tau}}^{j \in J_3}, \hspace{0.03in}
\overbrace{{\bar A} \cdot \bar \partial_z a}^{j \in J_4},\hspace{0.03in}
\overbrace{a \cdot \bar\partial_z {\bar A}}^{j \in J_5},\hspace{0.03in}
\overbrace{\bar\partial_z Q}^{j \in J_6}
\right] \wedge\omega(\zeta)
\eqno(\arabic{equation})
\newcounter{phiDet}
\setcounter{phiDet}{\value{equation}}
\addtocounter{equation}{1}$$
and
$${\widetilde \psi}^{i,J}_{r}(\zeta, z) =\frac{1}
{{\Phi(\zeta,z)}^{n}} \times \mbox{Det} \left[
a_i {\bar A_i},\hspace{0.03in}
\overbrace{{\bar A} \cdot \bar\partial_{\zeta}a}^{j \in J_1},\hspace{0.03in}
\overbrace{a \cdot \mu_{\nu},\hspace{0.03in}\chi}^{j \in J_2},\hspace{0.03in}
\overbrace{a \cdot \mu_{\tau}}^{j \in J_3},\hspace{0.03in}
\overbrace{{\bar A} \cdot \bar \partial_z a}^{j \in J_4},\hspace{0.03in}
\overbrace{a \cdot \bar\partial_z {\bar A}}^{j \in J_5},\hspace{0.03in}
\overbrace{\bar\partial_z Q}^{j \in J_6}
\right] \wedge\omega(\zeta).
\eqno(\arabic{equation})
\newcounter{psiDet}
\setcounter{psiDet}{\value{equation}}
\addtocounter{equation}{1}$$
\indent
Multiindices of ${\widetilde \phi}^{i,J}_r$ and ${\widetilde \psi}^{i,J}_r$
satisfy the following conditions
$$\begin{array}{ll}
\sum_{i=1}^3|J_i|=n-r-1,\vspace{0.1in}\\
|J_1|+|J_2|\leq m-1,
\end{array}
\eqno(\arabic{equation})
\newcounter{phiIndices}
\setcounter{phiIndices}{\value{equation}}
\addtocounter{equation}{1}$$
therefore, if $r < q$ then
$$|J_3| = n-r-1-|J_1|-|J_2|\geq n-r-m > n-q-m,$$
which is impossible.\qed

\section{Homotopy formula.}\label{Homotopy}

\indent
Our proof of Theorem~\ref{HomotopyTheorem} is based on a transformation
of operators ${\bold H}^r_{\bold M}$, ${\bold R}^r_{\bold M}$ and
${\bold R}^{r+1}_{\bold M}$ on manifolds ${\bold M}$ close
to a fixed manifold ${\bold M}_0$ such that $\dim H^r\left({\bold M}_0,
{\cal B}\big|_{{\bold M}_0}\right)=0$. We start with the construction of a
homotopy formula on ${\bold M}_0$.

\begin{proposition}\label{M0Formula}
Let ${\bold M}_0 \subset {\bold G}$ be a compact, regular
q-pseudoconcave submanifold of the class $C^p$ and let ${\cal B}$
be a holomorphic vector bundle on ${\bold G}$.
Let for some $r < q$ equality $\dim H^r\left({\bold M}_0,
{\cal B}\big|_{{\bold M}_0}\right)=0$ be satisfied.
Then for $k \leq p-4$ there exist $s,l \in \Z^{+}$,
linear continuous functionals
$\left\{\alpha_j\right\}_1^s$ on $C^{1}_{(0,r)}\left({\bold M}_0,
{\cal B}\big|_{{\bold M}_0}\right)$, linear continuous functionals
$\left\{\beta_j\right\}_1^l$ on $C^{1}_{(0,r+1)}\left({\bold M}_0,
{\cal B}\big|_{{\bold M}_0}\right)$, collections
$\left\{g_j\right\}_1^s \in C^{p-4}_{(0,r-1)}\left({\bold M}_0,
{\cal B}\big|_{{\bold M}_0}\right)$ and $\left\{f_j\right\}_1^l
\in C^{p-4}_{(0,r)}\left({\bold M}_0,
{\cal B}\big|_{{\bold M}_0}\right)$, and operators
$$\begin{array}{ll}
{\bold P}^r_{{\bold M}_0}(u)={\bold R}^r_{{\bold M}_0}
\left(u-\sum_{i=1}^l \beta_i(\bar\partial_{{\bold M}_0}u)f_i
-\sum_{i=1}^s \alpha_i(u)\bar\partial_{{\bold M}_0}g_i\right)
+\sum_{i=1}^s \alpha_i(u)g_i,\vspace{0.1in}\\
{\bold P}^{r+1}_{{\bold M}_0}(v)={\bold R}^{r+1}_{{\bold M}_0}
\left(v-\sum_{i=1}^l \beta_i(v)\bar\partial_{{\bold M}_0}f_i\right)
+\sum_{i=1}^l \beta_i(v)f_i,
\end{array}\eqno(\arabic{equation})
\newcounter{Poperators}
\setcounter{Poperators}{\value{equation}}
\addtocounter{equation}{1}$$
such that
$$\left|{\bold P}^r_{{\bold M}_0}\right|_k,
\left|{\bold P}^{r+1}_{{\bold M}_0}\right|_k \leq C(k),
\eqno(\arabic{equation})
\newcounter{boldPEstimates}
\setcounter{boldPEstimates}{\value{equation}}
\addtocounter{equation}{1}$$
and operator ${\bold F}_{{\bold M}_0}
=\bar\partial_{{\bold M}_0}{\bold P}^r_{{\bold M}_0} +
{\bold P}^{r+1}_{{\bold M}_0}\bar\partial_{{\bold M}_0}$
defines an isomorphism on $C^k_{(0,r)}\left({\bold M}_0,
{\cal B}\big|_{{\bold M}_0}\right)$.
\end{proposition}
\vspace{0.2in}

\indent
{\bf Proof.}\hspace{0.05in}
We consider the Fredholm operator from (\arabic{AlmostFormula})
$${\bold A}_{{\bold M}_0}={\bold I}-{\bold H}^r_{{\bold M}_0}
=\bar\partial_{{\bold M}_0}{\bold R}^r_{{\bold M}_0}
+{\bold R}^{r+1}_{{\bold M}_0}\bar\partial_{{\bold M}_0}:
C^k_{(0,r)}\left({\bold M}_0,{\cal B}\big|_{{\bold M}_0}\right)
\to C^k_{(0,r)}\left({\bold M}_0,{\cal B}\big|_{{\bold M}_0}\right)$$
and construct a subspace ${\cal N}^{(k)}_{{\bold M}_0}$ of finite
codimension in $C^k_{(0,r)}\left({\bold M}_0,
{\cal B}\big|_{{\bold M}_0}\right)$ such that the restriction of
${\bold A}_{{\bold M}_0}$ to ${\cal N}^{(k)}_{{\bold M}_0}$
is invertible. To construct such a subspace we will use the lemma below.

\begin{lemma} \label{Stabilization}
Let $A: {\cal B} \longrightarrow {\cal B}$ be a Fredholm operator
on a Banach space ${\cal B}$ of the form $A=I-H$ with compact $H$.
Then the sequence of subspaces
$$\cdots \supseteq \mbox{Ker}\left(A^n\right)\supseteq
\mbox{Ker}\left(A^{n-1}\right) \cdots
\supseteq \mbox{Ker}A
\eqno(\arabic{equation})
\newcounter{sequence}
\setcounter{sequence}{\value{equation}}
\addtocounter{equation}{1}$$
stabilizes on a finite step.
\end{lemma}

\indent
{\bf Proof.}\hspace{0.1in}
Let us assume that sequence (\arabic{sequence}) doesn't stabilize on
a finite step. Then there exist $\delta > 0$ and a sequence
$\left\{ x_n \right\} \in {\cal B}$ such that
$$\|x_n\|=1,\hspace{0.1in}x_n \in K_n \hspace{0.1in}\mbox{and}
\hspace{0.1in}\mbox{dist}\left(x_n, K_{n-1} \right) > \delta,$$
where $K_n = \mbox{Ker}\left(A^n\right)$ are finite-dimensional subspaces
in ${\cal B}$.\\
\indent
Using then Hahn-Banach theorem we can construct a sequence of
linear continuous functionals $\left\{l_n\right\}$ on ${\cal B}$ such that
$$\|l_n\| \leq 1/\delta,\hspace{0.1in}l_n(x_n) = 1,
\hspace{0.1in}\mbox{and}\hspace{0.1in}l_n(x) = 0
\hspace{0.1in}\mbox{for}\hspace{0.1in} x \in K_{n-1}.
\eqno(\arabic{equation})
\newcounter{lproperties}
\setcounter{lproperties}{\value{equation}}
\addtocounter{equation}{1}$$
\indent
Since $H$ is compact we can assume that sequence $\{ H(x_n) \}$
converges in ${\cal B}$, and therefore
$\lim_{n \to \infty}\left\| H(x_n -x_{n-1}) \right\|=0$.\\
\indent
Using then the estimate
$$\left| l_n \left( H(x_n -x_{n-1}) \right) \right| \leq \|l_n\| \cdot
\left\| H(x_n -x_{n-1}) \right\|$$
we obtain that $\left| l_n \left( H(x_n -x_{n-1}) \right) \right| \to 0$
when $n \to \infty$.\\
\indent
But on the other hand
$$l_n \left( H(x_n - x_{n-1}) \right)
= l_n \left( x_n - A(x_n) - x_{n-1} + A(x_{n-1}) \right)
= l_n \left( x_n \right) = 1,$$
where we used properties (\arabic{lproperties}) and
$$A\left( K_n \right) \subseteq K_{n-1},
\hspace{0.1in}A\left( K_{n-1} \right) \subseteq K_{n-1}.$$
\indent
Obtained contradiction proves the lemma.
\qed

Applying Lemma~\ref{Stabilization} to ${\bold A}_{{\bold M}_0}$
on $C^1_{(0,r)}\left({\bold M}_0,{\cal B}\big|_{{\bold M}_0}\right)$
we find $n_0$ such that
$$\mbox{Ker}\left({\bold A}_{{\bold M}_0}^{n_0+i}\right)
\cap C^1_{(0,r)}\left({\bold M}_0,{\cal B}\big|_{{\bold M}_0}\right)
=\mbox{Ker}\left({\bold A}_{{\bold M}_0}^{n_0}\right)
\cap C^1_{(0,r)}\left({\bold M}_0,{\cal B}\big|_{{\bold M}_0}\right)
\hspace{0.05in}\mbox{for}\hspace{0.05in}i\geq 1$$
and therefore
$$\mbox{Ker}\left({\bold A}_{{\bold M}_0}^{n_0}\right)
\cap\mbox{Im}\left({\bold A}_{{\bold M}_0}^{n_0}\right)
=\emptyset$$
for the restriction of ${\bold A}_{{\bold M}_0}^{n_0}$ to
$C^1_{(0,r)}\left({\bold M}_0,{\cal B}\big|_{{\bold M}_0}\right)$,
and hence on any $C^k_{(0,r)}\left({\bold M}_0,
{\cal B}\big|_{{\bold M}_0}\right)$ with $1\leq k\leq p-4$.\\
\indent
Using the estimates for operator ${\bold H}^r_{{\bold M}_0}$
we conclude that the restriction of the identity operator to
$\mbox{Ker}\left({\bold A}_{{\bold M}_0}^{n_0}\right)
\cap C^k_{(0,r)}\left({\bold M},{\cal B}\big|_{\bold M}\right)$
is smoothing from
$C^k_{(0,r)}\left({\bold M}_0,{\cal B}\big|_{{\bold M}_0}\right)$ to
$C^{k+1/2}_{(0,r)}\left({\bold M}_0,{\cal B}\big|_{{\bold M}_0}\right)$
for any $k\leq p-4$. To prove higher smoothness of the elements of
$\mbox{Ker}\left({\bold A}_{{\bold M}_0}^{n_0}\right)
\cap C^k_{(0,r)}\left({\bold M},{\cal B}\big|_{\bold M}\right)$
we use the following interpolation result, which follows from \cite{Kr},
(cf. \cite{LP}).

\begin{proposition} \label{Interpolation}
Let $0<\epsilon<1/2$ and let
$$L: C^{\infty}_{(0,r)}\left({\bold M}_0,{\cal B}\big|_{{\bold M}_0}
\right)\to C^{\infty}_{(0,r)}\left({\bold M}_0,
{\cal B}\big|_{{\bold M}_0}\right)$$
be a linear operator, satisfying
$$\begin{array}{ll}
|L(f)|_{k+1/2} \leq C\cdot|f|_{k+\epsilon},\vspace{0.1in}\\
|L(f)|_{k+3/2} \leq C\cdot|f|_{k+1+\epsilon}.
\end{array}$$
for $k\in\Z^+$.\\
\indent
Then for $0\leq \alpha \leq 1$ the estimate
$$|L(f)|_{k+\alpha+1/2} \leq C\cdot|f|_{k+\alpha+\epsilon}$$
holds.
\end{proposition}
\qed

\indent
Starting with $C^{k+1/2}_{(0,r)}\left({\bold M}_0,
{\cal B}\big|_{{\bold M}_0}\right)$ and consequtively $2(p-k-4)$ times
applying Proposition~\ref{Interpolation}
for small enough $\epsilon$ we obtain that the identity operator on
$\mbox{Ker}\left({\bold A}_{{\bold M}_0}^{n_0}\right)
\cap C^k_{(0,r)}\left({\bold M},{\cal B}\big|_{\bold M}\right)$
is smoothing from $C^k_{(0,r)}\left({\bold M}_0,
{\cal B}\big|_{{\bold M}_0}\right)$ to
$C^{p-4}_{(0,r)}\left({\bold M}_0,{\cal B}\big|_{{\bold M}_0}\right)$.
Defining then for $k\leq p-4$
$${\cal N}^{(k)}_{{\bold M}_0}={\bold A}_{{\bold M}_0}^{n_0}
\left[C^k_{(0,r)}\left({\bold M}_0,{\cal B}\big|_{{\bold M}_0}\right)
\right]\hspace{0.05in}\mbox{and}\hspace{0.05in}
{\cal K}_{{\bold M}_0}
=\mbox{Ker}\left({\bold A}_{{\bold M}_0}^{n_0}\right)
\cap C^{p-4}_{(0,r)}\left({\bold M},{\cal B}\big|_{\bold M}\right)$$
we obtain that for $k\leq p-4$
\begin{itemize}
\item[(i)]
$C^k_{(0,r)}\left({\bold M}_0,{\cal B}\big|_{{\bold M}_0}\right)
={\cal N}^{(k)}_{{\bold M}_0}\oplus{\cal K}_{{\bold M}_0}$,
\item[(ii)]
$\dim\hspace{0.02in}{\cal K}_{{\bold M}_0}<\infty$,
\item[(iii)]
restriction of ${\bold A}_{{\bold M}_0}$ to ${\cal N}^{(k)}_{{\bold M}_0}$
is an isomorphism.
\end{itemize}
\indent
Using condition $\dim H^r\left({\bold M}_0,
{\cal B}\big|_{{\bold M}_0}\right)=0$ we can find finitely many
forms $g_1,\dots,g_s$ and construct a finite basis
$$\left\{\bar\partial_{{\bold M}_0}g_1,\dots,
\bar\partial_{{\bold M}_0}g_s,f_1,\dots,f_l\right\}$$
of ${\cal K}_{{\bold M}_0}$ with
$\left\{g_j\right\}_1^s \in
C^{p-4}_{(0,r-1)}\left({\bold M}_0,{\cal B}\big|_{{\bold M}_0}\right)$
and $\left\{f_j\right\}_1^l \in
C^{p-4}_{(0,r)}\left({\bold M}_0,{\cal B}\big|_{{\bold M}_0}\right)$
so that
$$\begin{array}{lll}
{\cal S}_{{\bold M}_0}=\left\{h\in{\cal K}_{{\bold M}_0}:
\bar\partial_{{\bold M}_0}h=0 \right\}
=\mbox{Span}\left\{\bar\partial_{{\bold M}_0}g_1,\dots,
\bar\partial_{{\bold M}_0}g_s\right\},\vspace{0.1in}\\
{\cal P}_{{\bold M}_0}=\mbox{Span}\left\{f_1,\dots,
f_l\right\},\vspace{0.1in}\\
{\cal K}_{{\bold M}_0}={\cal S}_{{\bold M}_0}\oplus
{\cal P}_{{\bold M}_0}.
\end{array}
\eqno(\arabic{equation})
\newcounter{KandS}
\setcounter{KandS}{\value{equation}}
\addtocounter{equation}{1}$$
Applying then Hahn-Banach theorem we construct linear continuous
functionals
$$\left\{\alpha_j\right\}_1^s\hspace{0.05in}\mbox{on}
\hspace{0.05in}C^1_{(0,r)}\left({\bold M}_0,
{\cal B}\big|_{{\bold M}_0}\right)$$
such that
$$\begin{array}{ll}
\alpha_j\left(\bar\partial_{{\bold M}_0}g_i\right)= \delta_i^j,
\vspace{0.1in}\\
\alpha_j\left(h\right\}=0\hspace{0.05in}\mbox{for}\hspace{0.05in}
h \in {\cal N}^{(1)}_{{\bold M}_0}\oplus{\cal P}_{{\bold M}_0}.
\end{array}
\eqno(\arabic{equation})
\newcounter{alphaFunctionals}
\setcounter{alphaFunctionals}{\value{equation}}
\addtocounter{equation}{1}$$
\indent
In the construction of linear continuous functionals
$\left\{\beta_j\right\}_1^l$ from (\arabic{Poperators}) we will need
the following lemma.

\begin{lemma} \label{Separation}
$$\bar\partial_{{\bold M}_0}\left({\cal P}_{{\bold M}_0}\right)
\cap \mbox{cl}\hspace{0.02in}\left\{\bar\partial_{{\bold M}_0}
\left({\cal N}^{(2)}_{{\bold M}_0}\oplus {\cal S}_{{\bold M}_0}
\right)\right\}=\left\{0\right\},
\eqno(\arabic{equation})
\newcounter{KcapN}
\setcounter{KcapN}{\value{equation}}
\addtocounter{equation}{1}$$
where $\mbox{cl}\hspace{0.02in}\left\{\bar\partial_{{\bold M}_0}
\left({\cal N}^{(2)}_{{\bold M}_0}\oplus{\cal S}_{{\bold M}_0}
\right)\right\}$ is the closure in $C^{1}_{(0,r+1)}\left({\bold M}_0,
{\cal B}\big|_{{\bold M}_0}\right)$ of the image of
${\cal N}^{(2)}_{{\bold M}_0}\oplus{\cal S}_{{\bold M}_0}$ under
$\bar\partial_{{\bold M}_0}$.
\end{lemma}

\indent
{\bf Proof.}\hspace{0.1in}
Using (\arabic{KandS}) and (\arabic{alphaFunctionals}) we redefine
operator ${\bold R}^r_{{\bold M}_0}$ by the formula
$${\bold R}^r_{{\bold M}_0}\left(u-\sum_{i=1}^s\alpha_i(u)
\bar\partial_{{\bold M}_0}g_i\right)
+\sum_{i=1}^s\alpha_i(u)g_i.$$
Then the restriction of a new ${\bold A}_{{\bold M}_0}
=\bar\partial_{{\bold M}_0}{\bold R}^r_{{\bold M}_0}
+{\bold R}^{r+1}_{{\bold M}_0}\bar\partial_{{\bold M}_0}$ to
${\cal N}^{(k)}_{{\bold M}_0}\oplus{\cal S}_{{\bold M}_0}$
coincides with the old one on ${\cal N}^{(k)}_{{\bold M}_0}$
and with identity on ${\cal S}_{{\bold M}_0}$.\\
\indent
Therefore there exist bounded linear operators
$${\bold B}^{(k)}_{{\bold M}_0}:{\cal N}^{(k)}_{{\bold M}_0}\oplus
{\cal S}_{{\bold M}_0} \to
{\cal N}^{(k)}_{{\bold M}_0}\oplus{\cal S}_{{\bold M}_0}$$
such that
$$\begin{array}{lll}
{\bold B}^{(k)}_{{\bold M}_0}\circ {\bold A}_{{\bold M}_0}
= {\bold A}_{{\bold M}_0}\circ{\bold B}^{(k)}_{{\bold M}_0} = {\bold I},
\vspace{0.1in}\\
{\bold B}^{(k)}_{{\bold M}_0}\left({{\cal N}^{(k)}_{{\bold M}_0}
\oplus {\cal S}_{{\bold M}_0}}\right)
= {{\cal N}^{(k)}_{{\bold M}_0}\oplus {\cal S}_{{\bold M}_0}},
\vspace{0.1in}\\
\left|{\bold B}^{(k)}_{{\bold M}_0}\right| \leq C(k).
\end{array}
\eqno(\arabic{equation})
\newcounter{BProperties}
\setcounter{BProperties}{\value{equation}}
\addtocounter{equation}{1}$$
\indent
We denote
$${\cal Z}^{(k)}_{{\bold M}_0}=\left\{h \in
C^k_{(0,r)}\left({\bold M}_0,{\cal B}\big|_{{\bold M}_0}\right):
\bar\partial_{{\bold M}_0}h=0\right\},$$
and consider actions of ${\bold A}_{{\bold M}_0}$ and
${\bold B}^{(k)}_{{\bold M}_0}$ on ${\cal Z}^{(k)}_{{\bold M}_0}$.
For any $h \in {\cal Z}^{(k)}_{{\bold M}_0}$ we have
$${\bold A}_{{\bold M}_0}(h)
=\bar\partial_{{\bold M}_0}{\bold R}^r_{{\bold M}_0}(h)
\in {\cal Z}^{(k)}_{{\bold M}_0},$$
and thus
$${\bold A}_{{\bold M}_0}\left({\cal Z}^{(k)}_{{\bold M}_0}\right)
\subseteq {\cal Z}^{(k)}_{{\bold M}_0}\hspace{0.1in}\mbox{and}
\hspace{0.1in}{\bold A}^{n_0}_{{\bold M}_0}
\left({\cal Z}^{(k)}_{{\bold M}_0}\right)
\subseteq {\cal Z}^{(k)}_{{\bold M}_0}.$$
Using (\arabic{KandS}) we represent $h \in {\cal Z}^{(k)}_{{\bold M}_0}$ as
$$h=h_{\cal P}+h_{\cal S}+h_{\cal N}$$
with $h_{\cal P}\in {\cal P}_{{\bold M}_0}$,
$h_{\cal S}\in {\cal S}_{{\bold M}_0}$ and
$h_{\cal N}\in {\cal N}^{(k)}_{{\bold M}_0}$,
and obtain
$${\bold A}^{n_0}_{{\bold M}_0}(h)={\bold A}^{n_0}_{{\bold M}_0}
\left(h_{\cal N}\right)+h_{\cal S},$$
with
$${\bold A}^{n_0}_{{\bold M}_0}\left(h_{\cal N}\right)
\in {\cal N}^{(k)}_{{\bold M}_0}\cap {\cal Z}^{(k)}_{{\bold M}_0}.$$
\indent
Since ${\bold A}^{n_0}_{{\bold M}_0}$ defines an isomorphism on
${\cal N}^{(k)}_{{\bold M}_0}$ and ${\bold A}^{n_0}_{{\bold M}_0}
\left({\cal Z}^{(k)}_{{\bold M}_0}\right)
\subseteq {\cal Z}^{(k)}_{{\bold M}_0}$, we conclude that
$h_{\cal N}\in {\cal N}^{(k)}_{{\bold M}_0}
\cap {\cal Z}^{(k)}_{{\bold M}_0}$ and, therefore,
$${\cal Z}^{(k)}_{{\bold M}_0}
=\left({\cal N}^{(k)}_{{\bold M}_0}\cap {\cal Z}^{(k)}_{{\bold M}_0}\right)
\oplus{\cal S}_{{\bold M}_0}.
\eqno(\arabic{equation})
\newcounter{ZandN}
\setcounter{ZandN}{\value{equation}}
\addtocounter{equation}{1}$$
\indent
From (\arabic{ZandN}) we obtain that for $h \in {\cal Z}^{(k)}_{{\bold M}_0}$
$${\bold A}_{{\bold M}_0}(h)
=\bar\partial_{{\bold M}_0}{\bold R}^r_{{\bold M}_0}(h)
\in {\cal Z}^{(k)}_{{\bold M}_0}
\subset {\cal N}^{(k)}_{{\bold M}_0}\oplus{\cal S}_{{\bold M}_0},$$
and thus
$${\bold A}_{{\bold M}_0}\left({\cal Z}^{(k)}_{{\bold M}_0}\right)
= {\cal Z}^{(k)}_{{\bold M}_0},\hspace{0.1in}\mbox{and}
\hspace{0.1in}{\bold B}^{(k)}_{{\bold M}_0}
\left({\cal Z}^{(k)}_{{\bold M}_0}\right)
= {\cal Z}^{(k)}_{{\bold M}_0}.
\eqno(\arabic{equation})
\newcounter{ABofZ}
\setcounter{ABofZ}{\value{equation}}
\addtocounter{equation}{1}$$
\indent
Then for any
$h \in {\cal N}^{(k)}_{{\bold M}_0}\oplus{\cal S}_{{\bold M}_0}$
we will have
$${\bold A}_{{\bold M}_0}(h)
=\bar\partial_{{\bold M}_0}{\bold R}^r_{{\bold M}_0}(h)
+{\bold R}^{r+1}_{{\bold M}_0}(\bar\partial_{{\bold M}_0}h)
\in {\cal N}^{(k)}_{{\bold M}_0}\oplus{\cal S}_{{\bold M}_0}
\subset {\cal N}^{(k-1)}_{{\bold M}_0}\oplus{\cal S}_{{\bold M}_0}
\eqno(\arabic{equation})
\newcounter{AonN}
\setcounter{AonN}{\value{equation}}
\addtocounter{equation}{1}$$
with
$$\bar\partial_{{\bold M}_0}{\bold R}^r_{{\bold M}_0}(h)
\in {\cal Z}^{(k-1)}_{{\bold M}_0} \subset
{\cal N}^{(k-1)}_{{\bold M}_0}\oplus {\cal S}_{{\bold M}_0},$$
and thus
$${\bold R}^{r+1}_{{\bold M}_0}(\bar\partial_{{\bold M}_0}h)
\in {\cal N}^{(k-1)}_{{\bold M}_0}\oplus {\cal S}_{{\bold M}_0}.$$
\indent
Applying ${\bold B}^{(k-1)}_{{\bold M}_0}$ to both parts of
(\arabic{AonN}) and using the last two inclusions we obtain for
$h \in {\cal N}^{(k)}_{{\bold M}_0}\oplus {\cal S}_{{\bold M}_0}$
$$h={\bold B}^{(k-1)}_{{\bold M}_0}\left(\bar\partial_{{\bold M}_0}
{\bold R}^r_{{\bold M}_0}
+{\bold R}^{r+1}_{{\bold M}_0}\bar\partial_{{\bold M}_0}\right)(h)
\eqno(\arabic{equation})
\newcounter{BA}
\setcounter{BA}{\value{equation}}
\addtocounter{equation}{1}$$
$$= {\bold B}^{(k-1)}_{{\bold M}_0}\bar\partial_{{\bold M}_0}
{\bold R}^r_{{\bold M}_0}(h)
+ {\bold B}^{(k-1)}_{{\bold M}_0}{\bold R}^{r+1}_{{\bold M}_0}
(\bar\partial_{{\bold M}_0}h).$$
\indent
From (\arabic{ABofZ}) we conclude that
$$\bar\partial_{{\bold M}_0}{\bold B}^{(k-1)}_{{\bold M}_0}
\bar\partial_{{\bold M}_0}{\bold R}^r_{{\bold M}_0}(h)=0,$$
and
$${\bold B}^{(k-1)}_{{\bold M}_0}\bar\partial_{{\bold M}_0}
{\bold R}^r_{{\bold M}_0}(h)={\bold A}_{{\bold M}_0}
{\bold B}^{(k-1)}_{{\bold M}_0}\left[{\bold B}^{(k-1)}_{{\bold M}_0}
\bar\partial_{{\bold M}_0}{\bold R}^r_{{\bold M}_0}(h)\right]
\eqno(\arabic{equation})
\newcounter{BinZ}
\setcounter{BinZ}{\value{equation}}
\addtocounter{equation}{1}$$
$$=\bar\partial_{{\bold M}_0}{\bold R}^r_{{\bold M}_0}
{\bold B}^{(k-1)}_{{\bold M}_0}{\bold B}^{(k-1)}_{{\bold M}_0}
\bar\partial_{{\bold M}_0}{\bold R}^r_{{\bold M}_0}(h)
=\bar\partial_{{\bold M}_0}{\bold C}_{{\bold M}_0}(h),$$
with
$${\bold C}_{{\bold M}_0}={\bold R}^r_{{\bold M}_0}
{\bold B}^{(k-1)}_{{\bold M}_0}{\bold B}^{(k-1)}_{{\bold M}_0}
\bar\partial_{{\bold M}_0}{\bold R}^r_{{\bold M}_0}:
{\cal N}^{(k)}_{{\bold M}_0}\oplus {\cal S}_{{\bold M}_0}\to
C^{k-1}_{(0,r)}\left({\bold M}_0,{\cal B}\big|_{{\bold M}_0}\right).$$
\indent
Denoting then
$${\bold D}_{{\bold M}_0}={\bold B}^{(k-1)}_{{\bold M}_0}
{\bold R}^{r+1}_{{\bold M}_0},$$
we rewrite (\arabic{BA}) for
$h \in {\cal N}^{(k)}_{{\bold M}_0}\oplus{\cal S}_{{\bold M}_0}$
as
$$h=\bar\partial_{{\bold M}_0}{\bold C}_{{\bold M}_0}(h)
+{\bold D}_{{\bold M}_0}(\bar\partial_{{\bold M}_0}h).
\eqno(\arabic{equation})
\newcounter{CandD}
\setcounter{CandD}{\value{equation}}
\addtocounter{equation}{1}$$
\indent
To complete now the proof of the Lemma let us assume that there
exists $\psi \neq 0$, such that $\psi \in {\cal P}_{{\bold M}_0}$,
and a sequence $\left\{\bar\partial_{{\bold M}_0}\phi_i\right\}_1^{\infty}
\in C^{1}_{(0,r+1)}\left({\bold M}_0,{\cal B}\big|_{{\bold M}_0}\right)$
such that
$$\begin{array}{ll}
\phi_i \in {\cal N}^{(2)}_{{\bold M}_0}\oplus{\cal S}_{{\bold M}_0},
\vspace{0.1in}\\
\lim_{i\to\infty}\bar\partial_{{\bold M}_0}\phi_i
=\bar\partial_{{\bold M}_0}\psi\hspace{0.05in}\mbox{in}\hspace{0.05in}
C^{1}_{(0,r+1)}\left({\bold M}_0,{\cal B}\big|_{{\bold M}_0}\right).
\end{array}$$
From the continuity of operator ${\bold D}_{{\bold M}_0}$ we
conclude that sequence $\left\{{\bold D}_{{\bold M}_0}
(\bar\partial_{{\bold M}_0}\phi_i)\right\}_1^{\infty}$
also converges to some $\gamma \in C^{1}_{(0,r)}\left({\bold M}_0,
{\cal B}\big|_{{\bold M}_0}\right)$.\\
\indent
Using inclusions
$\phi_i \in {\cal N}^{(2)}_{{\bold M}_0}\oplus{\cal S}_{{\bold M}_0}$
and ${\bold D}_{{\bold M}_0}(\bar\partial_{{\bold M}_0}\phi_i)
\in C^{1}_{(0,r)}\left({\bold M}_0,{\cal B}\big|_{{\bold M}_0}\right)$
and representation (\arabic{CandD}) we obtain that
$\bar\partial_{{\bold M}_0}{\bold C}_{{\bold M}_0}(\phi_i) \in
{\cal Z}^{(1)}_{{\bold M}_0} \subset
{\cal N}^{(1)}_{{\bold M}_0}\oplus {\cal S}_{{\bold M}_0}$
and therefore
$${\bold D}_{{\bold M}_0}(\bar\partial_{{\bold M}_0}\phi_i)
\in {\cal N}^{(1)}_{{\bold M}_0}\oplus{\cal S}_{{\bold M}_0}.$$
\indent
Using then closedness of ${\cal N}^{(1)}_{{\bold M}_0}
\oplus{\cal S}_{{\bold M}_0}$ we obtain that
$\gamma \in {\cal N}^{(1)}_{{\bold M}_0}\oplus{\cal S}_{{\bold M}_0}$
and $\bar\partial_{{\bold M}_0}\gamma = \bar\partial_{{\bold M}_0}\psi$.
Condition $\bar\partial_{{\bold M}_0}\left(\gamma-\psi\right)=0$
implies $\gamma -\psi \in {\cal Z}^{(1)}_{{\bold M}_0} \subset
{\cal N}^{(1)}_{{\bold M}_0}\oplus{\cal S}_{{\bold M}_0}$
and therefore $\psi \in {\cal N}^{(1)}_{{\bold M}_0}
\oplus{\cal S}_{{\bold M}_0}$, which contradicts inclusion
$\psi \in {\cal P}_{{\bold M}_0}$.\qed\\

\indent
To conclude the proof of Proposition~\ref{M0Formula} we use condition
(\arabic{KcapN}), apply the Hahn-Banach
theorem to $C^{1}_{(0,r+1)}\left({\bold M}_0,
{\cal B}\big|_{{\bold M}_0}\right)$ and construct linear continuous
functionals $\left\{\beta_i\right\}_1^l$ such that
$$\beta_i(\bar\partial_{{\bold M}_0}f_j)=\delta_i^j,
\hspace{0.05in}\mbox{and}\hspace{0.05in}
\beta_i\left(\bar\partial_{{\bold M}_0}h\right)=0
\hspace{0.05in}\mbox{for}\hspace{0.05in}h \in
{\cal N}^{(2)}_{{\bold M}_0}\oplus{\cal S}_{{\bold M}_0}.$$
Then operator ${\bold F}_{{\bold M}_0}=\bar\partial_{{\bold M}_0}
{\bold P}^r_{{\bold M}_0}+{\bold P}^{r+1}_{{\bold M}_0}
\bar\partial_{{\bold M}_0}$ with ${\bold P}^r_{{\bold M}_0}$
and ${\bold P}^{r+1}_{{\bold M}_0}$ from (\arabic{Poperators}) will
satisfy the following conditions:
\begin{itemize}

\item[(i)]\hspace{0.05in}${\bold F}_{{\bold M}_0}(h)
={\bold A}_{{\bold M}_0}(h)$ for $h \in {\cal N}^{(k)}_{{\bold M}_0}$,
$(k=2,\dots,p-4)$,
\vspace{0.1in}

\item[(ii)]\hspace{0.05in}${\bold F}_{{\bold M}_0}
(\bar\partial_{{\bold M}_0}g_i)=\bar\partial_{{\bold M}_0}g_i$,
\vspace{0.1in}

\item[(iii)]\hspace{0.05in}${\bold F}_{{\bold M}_0}(f_i)=f_i$,
\end{itemize}
and will therefore define an isomorphism on
$C^k_{(0,r)}\left({\bold M}_0,{\cal B}\big|_{{\bold M}_0}\right)$.\\
\indent
Estimates (\arabic{boldPEstimates}) follow from the corresponding
estimates for operators ${\bold R}^r_{{\bold M}_0}$ and
${\bold R}^{r+1}_{{\bold M}_0}$ and from the boundedness of linear functionals
$\left\{\alpha_i\right\}_1^s$ and $\left\{\beta_i\right\}_1^l$.\qed

\indent
The proposition below, which can be considered as an
analogue of the Hodge-Kohn decomposition for $C^{p-4}_{(0,r)}\left({\bold M}_0,
{\cal B}\big|_{{\bold M}_0}\right)$, is proved using basically the same approach
as in the proof of Proposition~\ref{M0Formula}.

\begin{proposition}\label{M0Decomposition}
Let ${\bold M}_0 \subset {\bold G}$ be a compact, regular
q-pseudoconcave submanifold of the class $C^p$ and let ${\cal B}$
be a holomorphic vector bundle on ${\bold G}$.
Then for fixed $r < q$ and $1< k \leq p-4$ there exist a finite-dimensional
linear operator
$${\bold K}^r_{{\bold M}_0}:C^k_{(0,r)}\left({\bold M}_0,
{\cal B}\big|_{{\bold M}_0}\right)\to C^k_{(0,r)}\left({\bold M}_0,
{\cal B}\big|_{{\bold M}_0}\right),$$
and linear operators
${\bold Q}^r_{{\bold M}_0}$ and ${\bold Q}^{r+1}_{{\bold M}_0}$
such that
$$\begin{array}{ll}
\left|{\bold Q}^r_{{\bold M}_0}(h)\right|_{k-1} \leq C|h|_k,\vspace{0.1in}\\
\left|{\bold Q}^{r+1}_{{\bold M}_0}(h)\right|_{k-1} \leq C|h|_k,
\end{array}
\eqno(\arabic{equation})
\newcounter{boldQEstimates}
\setcounter{boldQEstimates}{\value{equation}}
\addtocounter{equation}{1}$$
and equality
$$h=\bar\partial_{{\bold M}_0}{\bold Q}^r_{{\bold M}_0}(h)
+{\bold Q}^{r+1}_{{\bold M}_0}\left(\bar\partial_{{\bold M}_0}h\right)
+{\bold K}^r_{{\bold M}_0}(h)
\eqno(\arabic{equation})
\newcounter{QKequality}
\setcounter{QKequality}{\value{equation}}
\addtocounter{equation}{1}$$
is satisfied for any $h \in C^k_{(0,r)}\left({\bold M}_0,
{\cal B}\big|_{{\bold M}_0}\right)$.\\
\indent
If for $h \in C^k_{(0,r)}\left({\bold M}_0,
{\cal B}\big|_{{\bold M}_0}\right)$ there exists
$g \in C^k_{(0,r-1)}\left({\bold M}_0,{\cal B}\big|_{{\bold M}_0}\right)$
such that
$\bar\partial_{{\bold M}_0}g =h$, then ${\bold K}^r_{{\bold M}_0}(h)=0$.
\end{proposition}
\vspace{0.1in}

\indent
{\bf Proof.}\ As in the proof of Proposition~\ref{M0Formula} we use
Lemma~\ref{Stabilization} and Proposition~\ref{Interpolation} and
obtain a more general decomposition
$$\begin{array}{llllll}
{\cal S}^{(k)}_{{\bold M}_0}=\left\{h \in {\cal K}_{{\bold M}_0}:
\exists g \in C^k_{(0,r-1)}\left({\bold M}_0,{\cal B}\big|_{{\bold M}_0}\right)
\ \mbox{s. t.}\ \bar\partial_{{\bold M}_0}g=h\right\}
=\mbox{Span}\left\{\bar\partial_{{\bold M}_0}g_1,\dots,
\bar\partial_{{\bold M}_0}g_s\right\},\vspace{0.1in}\\
{\cal L}^{(k)}_{{\bold M}_0}=\mbox{Span}\left\{v_1,\dots,v_t\right\},
\vspace{0.1in}\\
{\cal S}^{(k)}_{{\bold M}_0}\oplus {\cal L}^{(k)}_{{\bold M}_0}
=\left\{h\in{\cal K}_{{\bold M}_0}:\bar\partial_{{\bold M}_0}h=0 \right\},
\vspace{0.1in}\\
{\cal P}_{{\bold M}_0}=\mbox{Span}\left\{f_1,\dots,
f_l\right\},\vspace{0.1in}\\
{\cal K}_{{\bold M}_0}={\cal S}^{(k)}_{{\bold M}_0}
\oplus{\cal L}^{(k)}_{{\bold M}_0}
\oplus{\cal P}_{{\bold M}_0},\vspace{0.1in}\\
C^{k}_{(0,r)}\left({\bold M}_0,{\cal B}\big|_{{\bold M}_0}\right)
={\cal N}^{(k)}_{{\bold M}_0}\oplus{\cal K}_{{\bold M}_0},
\end{array}
\eqno(\arabic{equation})
\newcounter{KLandS}
\setcounter{KLandS}{\value{equation}}
\addtocounter{equation}{1}$$
with
$\left\{g_j\right\}_1^s \in
C^k_{(0,r-1)}\left({\bold M}_0,{\cal B}\big|_{{\bold M}_0}\right)$.\\
\indent
Applying then Hahn-Banach theorem we construct linear continuous
functionals
$$\left\{\alpha^{(k)}_j\right\}_1^s\ \mbox{and}\ \left\{\gamma^{(k)}_j\right\}_1^l\
\mbox{on}\ C^k_{(0,r)}\left({\bold M}_0,
{\cal B}\big|_{{\bold M}_0}\right)$$
such that
$$\begin{array}{lll}
\alpha^{(k)}_j\left(\bar\partial_{{\bold M}_0}g_i\right)= \delta_i^j,\
\gamma^{(k)}_j\left(v_i\right)= \delta_i^j,
\vspace{0.1in}\\
\alpha^{(k)}_j\left(h\right\}=0\ \mbox{for}\
h \in {\cal N}^{k)}_{{\bold M}_0}\oplus{\cal L}^{(k)}_{{\bold M}_0}
\oplus{\cal P}_{{\bold M}_0},\vspace{0.1in}\\
\gamma^{(k)}_j\left(h\right\}=0\ \mbox{for}\
h \in {\cal N}^{(k)}_{{\bold M}_0}\oplus{\cal S}^{(k)}_{{\bold M}_0}
\oplus{\cal P}_{{\bold M}_0}.
\end{array}
\eqno(\arabic{equation})
\newcounter{alphagammaFunctionals}
\setcounter{alphagammaFunctionals}{\value{equation}}
\addtocounter{equation}{1}$$
\indent
Redefining as in Lemma~\ref{Separation} operator
${\bold R}^r_{{\bold M}_0}$ for a fixed $1 < k \leq p-4$ by the formula
$${\bold R}^r_{{\bold M}_0}\left(u-\sum_{i=1}^s\alpha^{(k)}_i(u)
\bar\partial_{{\bold M}_0}g_i-\sum_{i=1}^l\gamma^{(k)}_i(u)v_i\right)
+\sum_{i=1}^s\alpha^{(k)}_i(u)g_i.$$
we obtain that the new ${\bold A}_{{\bold M}_0}
=\bar\partial_{{\bold M}_0}{\bold R}^r_{{\bold M}_0}
+{\bold R}^{r+1}_{{\bold M}_0}\bar\partial_{{\bold M}_0}$
coincides with the old one on ${\cal N}^{(k)}_{{\bold M}_0}$ for
$1 < k \leq p-4$, is identity on ${\cal S}^{(k)}_{{\bold M}_0}$, and
is zero on ${\cal L}^{(k)}_{{\bold M}_0}$. Thus operator ${\bold A}_{{\bold M}_0}$
can be extended to an isomorphism ${\bold F}_{{\bold M}_0}$ on
${\cal N}^{(k)}_{{\bold M}_0}\oplus{\cal S}^{(k)}_{{\bold M}_0}
\oplus{\cal L}^{(k)}_{{\bold M}_0}$ by setting it as identity on
${\cal L}^{(k)}_{{\bold M}_0}$.\\
\indent
Denoting
$${\cal Z}^{(k)}_{{\bold M}_0}=\left\{h \in C^k_{(0,r)}\left({\bold M}_0,
{\cal B}\big|_{{\bold M}_0}\right):\bar\partial_{{\bold M}_0}h=0\right\},$$
we obtain that ${\bold F}^{(k)}_{{\bold M}_0}$ defines an isomorphism of
${\cal N}^{(k)}_{{\bold M}_0}\oplus{\cal S}^{(k)}_{{\bold M}_0}
\oplus{\cal L}^{(k)}_{{\bold M}_0}$ into itself, preserving
${\cal Z}^{(k)}_{{\bold M}_0}$.\\
\indent
Therefore, there exists a bounded operator
$${\bold T}^{(k)}_{{\bold M}_0}:{\cal N}^{(k)}_{{\bold M}_0}
\oplus{\cal S}^{(k)}_{{\bold M}_0}\oplus{\cal L}^{(k)}_{{\bold M}_0}
\rightarrow {\cal N}^{(k)}_{{\bold M}_0}\oplus{\cal S}^{(k)}_{{\bold M}_0}
\oplus{\cal L}^{(k)}_{{\bold M}_0},$$
preserving ${\cal N}^{(k)}_{{\bold M}_0}$ and ${\cal Z}^{(k)}_{{\bold M}_0}$.
Using then that for 
$h \in {\cal N}^{(k)}_{{\bold M}_0} \subset {\cal N}^{(k-1)}_{{\bold M}_0}
\oplus{\cal S}^{(k-1)}_{{\bold M}_0}\oplus{\cal L}^{(k-1)}_{{\bold M}_0}$
we have $\bar\partial_{{\bold M}_0}{\bold R}^r_{{\bold M}_0}(h) \in
{\cal Z}^{(k-1)}_{{\bold M}_0}$, we obtain as in (\arabic{BA}) for
$h \in {\cal N}^{(k)}_{{\bold M}_0}$
$$h={\bold T}^{(k-1)}_{{\bold M}_0}\circ{\bold F}^{(k-1)}_{{\bold M}_0}(h)
={\bold T}^{(k-1)}_{{\bold M}_0}\left(\bar\partial_{{\bold M}_0}
{\bold R}^r_{{\bold M}_0}(h)\right)
+{\bold T}^{(k-1)}_{{\bold M}_0}\circ{\bold R}^{r+1}_{{\bold M}_0}
\left(\bar\partial_{{\bold M}_0}h\right).
\eqno(\arabic{equation})
\newcounter{TF}
\setcounter{TF}{\value{equation}}
\addtocounter{equation}{1}$$
\indent
For $h \in {\cal N}^{(k)}_{{\bold M}_0}$
we consider a representation of
$f=\bar\partial_{{\bold M}_0}{\bold R}^r_{{\bold M}_0}(h)
\in C^{k-1}_{(0,r-1)}\left({\bold M}_0,{\cal B}\big|_{{\bold M}_0}\right)$ as
$$f=f_{\cal N}+f_{\cal S}+f_{\cal L},$$
with $f_{\cal N} \in {\cal N}^{(k-1)}_{{\bold M}_0}$,
$f_{\cal S} \in {\cal S}^{(k-1)}_{{\bold M}_0}$, and
$f_{\cal L} \in {\cal L}^{(k-1)}_{{\bold M}_0}$, and from (\arabic{TF})
conclude that $g={\bold R}^{r+1}_{{\bold M}_0}
\left(\bar\partial_{{\bold M}_0}h\right)$ also admits analogous
representation
$$g=g_{\cal N}+g_{\cal S}+g_{\cal L}.$$
\indent
From the definition of operators ${\bold T}^{(k)}_{{\bold M}_0}$
we conclude that equality (\arabic{TF}) holds after projecting all terms
of this equality on ${\cal N}^{(k-1)}_{{\bold M}_0}$
$$h=f_{\cal N}+g_{\cal N}=\mbox{pr}_{\cal N}\circ{\bold T}^{(k)}_{{\bold M}_0}
\left(\bar\partial_{{\bold M}_0}{\bold R}^r_{{\bold M}_0}(h)\right)
+\mbox{pr}_{\cal N}\circ{\bold T}^{(k-1)}_{{\bold M}_0}
\circ{\bold R}^{r+1}_{{\bold M}_0}\left(\bar\partial_{{\bold M}_0}h\right),$$
where $\mbox{pr}_{\cal N}$ is the projection on ${\cal N}^{(k-1)}_{{\bold M}_0}$.\\
\indent
As in the proof of Lemma~\ref{Separation} we have
$$\bar\partial_{{\bold M}_0}\left(\mbox{pr}_{\cal N}\circ
{\bold T}^{(k-1)}_{{\bold M}_0}\left(\bar\partial_{{\bold M}_0}
{\bold R}^r_{{\bold M}_0}(h)\right)\right)=0,$$
and
$$\mbox{pr}_{\cal N}\circ{\bold T}^{(k-1)}_{{\bold M}_0}
\bar\partial_{{\bold M}_0}{\bold R}^r_{{\bold M}_0}(h)={\bold F}_{{\bold M}_0}
{\bold T}^{(k-1)}_{{\bold M}_0}\left[\mbox{pr}_{\cal N}\circ
{\bold T}^{(k-1)}_{{\bold M}_0}
\bar\partial_{{\bold M}_0}{\bold R}^r_{{\bold M}_0}(h)\right]
\eqno(\arabic{equation})
\newcounter{TinZ}
\setcounter{TinZ}{\value{equation}}
\addtocounter{equation}{1}$$
$$=\bar\partial_{{\bold M}_0}{\bold R}^r_{{\bold M}_0}
{\bold T}^{(k-1)}_{{\bold M}_0}\left[\mbox{pr}_{\cal N}\circ
{\bold T}^{(k-1)}_{{\bold M}_0}\bar\partial_{{\bold M}_0}
{\bold R}^r_{{\bold M}_0}(h)\right]
=\bar\partial_{{\bold M}_0}{\bold C}^r_{{\bold M}_0}(h),$$
with
$${\bold C}^r_{{\bold M}_0}={\bold R}^r_{{\bold M}_0}
{\bold T}^{(k-1)}_{{\bold M}_0}\mbox{pr}_{\cal N}{\bold T}^{(k-1)}_{{\bold M}_0}
\bar\partial_{{\bold M}_0}{\bold R}^r_{{\bold M}_0}:
C^{k}_{(0,r)}\left({\bold M}_0,{\cal B}\big|_{{\bold M}_0}\right)\to
C^{k-1}_{(0,r)}\left({\bold M}_0,{\cal B}\big|_{{\bold M}_0}\right).$$
\indent
Denoting then ${\bold D}^{r+1}_{{\bold M}_0}
=\mbox{pr}_{\cal N}{\bold T}^{(k-1)}_{{\bold M}_0}\circ
{\bold R}^{r+1}_{{\bold M}_0}$, we can rewrite (\arabic{TF})
for $h \in {\cal N}^{(k)}_{{\bold M}_0}$ as
$$h=\bar\partial_{{\bold M}_0}{\bold C}^r_{{\bold M}_0}(h)
+{\bold D}^{r+1}_{{\bold M}_0}\left(\bar\partial_{{\bold M}_0}h\right).
\eqno(\arabic{equation})
\newcounter{Nhomotopy}
\setcounter{Nhomotopy}{\value{equation}}
\addtocounter{equation}{1}$$
\indent
Using (\arabic{Nhomotopy}) as in the proof of Lemma~\ref{Separation} we
obtain
$$\bar\partial_{{\bold M}_0}\left({\cal P}_{{\bold M}_0}\right)
\cap \mbox{cl}\hspace{0.02in}\left\{\bar\partial_{{\bold M}_0}
\left({\cal N}^{(k)}_{{\bold M}_0}\oplus{\cal S}^{(k)}_{{\bold M}_0}
\oplus{\cal L}^{(k)}_{{\bold M}_0}\right)\right\}=\left\{0\right\},
\eqno(\arabic{equation})
\newcounter{PcapNSL}
\setcounter{PcapNSL}{\value{equation}}
\addtocounter{equation}{1}$$
where $\mbox{cl}\hspace{0.02in}\left\{\bar\partial_{{\bold M}_0}
\left({\cal N}^{(k)}_{{\bold M}_0}\oplus{\cal S}^{(k)}_{{\bold M}_0}
\oplus{\cal L}^{(k)}_{{\bold M}_0}\right)\right\}$ is the closure in
$C^{k-1}_{(0,r+1)}\left({\bold M}_0,{\cal B}\big|_{{\bold M}_0}\right)$
of the image of ${\cal N}^{(k)}_{{\bold M}_0}\oplus{\cal S}^{(k)}_{{\bold M}_0}
\oplus{\cal L}^{(k)}_{{\bold M}_0}$ under $\bar\partial_{{\bold M}_0}$.
Equality (\arabic{PcapNSL}) allows us to construct linear continuous
functionals $\left\{\beta^{(k-1)}_i\right\}_1^l$ on
$C^{k-1}_{(0,r+1)}\left({\bold M}_0,{\cal B}\big|_{{\bold M}_0}\right)$
such that
$$\beta^{(k-1)}_i(\bar\partial_{{\bold M}_0}f_j)=\delta_i^j,
\hspace{0.05in}\mbox{and}\hspace{0.05in}
\beta^{(k-1)}_i\left(\bar\partial_{{\bold M}_0}h\right)=0
\ \mbox{for}\ h \in {\cal N}^{(k)}_{{\bold M}_0}\oplus
{\cal S}^{(k)}_{{\bold M}_0}\oplus{\cal L}^{(k)}_{{\bold M}_0}.$$
\indent
We consider then operators
$$\begin{array}{lll}
{\bold P}^r_{{\bold M}_0}(u)={\bold R}^r_{{\bold M}_0}
\left(u-\sum_{i=1}^l \beta^{(k-1)}_i(\bar\partial_{{\bold M}_0}u)f_i
-\sum_{i=1}^s \alpha^{(k)}_i(u)\bar\partial_{{\bold M}_0}g_i
-\sum_{i=1}^t\gamma^{(k)}_i(u)v_i\right)\\
\hspace{0.7in}+\sum_{i=1}^s \alpha^{(k)}_i(u)g_i,\vspace{0.1in}\\
{\bold P}^{r+1}_{{\bold M}_0}(v)={\bold R}^{r+1}_{{\bold M}_0}
\left(v-\sum_{i=1}^l \beta^{(k-1)}_i(v)\bar\partial_{{\bold M}_0}f_i\right)
+\sum_{i=1}^l \beta^{(k-1)}_i(v)f_i,\vspace{0.1in}\\
{\bold K}^r_{{\bold M}_0}(u)=\sum_{i=1}^t \gamma^{(k)}_i(u)v_i
\end{array}\eqno(\arabic{equation})
\newcounter{AlmostOperators}
\setcounter{AlmostOperators}{\value{equation}}
\addtocounter{equation}{1}$$
and operator
$${\bold L}_{{\bold M}_0}=\bar\partial_{{\bold M}_0}{\bold P}^r_{{\bold M}_0}
+{\bold P}^{r+1}_{{\bold M}_0}\bar\partial_{{\bold M}_0},$$
satisfying the following conditions
$$\begin{array}{ll}
{\bold L}_{{\bold M}_0}\Big|_{{\cal N}^{(k)}}
={\bold A}_{{\bold M}_0},\ {\bold L}_{{\bold M}_0}\Big|_{{\cal S}^{(k)}}
={\bold I},\ {\bold L}_{{\bold M}_0}\Big|_{{\cal L}^{(k)}}=0,\vspace{0.1in}\\
{\bold K}_{{\bold M}_0}\Big|_{{\cal N}^{(k)}\oplus{\cal S}^{(k)}}=0,\
{\bold K}_{{\bold M}_0}\Big|_{{\cal L}^{(k)}}={\bold I}.
\end{array}\eqno(\arabic{equation})
\newcounter{LKproperties}
\setcounter{LKproperties}{\value{equation}}
\addtocounter{equation}{1}$$
\indent
From (\arabic{LKproperties}) we conclude that the restriction of
${\bold L}_{{\bold M}_0}$ to ${\cal N}^{(k)}_{{\bold M}_0}$ is
an isomorphism and in order to construct
necessary operators ${\bold Q}^r_{{\bold M}_0}$ and
${\bold Q}^{r+1}_{{\bold M}_0}$ we have to modify operators
${\bold P}^r_{{\bold M}_0}$ and ${\bold P}^{r+1}_{{\bold M}_0}$
only on ${\cal N}^{(k)}_{{\bold M}_0}$. Such a modification has been made
already for operator ${\bold A}_{{\bold M}_0}$ in (\arabic{TinZ}) and
(\arabic{Nhomotopy}). Namely, using (\arabic{TinZ}), we obtain for
$h \in {\cal N}^{(k)}_{{\bold M}_0}$
$$\mbox{pr}_{\cal N}\circ{\bold T}^{(k-1)}_{{\bold M}_0}
\bar\partial_{{\bold M}_0}{\bold P}^r_{{\bold M}_0}(h)={\bold F}_{{\bold M}_0}
{\bold T}^{(k-1)}_{{\bold M}_0}\left[\mbox{pr}_{\cal N}\circ
{\bold T}^{(k-1)}_{{\bold M}_0}
\bar\partial_{{\bold M}_0}{\bold P}^r_{{\bold M}_0}(h)\right]$$
$$=\bar\partial_{{\bold M}_0}{\bold P}^r_{{\bold M}_0}
{\bold T}^{(k-1)}_{{\bold M}_0}\left[\mbox{pr}_{\cal N}\circ
{\bold T}^{(k-1)}_{{\bold M}_0}\bar\partial_{{\bold M}_0}
{\bold P}^r_{{\bold M}_0}(h)\right]
=\bar\partial_{{\bold M}_0}{\bold Q}^r_{{\bold M}_0}(h),$$
with
$${\bold Q}^r_{{\bold M}_0}={\bold P}^r_{{\bold M}_0}
{\bold T}^{(k-1)}_{{\bold M}_0}\mbox{pr}_{\cal N}{\bold T}^{(k-1)}_{{\bold M}_0}
\bar\partial_{{\bold M}_0}{\bold P}^r_{{\bold M}_0}:
C^{k}_{(0,r)}\left({\bold M}_0,{\cal B}\big|_{{\bold M}_0}\right)\to
C^{k-1}_{(0,r)}\left({\bold M}_0,{\cal B}\big|_{{\bold M}_0}\right).$$
\indent
Defining then ${\bold Q}^{r+1}_{{\bold M}_0}
=\mbox{pr}_{\cal N}{\bold T}^{(k-1)}_{{\bold M}_0}\circ
{\bold P}^{r+1}_{{\bold M}_0}$ we obtain equality (\arabic{QKequality})
with operators ${\bold Q}^r_{{\bold M}_0}$ and ${\bold Q}^{r+1}_{{\bold M}_0}$
satisfying estimates (\arabic{boldQEstimates}).
\qed\\

\indent
In order to prove existence of homotopy operators ${\bold P}_{\bold M}$ for
a manifold ${\bold M}$ close enough to ${\bold M}_0$ we will make several
assumptions on the local structure
of such a manifold and its defining functions. We assume the existence
of a finite cover $\left\{{\cal U}^{\iota}\right\}_1^N$ of some
neighborhood of ${\bold M}_0$ in ${\bold G}$ such that in each
${\cal U}^{\iota}$ the manifold ${\bold M}_0 \cap {\cal U}^{\iota}$ has
the form (\arabic{manifold}) with defining functions
$\left\{\rho^{(0)}_{\iota,l}\right\}_{1\leq l \leq m}$. We also
assume the existence of $C^p$-diffeomorphisms
${\cal F}:{\bold M}_0\to {\bold M}$ such that
${\cal F}(z) = z+f^{\iota}(z)$ in $U^{\iota}_0
={\cal U}^{\iota}\cap {\bold M}_0$
with $f^{\iota} \in \left[C^p(U^{\iota}_0)\right]^n$.\\
\indent
For ${\cal F}$ close enough to identity or
$\left|f\right|_p = \max_{\iota}\{\left|f_{\iota}\right|_p\}$
small enough the inverse map ${\cal G} = {\cal F}^{-1}$
has the form ${\cal G}(z) = z+g^{\iota}(z)$ with
$g^{\iota} \in \left[C^p\left(V^{\iota}\right)\right]^n$
for some neighborhood
${\cal V}^{\iota} \subset {\cal F}\left(U^{\iota}\right)$ and
$V^{\iota} = {\cal V}^{\iota} \cap {\bold M}$.
The proposition below shows that the collection
$\left\{{\cal V}^{\iota} \subset {\cal U}^{\iota}\right\}_1^N$
may be chosen so that it also covers some neighborhood of ${\bold M}_0$.\\
\indent
We denote by $\B(r)$ the ball in $\C^n$ of radius $r$ centered at the origin
and for a function $f:\B(r) \rightarrow \C$ we denote
$$\left|f\right|_{r,p} \equiv \sup_{z \in \B(r), |J|\leq p}
\left|\partial^{J}f(z)\right|,$$
where $J=(j_1,\dots,j_{2n})$ is a multiindex, $|J|=j_1+\dots+j_{2n}$,
and $\partial^{J} \equiv
\partial^{j_1}_{x_1}\cdots \partial^{j_{2n}}_{x_{2n}}$ with coordinates
$\left\{x_j\right\}_1^{2n}$ such that $z_j=x_j+\sqrt{-1}x_{n+j}$.
For a $C^p$-smooth vector function $f:\B(r) \rightarrow \C^n$ and
$k \leq p$ we denote
$$\left|f\right|_{r,k} \equiv
\sup_{1\leq i \leq n}\left|f_i\right|_{r,k}.$$

\begin{proposition} \label{Inverse}(cf. \cite{W2})
Let $F_i(z) = z_i+f_i(z)$ for $i=1,\dots,n$, and let the
functions $\{ f_i\}_1^n \in C^p(\B(1))$ satisfy the estimate
$|f_i|_{1,1} < \epsilon$.\\
\indent
Then for small enough $\epsilon$ and fixed $s,k \in \Z$, such
that $0\leq s\leq k,\ k+s\leq p$ there exist a constant $C(k)$ and
a set of functions
$$\{ g_i\}_1^n \in C^p\left(\B(1-2\epsilon)\right)$$
such that $G(z) \equiv z+g(z)\in \B(1)$ for $z \in \B(1-2\epsilon)$,
$$F\circ G(z)=z
\eqno(\arabic{equation})
\newcounter{FG}
\setcounter{FG}{\value{equation}}
\addtocounter{equation}{1}$$
is satisfied on $\B(1-2\epsilon)$, and
$$\left|g_i\right|_{1-2\epsilon,k} \leq C(k)
\cdot\left(1+|f|_{1,k}\right)^{P(k)}|f|_{1,k},
\eqno(\arabic{equation})
\newcounter{gEstimates}
\setcounter{gEstimates}{\value{equation}}
\addtocounter{equation}{1}$$
$$\left|g_i\right|_{1-2\epsilon,k+s} \leq C(k)
\cdot\left(1+|f|_{1,k}\right)^{P(k)}|f|_{1,k+s},$$
with polynomial $P(k)$.
\end{proposition}

\indent
{\bf Proof.}\hspace{0.05in}
Substituting formulas for $F$ and $G$ into (\arabic{FG}) we obtain
the equation
$$g(z)+f(z+g(z))=0
\eqno(\arabic{equation})
\newcounter{fg}
\setcounter{fg}{\value{equation}}
\addtocounter{equation}{1}$$
for mapping $g$. To construct $g$ satisfying (\arabic{fg})
we consider the following sequence of mappings
$g_0,g_1,\dots$ defined inductively
$$\begin{array}{ll}
g_0(z)=0,\vspace{0.1in}\\
g_{l+1}(z)=-f\left(z+g_l(z)\right).
\end{array}
\eqno(\arabic{equation})
\newcounter{gSequence}
\setcounter{gSequence}{\value{equation}}
\addtocounter{equation}{1}$$
\indent
For $g_0$ we have the estimates
$$\begin{array}{ll}
|g_0(z)+f(z+g_0(z))|=|f(z)|\leq \epsilon,\vspace{0.1in}\\
|g_1(z)-g_0(z)| = |-f(z)| \leq \epsilon.
\end{array}$$
\indent
Assuming then the estimates
$$\begin{array}{ll}
|g_l(z)+f(z+g_l(z))| \leq (2n)^l \cdot \epsilon^{l+1},\vspace{0.1in}\\
|g_{l+1}(z)- g_l(z)| \leq (2n)^l \cdot \epsilon^{l+1},
\end{array}
\eqno(\arabic{equation})
\newcounter{gkEstimates}
\setcounter{gkEstimates}{\value{equation}}
\addtocounter{equation}{1}$$
for some $l$ we prove them for $l+1$. For the first estimate using the
mean value theorem we obtain
$$|g_{l+1}(z)+f(z+g_{l+1}(z))|=|-f(z+g_l(z))+f(z-f(z+g_l(z)))|$$
$$=|-f(z+g_l(z))+f(z+g_l(z)-g_l(z)-f(z+g_l(z)))|$$
$$\leq 2n \cdot |f|_{1,1}\cdot |g_l(z)+f(z+g_l(z))|
\leq (2n)^{l+1}\epsilon^{l+2}.$$
\indent
For $g_{l+2}-g_{l+1}$ using the first estimate for $l+1$ we obtain
$$|g_{l+2}(z)-g_{l+1}(z)| = |-f(z+g_{l+1}(z))-g_{l+1}(z)|
\leq (2n)^{l+1}\epsilon^{l+2}.$$
\indent
From estimates (\arabic{gkEstimates}) we obtain that for
$\epsilon < (4n)^{-1}$
$$|g_l(z)| = |g_0(z)+\left(g_1(z)-g_0(z)\right)+\dots
+\left(g_l(z)-g_{l-1}(z)\right)|
\leq \epsilon \cdot\left(1+\frac{1}{2}+\dots+\frac{1}{2^l}\right)
< 2\epsilon,$$
and therefore sequence (\arabic{gSequence}) is well defined for
$z \in \B(0,1-2\epsilon)$ and converges uniformly
on $\B(0,1-2\epsilon)$ to a continuous mapping $g$
satisfying (\arabic{fg}).\\
\indent
To estimate $|g|_k$ and $|g|_{k+s}$ we consider variables
$\left\{x_j\right\}_1^{2n}$ such that $z_j=x_j+\sqrt{-1}x_{n+j}$
and functions $\left\{y_j\right\}_1^{2n}$ and $\left\{u_j\right\}_1^{2n}$
such that $f_j(x)=y_j(x)+\sqrt{-1}y_{n+j}(x)$ and
$g_j(x)=u_j(x)+\sqrt{-1}u_{n+j}(x)$.
Then we differentiate (\arabic{fg}) with respect to variables $x_j$
and as in (\cite{St}, Appendix 1) obtain the following system of
differential equations
$$\left[I + \frac{\partial y}{\partial x}\left(x+u(x)\right)\right]
\cdot\frac{\partial u}{\partial x}(x)
= -\frac{\partial y}{\partial x}(x),$$
or
$$\frac{\partial u}{\partial x}(x)
= -\left[I + \frac{\partial y}{\partial x}\left(x+u(x)\right)\right]^{-1}
\frac{\partial y}{\partial x}(x).
\eqno(\arabic{equation})
\newcounter{uEquation}
\setcounter{uEquation}{\value{equation}}
\addtocounter{equation}{1}$$
\indent
Successively differentiating (\arabic{uEquation}) we obtain
(\arabic{gEstimates}).\qed\\
\indent
Using proposition~\ref{Inverse} we can assume now the existence of
finite covers
$\left\{{\cal V}^{\iota} \subset {\cal U}^{\iota}\right\}_1^N$
of some neighborhood of ${\bold M}_0$ in ${\bold G}$ and
of $C^p$-diffeomorphisms
${\cal F}:{\bold M}_0\to {\bold M}$ and
${\cal G} = {\cal F}^{-1}:{\bold M}\to {\bold M}_0$ such that
${\cal F}(z) = z+f^{\iota}(z)$ in $U^{\iota}_0
={\cal U}^{\iota}\cap {\bold M}_0$
with $f^{\iota} \in \left[C^p(U^{\iota}_0)\right]^n$
and ${\cal G}(z) = z+g^{\iota}(z)$ in $V^{\iota}
={\cal V}^{\iota}\cap {\bold M}$
with $g^{\iota} \in \left[C^p\left(V^{\iota}\right)\right]^n$.\\
\indent
In what follows we will use special coordinates in
$\left\{{\cal U}^{\iota}\right\}$.
Namely, for a fixed neighborhood ${\cal U}$
we choose coordinates $\{Z_j=X_j+iY_j\}_1^m$ and
$\{W_j=U_j+iV_j\}_1^{n-m}$ so that
$$\begin{array}{ll}
\rho^{(0)}_l(z) = X_l - \phi^{(0)}_l\left(Y_1,\dots,Y_m,W_1,
\dots,W_{n-m}\right)\hspace{0.1in}\mbox{for}\hspace{0.1in}l=1,\dots,m,
\vspace{0.1in}\\
\rho_l(z) = \rho^{(0)}_l\left(z+g(z)\right)
\end{array}
\eqno(\arabic{equation})
\newcounter{DefiningFunctions}
\setcounter{DefiningFunctions}{\value{equation}}
\addtocounter{equation}{1}$$
are defining functions for ${\bold M}_0$ and ${\bold M}$ respectively.\\
\indent
For $\delta$ small enough according to implicit function theorem
manifold ${\bold M}$ with $|g|_p < \delta$ can be defined as
$${\bold M} = \left\{z \in {\cal U}:X_l - \phi_l\left(Y_1,\dots,Y_m,W_1,
\dots,W_{n-m}\right)=0\hspace{0.1in}\mbox{for}\hspace{0.1in}l=1,
\dots,m\right\}.$$
\indent
In the lemma below we prove necessary estimates for
$\sup_l\left|\phi_l - \phi^{(0)}_l\right|_p$ in terms of $|g|_p$
for $|g|_p$ small enough.

\begin{lemma}\label{phi-phi0}(cf. \cite{W2})
Let functions $\rho^{(0)}_l$ and $\rho_l$ be defined in ${\cal U}$
by formulas (\arabic{DefiningFunctions}). Then there exists $\delta > 0$
such that for $g$ with $|g|_p < \delta$ and $k \leq p-1$ the following
estimate holds:
$$\left|\phi_l-\phi^{(0)}_l\right|_k\leq C(k)
\left(1+\left|\rho^{(0)}\right|_{k+1}\right)^{P(k)}\cdot|g|_k.
\eqno(\arabic{equation})
\newcounter{phi0-phiestimate}
\setcounter{phi0-phiestimate}{\value{equation}}
\addtocounter{equation}{1}$$
\end{lemma}

\indent
{\bf Proof.}\hspace{0.05in}
We consider for $z \in {\bold M}$ the following equality
$$\phi^{(0)}_l\left(Y_1,\dots,Y_m,W_1,\dots,W_{n-m}\right)
-\phi_l\left(Y_1,\dots,Y_m,W_1,\dots,W_{n-m}\right)
\eqno(\arabic{equation})
\newcounter{phi0-phi}
\setcounter{phi0-phi}{\value{equation}}
\addtocounter{equation}{1}$$
$$= \mbox{Re}g_l(z)+\phi^{(0)}_l\left(Y_1,\dots,Y_m,W_1,\dots,W_{n-m}
\right)$$
$$-\phi_l^{(0)}\left(Y_1+\mbox{Im}g_1(z),\dots,Y_m+\mbox{Im}g_m(z),
W_1+g_{m+1}(z),\dots,W_{n-m}+g_n(z)\right),$$
which is a corollary of equalities
$$X_l - \phi_l\left(Y_1,\dots,Y_m,W_1,\dots,W_{n-m}\right)=0$$
and
$$\rho^{(0)}_l\left(z+g(z)\right)$$
$$=X_l + \mbox{Re}g_l(z)
-\phi^{(0)}_l\left(Y_1+\mbox{Im}g_1(z),\dots,Y_m+\mbox{Im}g_m(z),
W_1+g_{m+1}(z),\dots,W_{n-m}+g_n(z)\right)=0.$$
\indent
From equality (\arabic{phi0-phi}) we obtain the estimate
$$\left|\phi_l - \phi^{(0)}_l\right|_0 \leq C\left(|g|_0
+\left|\rho^{(0)}\right|_1|g|_0\right).
\eqno(\arabic{equation})
\newcounter{phi0-phiestimate0}
\setcounter{phi0-phiestimate0}{\value{equation}}
\addtocounter{equation}{1}$$
\indent
To estimate the derivatives we differentiate with respect
to $Y,W$ equalities
$$\rho^{(0)}_l\left(z+g(z)\right)$$
$$=\rho^{(0)}_l\left(\phi_1(Y,W)+\mbox{Re}g_1\left(\phi(Y,W),Y,W\right),
\dots,\phi_m(Y,W)+\mbox{Re}g_m\left(\phi(Y,W),Y,W\right),\right.$$
$$Y_1+\mbox{Im}g_1\left(\phi(Y,W),Y,W\right),\dots,
Y_m+\mbox{Im}g_m\left(\phi(Y,W),Y,W\right),$$
$$\left.W_1+g_{m+1}\left(\phi(Y,W),Y,W\right),\dots,
W_{n-m}+g_n\left(\phi(Y,W),Y,W\right)\right)=0$$
for $l=1,\dots,m$. Then, using condition
$$\left[\frac{\partial\rho^{(0)}}{\partial X}\right]=I,$$
we obtain equality in a matrix form
$$\left[\frac{\partial\phi}{\partial\left(Y,W\right)}\right]
+\left[\frac{\partial\mbox{Re}g^{\prime}}{\partial X}\right]
\left[\frac{\partial\phi}{\partial\left(Y,W\right)}\right]
+\left[\frac{\partial\mbox{Re}g^{\prime}}{\partial\left(Y,W\right)}\right]$$
$$+\left[\frac{\partial\rho^{(0)}}{\partial\left(Y,W\right)}
\left(z+g(z)\right)\right]\times
\left(I+\left[\frac{\partial\left(\mbox{Im}g^{\prime},
g^{\prime\prime}\right)}{\partial\left(Y,W\right)}\right]
+\left[\frac{\partial\left(\mbox{Im}g^{\prime},g^{\prime\prime}\right)}
{\partial X}\right]
\left[\frac{\partial\phi}{\partial\left(Y,W\right)}\right]\right)=0,$$
where $g^{\prime}=\left(g_1,\dots,g_m\right)$ and
$g^{\prime\prime}=\left(g_{m+1},\dots,g_n\right)$.
We transform this equality into
$$\left(I+\left[\frac{\partial\mbox{Re}g^{\prime}}{\partial X}\right]
+\left[\frac{\partial\rho^{(0)}}{\partial\left(Y,W\right)}(z+g(z))\right]
\times
\left[\frac{\partial\left(\mbox{Im}g^{\prime},g^{\prime\prime}\right)}
{\partial X}\right]\right)
\times\left[\frac{\partial\phi}{\partial\left(Y,W\right)}\right]$$
$$=-\left[\frac{\partial\mbox{Re}g^{\prime}}{\partial\left(Y,W\right)}
\right]-\left[\frac{\partial\rho^{(0)}}{\partial\left(Y,W\right)}(z+g(z))
\right]\times
\left(I+\left[\frac{\partial\left(\mbox{Im}g^{\prime},
g^{\prime\prime}\right)}{\partial\left(Y,W\right)}\right]\right),$$
and then into
$$\left[\frac{\partial\phi}{\partial\left(Y,W\right)}\right]
=\left(I+\left[\frac{\partial\mbox{Re}g^{\prime}}{\partial X}\right]
+\left[\frac{\partial\rho^{(0)}}{\partial\left(Y,W\right)}(z+g(z))
\right]\times
\left[\frac{\partial\left(\mbox{Im}g^{\prime},g^{\prime\prime}\right)}
{\partial X}\right]\right)^{-1}$$
$$\times
\left(-\left[\frac{\partial\mbox{Re}g^{\prime}}{\partial\left(Y,W\right)}
\right]-
\left[\frac{\partial\rho^{(0)}}{\partial\left(Y,W\right)}(z+g(z))\right]
\times
\left(I+\left[\frac{\partial\left(\mbox{Im}g^{\prime},
g^{\prime\prime}\right)}{\partial\left(Y,W\right)}\right]\right)\right).$$
\indent
Using equality
$$\frac{\partial}{\partial\left(Y,W\right)}\rho^{(0)}(z+g(z))
=\frac{\partial}{\partial\left(Y,W\right)}\left[\rho^{(0)}(z)
+\int_0^1\Bigl\langle\nabla\rho^{(0)}\left(z+tg(z)\right),
g(z)\Bigr\rangle dt\right]$$
we can rewrite the previous equality for $|g|_{k+1}$ small enough as
$$\left[\frac{\partial\phi}{\partial\left(Y,W\right)}\right]
=\left(I+A\right)\left(
-\left[\frac{\partial\rho^{(0)}}{\partial\left(Y,W\right)}(z)\right]
+B\right),
\eqno(\arabic{equation})
\newcounter{dphidYW}
\setcounter{dphidYW}{\value{equation}}
\addtocounter{equation}{1}$$
with matrices $A, B$ satisfying estimates
$$\left|A\right|_k, \left|B\right|_k\leq C(k)
\left(1+\left|\rho^{(0)}\right|_{k+2}\right)^{P(k)}|g|_{k+1}.
\eqno(\arabic{equation})
\newcounter{matrixestimates}
\setcounter{matrixestimates}{\value{equation}}
\addtocounter{equation}{1}$$
\indent
For $g \equiv 0$ (\arabic{dphidYW}) becomes
$$\left[\frac{\partial\phi^{(0)}}{\partial\left(Y,W\right)}\right]
=-\left[\frac{\partial\rho^{(0)}}{\partial\left(Y,W\right)}(z)\right],
\eqno(\arabic{equation})
\newcounter{dphi0}
\setcounter{dphi0}{\value{equation}}
\addtocounter{equation}{1}$$
and therefore we can rewrite (\arabic{dphidYW}) as
$$\left[\frac{\partial\phi^{(0)}}{\partial\left(Y,W\right)}\right]
-\left[\frac{\partial\phi}{\partial\left(Y,W\right)}\right]=E,$$
with a matrix $E$ satisfying estimate (\arabic{matrixestimates}).
Using this representation and estimates (\arabic{matrixestimates})
we obtain
$$\sup_l\left|\phi_l-\phi^{(0)}_l\right|_k
=\left|\left[\frac{\partial\phi^{(0)}}{\partial\left(Y,W\right)}\right]
-\left[\frac{\partial\phi}{\partial\left(Y,W\right)}\right]\right|_{k-1}$$
$$\leq C(k)\left(1+\left|\rho^{(0)}\right|_{k+1}\right)^{P(k)}|g|_{k}.$$
\qed

We assume from now on that the neighborhoods ${\cal U}^{\iota}$ and functions
$\phi^{(0)}_l$ are chosen so that conditions (\arabic{DefiningFunctions})
are satisfied. We assume also that for small enough $\delta > 0$
and a manifold ${\bold M}$ with $|g|_p < \delta$ condition
(\arabic{NonDegeneracy}) is satisfied on ${\cal U}^{\iota}$ for some fixed
$c > 0$.\\
\indent
As before our local extension operators $\left\{E^{\iota}\right\}$
for functions are defined by the formula
$$E^{\iota}(h)(z) = h\left(Y_1,\dots,Y_m,W_1,\dots,W_{n-m}\right),$$
and for a differential form
$$h = \sum_{I,J,K} h_{I,J,K}dY_I \wedge dW_J \wedge d{\overline W}_K$$
by extending coefficients as in the formula above.\\
\indent
A global extension operator for any ${\bold M}$ close
to ${\bold M}_0$ we define by the formula
$${\bold E}_{\bold M}(h)
= \sum_{\iota} E^{\iota}\left(\vartheta_{\iota}h\right),$$
where $h$ is a differential form on ${\bold M}$,
$\left\{\vartheta_{\iota}\right\}$ is a partition of unity
subordinate to the cover $\left\{ {\cal V}^{\iota} \right\}
\subset \left\{ {\cal U}^{\iota} \right\}$
and $\left\{E_{\iota}\right\}$ are local extension operators.\\

\indent
Our next goal is to estimate the operator
$$\left({\bold E}_{{\bold M}_0}\circ{\bold H}^r_{{\bold M}_0}
\circ{\bold E}_{\bold M}-{\bold H}^r_{\bold M}\right),$$
where ${\bold E}_{\bold M}$ and ${\bold E}_{{\bold M}_0}$
are extension operators from ${\bold M}$ and ${\bold M}_0$
respectively, $r<q$, and manifold ${\bold M}$ is close enough to
${\bold M}_0$. In the lemma below we prove a special representation
for operator ${\bold H}^r_{\bold M}$.\\

\begin{lemma}\label{HFinite}
For any $\epsilon > 0$ there exist $\delta,\alpha > 0$ such that
for ${\bold M}$ with $|{\cal E}|_p < \delta$
the following representation holds
$${\bold H}^r_{\bold M} = {\bold N}^r\circ{\bold E}_{\bold M}
+ {\bold L}^r_{\bold M},$$
where ${\bold E}_{\bold M}$ is the extension operator,
$${\bold N}^r:C^k_{(0,r)}\left({\bold G}(\alpha),
{\cal B}\right) \rightarrow
C^p_{(0,r)}\left({\bold M},{\cal B}\big|_{\bold M}\right)$$
is a finite-dimensional operator,
$${\bold G}(\alpha)= \bigcup_{\iota}\left\{z \in {\cal U}^{\iota}:
\rho^{(0)}_{\iota}(z) < \alpha\right\},$$
and ${\bold L}^r_{\bold M}$ admits an estimate
$$\left|{\bold L}^r_{\bold M}\right|_k \leq C(k)
\epsilon \left(1 + |\rho|_{k+4}\right)^{P(k)}
\eqno(\arabic{equation})
\newcounter{LMEstimate}
\setcounter{LMEstimate}{\value{equation}}
\addtocounter{equation}{1}$$
for $k \leq p-4$.
\end{lemma}

\indent
{\bf Proof.}\hspace{0.05in}
Using Lemma~\ref{HZero} in formula (\arabic{boldH}) we obtain
$${\bold H}^r_{\bold M}(h)(z) = \sum_{\iota} \left[ - \bar\partial_{\bold M}
\vartheta^{\prime}_{\iota}(z) \wedge R_r^{\iota}(\vartheta_{\iota}h)(z)
+ \vartheta^{\prime}_{\iota}(z)
R_{r+1}^{\iota}(\bar\partial_{\bold M} \vartheta_{\iota} \wedge h)(z)
\right]$$
and conclude that the statement of the lemma would follow from a
finite-dimensional approximation of operators $R_r$ with a remainder
admitting estimate (\arabic{LMEstimate}).
We fix a neighborhood ${\cal U}$ from the cover
$\left\{{\cal U}^{\iota}\right\}_1^N$ and consider an approximation of
the local solution operator $R_r$ by the operators
$$R^{(0)}_r(\alpha)(h)(z) = (-1)^{r} \frac{(n-1)!}{(2\pi i)^n}$$
$$\times\int_{U_0(\alpha)\times [0,1]}\vartheta (\zeta)E(h)(\zeta)
\wedge\omega^{\prime}_{r-1}\left((1-t)\frac{\bar\zeta - \bar z}
{{\mid \zeta - z \mid}^2} + t\frac{P(\zeta,z)}
{\Phi(\zeta,z)}\right) \wedge\omega(\zeta).$$
\indent
From estimate (\arabic{phi0-phiestimate}) we obtain
$$|\rho_l-\rho^{(0)}_l|_k \leq C(k)
\left(1 + |\rho^{(0)}_l|_{k+1}\right)|g|_k
\eqno(\arabic{equation})
\newcounter{rhoMestimate}
\setcounter{rhoMestimate}{\value{equation}}
\addtocounter{equation}{1}$$
for $l=1,\dots,m$. Therefore, for any $\alpha > 0$ there exists
$\delta(\alpha) > 0$ such that for $|g|_p < \delta(\alpha)$ the
following inclusions hold
$${\cal U}(\alpha/2) \subset {\cal U}_0(\alpha) \subset
{\cal U}(2\alpha).$$
Then
$$\omega^{\prime}_{r-1}\left((1-t)\frac{\bar\zeta - \bar z}
{{\mid \zeta - z \mid}^2} + t\frac{P(\zeta,z)}
{\Phi(\zeta,z)}\right) \wedge\omega(\zeta)
\in C^p_{\left\{(0,r-1),(n,n-r-1)\right\}}
\left({\cal U}_z(\alpha/2)\times U_{0,\zeta}(\alpha)\times [0,1]\right)$$
for such $g$ and therefore for any $\epsilon > 0$ we can find differential
forms
$$h_j \in C^p_{(0,r-1)}\left({\cal U}_z(\alpha/2)\right),
\hspace{0.05in}
N_j \in C^p_{(n,n-r-1)}\left(U_{0,\zeta}(\alpha)\times [0,1]\right),$$
such that
$$\left|\omega^{\prime}_{r-1}\left((1-t)\frac{\bar\zeta - \bar z}
{{\mid \zeta - z \mid}^2} + t\frac{P(\zeta,z)}
{\Phi(\zeta,z)}\right) \wedge\omega(\zeta)
- \sum_{j=1}^{l} N_j(\zeta,t)h_j(z)
\right|_{C^p\left({\cal U}_z(\alpha/2)\times
U_{0,\zeta}(\alpha)\times [0,1]\right)} < \epsilon.$$
\indent
Therefore, for any $\alpha > 0$ there exists a finite-dimensional operator
$$N^r(u)(z) = \sum_{j=1}^l h_j(z)\int_{U_0(\alpha)\times [0,1]}
u(\zeta)N_j(\zeta,t),$$
such that for any
$u \in C^k_{(0,r-1)}\left({\cal U}_0(\alpha)\right)$
with $\left|u \right|_k < 1$ we have
$$\left|R^{(0)}_r(\alpha)(u) - N^r(u) \right|_k \leq C\epsilon.
\eqno(\arabic{equation})
\newcounter{R-Nestimate}
\setcounter{R-Nestimate}{\value{equation}}
\addtocounter{equation}{1}$$
\indent
To estimate the $C^k$ norm of the operator $R^{(0)}_r(\alpha) - R_r$ we use
representation (\arabic{cauchymartinellikernel}) and
Lemmas~\ref{Smoothness} and \ref{ReductionTok0} and reduce the problem to
the $C^0$-estimate of integrals
$$\int_{U_0(\alpha)\times [0,1]} E(h)(\zeta)c_{\bold M}(\zeta, z, t)
{\cal K}^{S}_{{\bold M},a,b}(\zeta, z)dt
- \int_{U(\beta)\times [0,1]} E(h)(\zeta)c_{\bold M}(\zeta, z, t)
{\cal K}^{S}_{{\bold M},a,b}(\zeta, z) dt
\eqno(\arabic{equation})
\newcounter{Mdifference}
\setcounter{Mdifference}{\value{equation}}
\addtocounter{equation}{1}$$
with $\|c(\zeta, z, t)\|_1 \leq C(k)\left(1 + |\rho|_{k+4}\right)^{P(k)}$
and $\beta < \alpha/2$.\\
\indent
As in the proof of Proposition~\ref{AlmostHomotopy} we
use Lemma~\ref{RDoomed} to conclude that it suffices to consider
only integrals with kernels ${\cal K}^{S}_{{\bold M},a,b}$
satisfying (\arabic{FinalIndices}), since otherwise estimate
(\arabic{DoomedEstimate}) holds, and therefore for $\alpha$
small enough the estimate (\arabic{LMEstimate}) will also hold.\\
\indent
To estimate (\arabic{Mdifference}) we consider a function
$\phi(z,\zeta) \in C^{\infty}({\cal U}_{z}\times{\cal U}_{\zeta})$,
such that in coordinates
$$\begin{array}{ll}
z_j=X_j+\sqrt{-1}Y_j\hspace{0.05in}(j=1,\dots,m),
z_j=W_{j-m}\hspace{0.05in}(j=m+1,\dots,n),\vspace{0.1in}\\
\zeta_j=\xi_j+\sqrt{-1}\eta_j\hspace{0.05in}(j=1,\dots,m),
\zeta_j=\omega_{j-m}\hspace{0.05in}(j=m+1,\dots,n),
\end{array}$$
we have
$$\phi(z,\zeta)=\phi\left(Y_1,\dots,Y_m,W_1,\dots,W_{n-m},
\eta_1,\dots,\eta_m,\omega_1,\dots,\omega_{n-m}\right)$$
and
$$\begin{array}{lll}
\phi(z,\zeta)\equiv 1\hspace{0.1in}\mbox{for}\hspace{0.1in}
\zeta \in {\cal W}(z,\alpha,\sqrt{\alpha}),\vspace{0.1in}\\
\phi(z,\zeta)\equiv 0\hspace{0.1in}\mbox{for}\hspace{0.1in}
\zeta\notin{\cal W}(z,C\alpha,\sqrt{C\alpha})
\hspace{0.1in}\mbox{with some}\hspace{0.1in}C >1.
\end{array}$$
\indent
Then we represent for $\beta<\alpha/2$
$$\int_{U(\beta)\times [0,1]} E(h)(\zeta)c_{\bold M}(\zeta,z,t)
{\cal K}^{S}_{{\bold M},a,b}(\zeta,z)dt
\eqno(\arabic{equation})
\newcounter{intRepresentation}
\setcounter{intRepresentation}{\value{equation}}
\addtocounter{equation}{1}$$
$$=\int_{U(\beta)\times [0,1]}E(h)(\zeta)\phi(z,\zeta)
c_{\bold M}(\zeta,z,t){\cal K}^{S}_{{\bold M},a,b}(\zeta,z)dt$$
$$+\int_{U(\beta)\times[0,1]}E(h)(\zeta)\left(1-\phi(z,\zeta)\right)
c_{\bold M}(\zeta,z,t){\cal K}^{S}_{{\bold M},a,b}(\zeta,z)dt.$$
\indent
For the first term of the right hand side of
(\arabic{intRepresentation}) we obtain as in (\arabic{WEstimate})
$$\left|\int_{U(\beta)\times [0,1]}E(h)(\zeta)\phi(z,\zeta)
c_{\bold M}(\zeta,z,t){\cal K}^{S}_{{\bold M},a,b}(\zeta,z)dt\right|_0$$
$$\leq\int_{U(\beta) \cap {\cal W}(z,C\alpha,\sqrt{C\alpha})\times[0,1]}
|h(\zeta,z,t)||c_{\bold M}(\zeta,z,t)|
|{\cal K}^{S}_{{\bold M},a,b}(\zeta,z)|dt$$
$$\leq C \|c_{\bold M}\|_0|h|_0 \alpha^{l({\cal K})}\cdot
{\cal I}_1\left\{k({\cal K}), h({\cal K})\right\}
\left( \alpha,\sqrt{\alpha}\right)
\leq C(k)\alpha^{1/2}\left(1 + |\rho^{(0)}|_{k+4}\right)^{P(k)}|h|_0.$$
\indent
For the second term of the right hand side of (\arabic{intRepresentation})
using the fact that by construction form $E(h)$ and function
$\phi(z,\zeta)$ do not depend on variables $\xi_1,\dots,\xi_m$ and that
form $E(h){\cal K}^{S}_{{\bold M},a,b}$ contains differentials 
$d\eta_1,\dots,d\eta_m,d\omega_1,\dots,d\omega_{n-m},
d\bar\omega_1,\dots,d\bar\omega_{n-m}$ we apply Stokes' formula and obtain
$$\int_{U_0(\alpha)\times [0,1]}E(h)(\zeta)\left(1-\phi(z,\zeta)\right)
c_{\bold M}(\zeta,z,t){\cal K}^{S}_{{\bold M},a,b}(\zeta, z)dt
\eqno(\arabic{equation})
\newcounter{intRepresentation1}
\setcounter{intRepresentation1}{\value{equation}}
\addtocounter{equation}{1}$$
$$- \int_{U(\beta)\times [0,1]}E(h)(\zeta)\left(1-\phi(z,\zeta)\right)
c_{\bold M}(\zeta,z,t){\cal K}^{S}_{{\bold M},a,b}(\zeta, z)dt$$
$$=\int_{\left[{\cal U}_0(\alpha)\setminus{\cal U}(\beta)\right]\times [0,1]}
E(h)(\zeta)\left(1-\phi(z,\zeta)\right)d\left[c_{\bold M}(\zeta,z,t)
{\cal K}^{S}_{{\bold M},a,b}(\zeta, z)dt\right]$$
$$=\int_{\left[{\cal U}_0(\alpha)\setminus{\cal U}(\beta)\right]\times [0,1]}
E(h)(\zeta)\left(1-\phi(z,\zeta)\right)d\left[c_{\bold M}(\zeta,z,t)\right]
\wedge{\cal K}^{S}_{{\bold M},a,b}(\zeta, z)dt$$
$$+\int_{\left[{\cal U}_0(\alpha)\setminus{\cal U}(\beta)\right]\times [0,1]}
E(h)(\zeta)\left(1-\phi(z,\zeta)\right)c_{\bold M}(\zeta,z,t)
d\left[{\cal K}^{S}_{{\bold M},a,b}(\zeta, z)\right]dt.$$
\indent
For the first term of the right hand side of (\arabic{intRepresentation1})
using estimates for $c_{\bold M}$ and integrating in coordinates
$$\rho=\sqrt{\sum_1^m \rho_j^2},\ \theta_j = \frac{\rho_j}{\rho}
\ (j=1,\dots,m),\ \mbox{Im}F_j,\ (j=1,\dots,m),\ \omega_l\ (l=1,\dots,n-m),$$
we have
$$\left|\int_{\left[{\cal U}_0(\alpha)\setminus{\cal U}(\beta)\right]
\times [0,1]}
E(h)(\zeta)\left(1-\phi(z,\zeta)\right)d\left[c_{\bold M}(\zeta,z,t)\right]
\wedge{\cal K}^{S}_{{\bold M},a,b}(\zeta, z)dt\right|$$
$$\leq\int_{\left( {\cal U}(2\alpha)\setminus \left[{\cal U}(\beta)
\cup {\cal W}(z,\alpha,\sqrt{\alpha})\right] \right) \times [0,1]}
|E(h)(\zeta)|\left|\partial_{\rho}\left[c_{\bold M}(\zeta,z,t)\right]
\wedge{\cal K}^{S}_{{\bold M},a,b}(\zeta, z)\right|dt$$
$$\leq C(k)\cdot\left(1 + |\rho^{(0)}|_{k+4}\right)^{P(k)}\cdot |h|_0
\left[ {\cal I}_2\left\{ k({\cal K}),b-l({\cal K}) \right\}
\left( \alpha, \sqrt{\alpha} \right)\right]\int_0^{2\alpha} d\rho$$
$$\leq C(k)\alpha^{1/2}\cdot\left(1 + |\rho^{(0)}|_{k+4}\right)^{P(k)}
\cdot|h|_0,$$
where
$$\partial_{\rho}f=d\rho\wedge\frac{\partial f}{\partial\rho},$$
and in the last inequality we used Lemma~\ref{Integral} for
kernel ${\cal K}^{S}_{{\bold M},a,b}$, satisfying conditions
(\arabic{FinalIndices}).\\
\indent
Again using estimates for $c_{\bold M}$ we obtain for the second term
of the right hand side of (\arabic{intRepresentation1})
$$\left|\int_{\left[{\cal U}_0(\alpha)\setminus{\cal U}(\beta)\right]
\times [0,1]}
E(h)(\zeta)\left(1-\phi(z,\zeta)\right)c_{\bold M}(\zeta,z,t)
d\left[{\cal K}^{S}_{{\bold M},a,b}(\zeta, z)\right]dt\right|$$
$$\leq\int_{\left( {\cal U}(2\alpha)\setminus \left[{\cal U}(\beta)
\cup {\cal W}(z,\alpha,\sqrt{\alpha})\right] \right) \times [0,1]}
|E(h)(\zeta)|\left|c_{\bold M}(\zeta,z,t)
d\left[{\cal K}^{S}_{{\bold M},a,b}(\zeta, z)\right]\right|dt$$
$$\leq C(k)\cdot\left(1 + |\rho^{(0)}|_{k+4}\right)^{P(k)}\cdot |h|_0$$
$$\times\left[ {\cal I}_2\left\{ k({\cal K})+1,b-l({\cal K}) \right\}
\left( \alpha, \sqrt{\alpha} \right)
+ {\cal I}_2\left\{ k({\cal K}),b-l({\cal K})+1 \right\}
\left( \alpha,\sqrt{\alpha} \right)\right]\int_0^{2\alpha}d\rho$$
$$\leq C(k)\cdot\left(1 + |\rho^{(0)}|_{k+4}\right)^{P(k)}\alpha^{1/2}
\cdot |h|_0,$$
where the indices in integrals ${\cal I}_2$ are the
indices that appear after differentiation of the kernel
${\cal K}^{S}_{{\bold M},a,b}$, namely, either $k({\cal K})$
or $b$ increase by one.\\
\indent
In the last inequality we used inequalities
$$\begin{array}{lll}
k({\cal K})+1+b-l({\cal K}) \leq 2n-m,\vspace{0.1in}\\
k({\cal K})+1+2\left(b-l({\cal K})\right) \leq 2n-m+2,\vspace{0.1in}\\
k({\cal K})+2\left(b-l({\cal K})+1\right) \leq 2n-m+2,
\end{array}$$
which follow from inequalities (\arabic{FinalIndices}) for
${\cal K}^{S}_{{\bold M},a,b}(\zeta, z)$ and Lemma~\ref{Integral}.\\
\indent
Combining the estimates above we obtain that for any fixed $\alpha >0$
there exists $\delta(\alpha) >0$ such that for any submanifold ${\bold M}$
with $|{\cal E}|_p<\delta(\alpha)$ we have
$$\left|R^{(0)}_r(\alpha) - R_r\right|_k
\leq C(k)\cdot\left(1 + |\rho^{(0)}|_{k+4}\right)^{P(k)}\alpha^{1/2},$$
and combining this estimate with estimate (\arabic{R-Nestimate})
we obtain the statement of the Lemma.
\qed

\indent
Using Lemmas~\ref{phi-phi0} and ~\ref{HFinite} we obtain the following

\begin{proposition}\label{HM-HM0}
Let $\rho^{(0)}_l, \rho_l$ and $g$ satisfy conditions
(\arabic{DefiningFunctions}).Then for any $\epsilon >0$ there exists
$\delta > 0$ such that for an arbitrary compact, regular q-pseudoconcave CR
submanifold ${\bold M}$ with $\left|{\cal E}\right|_p < \delta$
and $k \leq p-4$ the following estimate holds
$$\left|{\bold E}_{{\bold M}_0}\circ{\bold H}_{{\bold M}_0}
\circ{\bold E}_{\bold M}-{\bold H}^r_{\bold M}\right|_k \leq C(k)\epsilon
\left(1+ |\rho^{(0)}|_{k+4}\right)^{P(k)}.
\eqno(\arabic{equation})
\newcounter{HM-HM0estimate}
\setcounter{HM-HM0estimate}{\value{equation}}
\addtocounter{equation}{1}$$
\end{proposition}
\vspace{0.2in}

\indent
{\bf Proof.}\hspace{0.05in}
Using Lemma~\ref{HFinite} we obtain equality
$${\bold E}_{{\bold M}_0}\circ{\bold H}_{{\bold M}_0}
\circ{\bold E}_{\bold M}(h)-{\bold H}^r_{\bold M}(h)
\eqno(\arabic{equation})
\newcounter{HM-HM0equality}
\setcounter{HM-HM0equality}{\value{equation}}
\addtocounter{equation}{1}$$
$$={\bold E}_{{\bold M}_0}\circ\left({\bold N}^r\circ
{\bold E}_{{\bold M}_0}+{\bold L}^r_{{\bold M}_0}\right)
\circ{\bold E}_{\bold M}(h)-{\bold N}^r\circ{\bold E}_{\bold M}(h)
-{\bold L}^r_{\bold M}(h)$$
$$=\left({\bold E}_{{\bold M}_0}\circ{\bold N}^r
\circ{\bold E}_{{\bold M}_0}
\circ{\bold E}_{\bold M}-{\bold N}^r\circ{\bold E}_{\bold M}\right)(h)
+\left({\bold E}_{{\bold M}_0}\circ{\bold L}^r_{{\bold M}_0}
\circ{\bold E}_{\bold M}-{\bold L}_{\bold M}\right)(h)$$
$$=\left({\bold E}_{{\bold M}_0}-{\bold I}\right)
\circ{\bold N}^r\circ{\bold E}_{{\bold M}_0}
\circ{\bold E}_{\bold M}(h)
+{\bold N}^r\circ\left({\bold E}_{{\bold M}_0}-{\bold I}\right)
\circ{\bold E}_{\bold M}(h)$$
$$+\left({\bold E}_{{\bold M}_0}\circ{\bold L}^r_{{\bold M}_0}
\circ{\bold E}_{\bold M}-{\bold L}_{\bold M}\right)(h).$$
\indent
For the first term of the right hand
side of (\arabic{HM-HM0equality}) we reduce the estimate to the
estimate of the operator ${\bold E}_{{\bold M}_0}-{\bold I}$ and
then further reduce it to the estimate of $E^{\iota}_{{\bold M}_0}-I$
for a local extension operator $E^{\iota}_{{\bold M}_0}$ on
$C^{k}_{(0,0)}\left({\cal U}^{\iota},{\cal B}\right).$\\
\indent
For any $f \in C^{p}_{(0,0)}\left({\cal U}^{\iota},{\cal B}\right)$
using Lemma~\ref{phi-phi0} we obtain
$$\left|E^{\iota}_{{\bold M}_0}(f)\big|_{\bold M}
-f\big|_{\bold M}\right|_k=\left|f\left(\phi^{(0)}(Y,W),Y,W\right)
-f\left(\phi(Y,W),Y,W\right)\right|_k
\eqno(\arabic{equation})
\newcounter{E-Iestimate}
\setcounter{E-Iestimate}{\value{equation}}
\addtocounter{equation}{1}$$
$$\leq C(k)|g^{\iota}|_k\left(1+ |\rho^{(0)}|_{k+1}\right)^{P(k)}
|f|_{k+1}.$$
Applying this estimate to
$$f(z)=N^{\iota}(E(h))(z)=\sum_{j=1}^l h_j(z)\int_{U_0(\alpha)\times [0,1]}
E(h)(\zeta)N_j(\zeta,t)$$
and using estimate
$$\left|N^{\iota}(E(h))\right|_k \leq C(k)|h|_0
\eqno(\arabic{equation})
\newcounter{NEestimate}
\setcounter{NEestimate}{\value{equation}}
\addtocounter{equation}{1}$$
we obtain the estimate
$$\left|E^{\iota}_{{\bold M}_0}(f)\big|_{\bold M}-f\big|_{\bold M}\right|_k
\leq C(k)|g^{\iota}|_k\left(1+ |\rho^{(0)}|_{k+1}\right)^{P(k)}|h|_0.$$
\indent
For the second term of the right hand side of (\arabic{HM-HM0equality}) we
apply estimates (\arabic{E-Iestimate}) and (\arabic{NEestimate}) in opposite
order and obtain for $u={\bold E}_{\bold M}(h)$
$$\left|N^{\iota}\left(E^{\iota}(u)-u\right)\right|_k
\leq C(k)\left|E^{\iota}(u)-u\right|_0$$
$$\leq C(k)|g^{\iota}|_0\left(1+ |\rho^{(0)}|_{1}\right)^{P}|h|_1.$$
\indent
Necessary estimate for the third term of the right hand side of
(\arabic{HM-HM0equality}) follows from estimate (\arabic{LMEstimate})
in Lemma~\ref{HFinite}.
\qed

\indent
{\bf Proof of Theorem~\ref{HomotopyTheorem}.}\hspace{0.1in}
We consider a compact, regular q-pseudoconcave CR submanifold ${\bold M}_0$
such that $\dim H^r\left({\bold M}_0,{\cal B}\big|_{{\bold M}_0}\right)=0$
for some fixed $r<q$.
We assume that a compact, regular q-pseudoconcave CR submanifold ${\bold M}$
is close to ${\bold M}_0$ with $\left|{\cal E}\right|_p$ small enough
so that functions $\rho^{(0)}_l, \rho_l$ and $g$ satisfy
(\arabic{DefiningFunctions}).\\
\indent
Abusing notation we will denote also by
$\left\{g_j\right\}_1^s \in C^{p-4}_{(0,r-1)}\left({\bold G},{\cal B}\right)$
and $\left\{f_j\right\}_1^l \in C^{p-4}_{(0,r)}\left({\bold G},
{\cal B}\right)$ extensions of the corresponding forms from
Proposition~\ref{M0Formula}.\\
\indent
We define for $h \in C^k_{(0,r)}\left({\bold M},
{\cal B}\big|_{\bold M}\right)$, $v \in C^k_{(0,r+1)}\left({\bold M},
{\cal B}\big|_{\bold M}\right)$
$${\bold P}^r_{\bold M}(h)={\bold R}^r_{\bold M}
\left(h-\sum_{i=1}^l \beta_i\left(\bar\partial_{{\bold M}_0}
{\bold E}_{\bold M}(h)\right)f_i-\sum_{i=1}^s
\alpha_i\left({\bold E}_{\bold M}(h)\right)
{\bold E}_{{\bold M}_0}\left(\bar\partial_{{\bold M}_0}g_i\right)\right)
+\sum_{i=1}^s \alpha_i\left({\bold E}_{\bold M}(h)\right)g_i,$$
$${\bold P}^{r+1}_{\bold M}(v)={\bold R}^{r+1}_{\bold M}
\left(v-\sum_{i=1}^l \beta_i\left({\bold E}_{\bold M}(v)\right)
\bar\partial_{\bold M}f_i\right)
+\sum_{i=1}^l \beta_i\left({\bold E}_{\bold M}(v)\right)f_i,$$
and consider the following equality
$${\bold E}_{{\bold M}_0}\left[\bar\partial_{{\bold M}_0}
{\bold P}^r_{{\bold M}_0}
\left({\bold E}_{\bold M}(h)\right)+{\bold P}^{r+1}_{{\bold M}_0}
\left(\bar\partial_{{\bold M}_0} {\bold E}_{\bold M}(h)\right)\right]
-\bar\partial_{\bold M}{\bold P}^r_{\bold M}(h)
-{\bold P}^{r+1}_{\bold M}(\bar\partial_{\bold M}h)
\eqno(\arabic{equation})
\newcounter{PM0-PM}
\setcounter{PM0-PM}{\value{equation}}
\addtocounter{equation}{1}$$
$$={\bold E}_{{\bold M}_0}\left(\bar\partial_{{\bold M}_0}
\left[{\bold R}^r_{{\bold M}_0}
\left({\bold E}_{\bold M}(h)-\sum_{i=1}^l \beta_i
\left(\bar\partial_{{\bold M}_0}{\bold E}_{\bold M}(h)\right)f_i
-\sum_{i=1}^s \alpha_i\left({\bold E}_{\bold M}(h)\right)
\bar\partial_{{\bold M}_0}g_i\right)\right.\right.$$
$$\left.+\sum_{i=1}^s \alpha_i\left({\bold E}_{\bold M}(h)
\right)g_i\right]$$
$$\left.+{\bold R}^{r+1}_{{\bold M}_0}
\left(\bar\partial_{{\bold M}_0}{\bold E}_{\bold M}(h)
-\sum_{i=1}^l \beta_i\left(\bar\partial_{{\bold M}_0}
{\bold E}_{\bold M}(h)\right)\bar\partial_{{\bold M}_0}f_i\right)
+\sum_{i=1}^l \beta_i\left(\bar\partial_{{\bold M}_0}
{\bold E}_{\bold M}(h)\right)f_i\right)$$
$$-\bar\partial_{\bold M}\left[{\bold R}^r_{\bold M}
\left(h-\sum_{i=1}^l \beta_i\left(\bar\partial_{{\bold M}_0}
{\bold E}_{\bold M}(h)\right)f_i
-\sum_{i=1}^s \alpha_i\left({\bold E}_{\bold M}(h)\right)
{\bold E}_{{\bold M}_0}\left(\bar\partial_{{\bold M}_0}g_i\right)\right)
+\sum_{i=1}^s \alpha_i\left({\bold E}_{\bold M}(h)\right)g_i\right]$$
$$-{\bold R}^{r+1}_{\bold M}\left(\bar\partial_{\bold M}h
-\sum_{i=1}^l \beta_i\left({\bold E}_{\bold M}
(\bar\partial_{\bold M}h)\right)\bar\partial_{\bold M}f_i\right)
-\sum_{i=1}^l \beta_i\left({\bold E}_{\bold M}
(\bar\partial_{\bold M}h)\right)f_i$$
$$=\left\{{\bold E}_{{\bold M}_0}\left[\bar\partial_{{\bold M}_0}
{\bold R}^r_{{\bold M}_0}\left({\bold E}_{\bold M}(h)\right)
+{\bold R}^{r+1}_{{\bold M}_0}\left(\bar\partial_{{\bold M}_0}
{\bold E}_{\bold M}(h)\right)\right]
-\bar\partial_{\bold M}{\bold R}^r_{\bold M}(h)
-{\bold R}^{r+1}_{\bold M}\left(\bar\partial_{\bold M}h\right)\right\}$$
$$+\left\{\bar\partial_{\bold M}{\bold R}^r_{\bold M}
\left(\sum_{i=1}^l \beta_i\left(\bar\partial_{{\bold M}_0}
{\bold E}_{\bold M}(h)\right)f_i
\right)+{\bold R}^{r+1}_{\bold M}\left(\bar\partial_{\bold M}
\left(\sum_{i=1}^l \beta_i\left(\bar\partial_{{\bold M}_0}
{\bold E}_{\bold M}(h)\right)f_i\right)\right)\right.$$
$$-\left.{\bold E}_{{\bold M}_0}\left[\bar\partial_{{\bold M}_0}
{\bold R}^r_{{\bold M}_0}
\left(\sum_{i=1}^l \beta_i\left(\bar\partial_{{\bold M}_0}
{\bold E}_{\bold M}(h)\right)f_i\right)+{\bold R}^{r+1}_{{\bold M}_0}
\left(\bar\partial_{{\bold M}_0}\left(\sum_{i=1}^l \beta_i
\left(\bar\partial_{{\bold M}_0}{\bold E}_{\bold M}(h)\right)f_i\right)\right)
\right]\right\}$$
$$+{\bold R}^{r+1}_{\bold M}\bar\partial_{\bold M}
\left(\sum_{i=1}^l \beta_i\left({\bold E}_{\bold M}
(\bar\partial_{\bold M}h)\right)f_i
-\sum_{i=1}^l \beta_i\left(\bar\partial_{{\bold M}_0}
{\bold E}_{\bold M}(h)\right)f_i\right)$$
$$+\left\{\bar\partial_{\bold M}{\bold R}^r_{\bold M}
\left(\sum_{i=1}^s \alpha_i\left({\bold E}_{\bold M}(h)\right)
{\bold E}_{{\bold M}_0}\left(\bar\partial_{{\bold M}_0}g_i\right)\right)
+{\bold R}^{r+1}_{\bold M}\bar\partial_{\bold M}
\left(\sum_{i=1}^s \alpha_i\left({\bold E}_{\bold M}(h)\right)
{\bold E}_{{\bold M}_0}\left(\bar\partial_{{\bold M}_0}g_i\right)\right)\right.$$
$$-\left.{\bold E}_{{\bold M}_0}\left[\bar\partial_{{\bold M}_0}
{\bold R}^r_{{\bold M}_0}
\left(\sum_{i=1}^s \alpha_i\left({\bold E}_{\bold M}(h)\right)
\bar\partial_{{\bold M}_0}g_i\right)
+{\bold R}^{r+1}_{{\bold M}_0}\bar\partial_{{\bold M}_0}
\left(\sum_{i=1}^s \alpha_i\left({\bold E}_{\bold M}(h)\right)
\bar\partial_{{\bold M}_0}g_i\right)\right]\right\}$$
$$-{\bold R}^{r+1}_{\bold M}\bar\partial_{\bold M}
\left[\sum_{i=1}^s \alpha_i\left({\bold E}_{\bold M}(h)\right)
\left({\bold E}_{{\bold M}_0}\left(\bar\partial_{{\bold M}_0}g_i\right)
-\bar\partial_{\bold M}g_i\right)\right]$$
$$+\sum_{i=1}^l \beta_i\left(\bar\partial_{{\bold M}_0}{\bold E}_{\bold M}(h)
-{\bold E}_{\bold M}\left(\bar\partial_{\bold M}h\right)\right)f_i$$
$$+\sum_{i=1}^s \alpha_i\left({\bold E}_{\bold M}(h)\right)
\left({\bold E}_{{\bold M}_0}\left(\bar\partial_{{\bold M}_0}g_i\right)
-\bar\partial_{\bold M}g_i\right)$$
$$=\left\{{\bold E}_{{\bold M}_0}\left[\left({\bold I}
-{\bold H}_{{\bold M}_0}\right)
\left({\bold E}_{\bold M}(h)\right)\right]-\left({\bold I}
-{\bold H}_{\bold M}\right)(h)\right\}$$
$$-\left\{{\bold E}_{{\bold M}_0}\left[\left({\bold I}
-{\bold H}_{{\bold M}_0}\right)
\left(\sum_{i=1}^l \beta_i\left(\bar\partial_{{\bold M}_0}
{\bold E}_{\bold M}(h)\right)f_i
\right)\right]-\left({\bold I}-{\bold H}_{\bold M}\right)
\left(\sum_{i=1}^l \beta_i\left(\bar\partial_{{\bold M}_0}
{\bold E}_{\bold M}(h)\right)f_i\right)\right\}$$
$$-\left\{{\bold E}_{{\bold M}_0}\left[\left({\bold I}
-{\bold H}_{{\bold M}_0}\right)
\left(\sum_{i=1}^s \alpha_i\left({\bold E}_{\bold M}(h)\right)
\bar\partial_{{\bold M}_0}g_i
\right)\right]-\left({\bold I}-{\bold H}_{\bold M}\right)
\left(\sum_{i=1}^s \alpha_i\left({\bold E}_{\bold M}(h)\right)
{\bold E}_{{\bold M}_0}
\left(\bar\partial_{{\bold M}_0}g_i\right)\right)\right\}$$
$$+{\bold R}^{r+1}_{\bold M}
\left(\sum_{i=1}^l \beta_i\left({\bold E}_{\bold M}(\bar\partial_{\bold M}h)
-\bar\partial_{{\bold M}_0}{\bold E}_{\bold M}(h)\right)
\bar\partial_{\bold M}f_i\right)$$
$$+\sum_{i=1}^l \beta_i\left(\bar\partial_{{\bold M}_0}{\bold E}_{\bold M}(h)
-{\bold E}_{\bold M}\left(\bar\partial_{\bold M}h\right)\right)f_i$$
$$-{\bold R}^{r+1}_{\bold M}\bar\partial_{\bold M}
\left[\sum_{i=1}^s \alpha_i\left({\bold E}_{\bold M}(h)\right)
\left({\bold E}_{{\bold M}_0}\left(\bar\partial_{{\bold M}_0}g_i\right)
-\bar\partial_{\bold M}g_i\right)\right]$$
$$+\sum_{i=1}^s \alpha_i\left({\bold E}_{\bold M}(h)\right)
\left({\bold E}_{{\bold M}_0}\left(\bar\partial_{{\bold M}_0}g_i\right)
-\bar\partial_{\bold M}g_i\right).$$
\indent
For the first term of the right hand side of equality above
using equality
$${\bold E}_{{\bold M}_0}{\bold E}_{\bold M}(h)=h,
\eqno(\arabic{equation})
\newcounter{EEh}
\setcounter{EEh}{\value{equation}}
\addtocounter{equation}{1}$$
and applying estimate (\arabic{HM-HM0estimate}) from
Proposition~\ref{HM-HM0} we obtain that for any $\epsilon >0$ there
exists $\delta > 0$ such that for ${\bold M}$ with $|{\cal E}|_k < \delta$
the following estimate holds
$$\left|\left\{{\bold E}_{{\bold M}_0}\left[\left({\bold I}
-{\bold H}_{{\bold M}_0}\right)
\left({\bold E}_{\bold M}(h)\right)\right]-\left({\bold I}
-{\bold H}_{\bold M}\right)(h)\right\}\right|_k
\leq C(k)\epsilon\left(1+ |\rho^{(0)}|_{k+4}\right)^{P(k)}|h|_k.$$
\indent
For the second term of the right hand side of (\arabic{PM0-PM})
using equality
$${\bold E}_{{\bold M}_0}f_i={\bold E}_{\bold M}f_i = f_i,
\eqno(\arabic{equation})
\newcounter{Ef}
\setcounter{Ef}{\value{equation}}
\addtocounter{equation}{1}$$
estimate (\arabic{HM-HM0estimate}) and boundedness of functionals
$\left\{\beta_i\right\}_1^l$ on $C^{1}_{(0,r+1)}\left({\bold M}_0,
{\cal B}\big|_{{\bold M}_0}\right)$,
we obtain that for any $\epsilon >0$ there
exists $\delta > 0$ such that for ${\bold M}$ with $|{\cal E}|_k < \delta$
the following estimate holds
$$\left|{\bold E}_{{\bold M}_0}\left[\left({\bold I}
-{\bold H}_{{\bold M}_0}\right)
\left(\sum_{i=1}^l \beta_i\left(\bar\partial_{{\bold M}_0}
{\bold E}_{\bold M}(h)\right)f_i
\right)\right]-\left({\bold I}-{\bold H}_{\bold M}\right)
\left(\sum_{i=1}^l \beta_i\left(\bar\partial_{{\bold M}_0}
{\bold E}_{\bold M}(h)\right)f_i\right)\right|_k$$
$$\leq C(k)\epsilon\left(1+ |\rho^{(0)}|_{k+4}\right)^{P(k)}|h|_2.$$
\indent
For the third term of the right hand side of (\arabic{PM0-PM})
using equality
$${\bold E}_{\bold M}{\bold E}_{{\bold M}_0}
\left(\bar\partial_{{\bold M}_0}g_i\right)
=\bar\partial_{{\bold M}_0}g_i,
\eqno(\arabic{equation})
\newcounter{Ebarg}
\setcounter{Ebarg}{\value{equation}}
\addtocounter{equation}{1}$$
estimate (\arabic{HM-HM0estimate}) and boundedness of functionals
$\left\{\alpha_i\right\}_1^s$ on $C^{1}_{(0,r)}\left({\bold M}_0,
{\cal B}\big|_{{\bold M}_0}\right)$,
we obtain that for any $\epsilon >0$ there
exists $\delta > 0$ such that for ${\bold M}$ with $|{\cal E}|_k < \delta$
the following estimate holds
$$\left|\left\{{\bold E}_{{\bold M}_0}\left[\left({\bold I}
-{\bold H}_{{\bold M}_0}\right)
\left(\sum_{i=1}^s \alpha_i\left({\bold E}_{\bold M}(h)\right)
\bar\partial_{{\bold M}_0}g_i
\right)\right]-\left({\bold I}-{\bold H}_{\bold M}\right)
\left(\sum_{i=1}^s \alpha_i\left({\bold E}_{\bold M}(h)\right)
{\bold E}_{{\bold M}_0}\left(\bar\partial_{{\bold M}_0}g_i
\right)\right)\right\}\right|_k$$
$$\leq C(k)\epsilon\left(1+ |\rho^{(0)}|_{k+4}\right)^{P(k)}|h|_1.$$
\indent
For the next two terms we use the definition of $\bar\partial_{\bold M}$ as
$$\bar\partial_{\bold M}={\bold r}_{\bold M}\circ\bar\partial
\circ{\bold E}_{\bold M}$$ where ${\bold r}_{\bold M}$ is the operator
of restriction onto ${\bold M}$. Then we obtain that
$${\bold E}_{\bold M}\left(\bar\partial_{\bold M}h\right)
={\bold E}_{\bold M}\circ{\bold r}_{\bold M}
\circ\bar\partial\circ{\bold E}_{\bold M}(h)$$
and
$$\bar\partial_{{\bold M}_0}{\bold E}_{\bold M}(h)={\bold r}_{{\bold M}_0}
\circ\bar\partial\circ{\bold E}_{\bold M}(h).$$
\indent
Using local representation
$$f=\bar\partial\circ{\bold E}_{\bold M}(h)
=\sum_{I,J,K,L}f_{I,J,K}(X,Y,W)dX^I\wedge dY^J\wedge d{\overline W}^K$$
we obtain locally
$$\begin{array}{ll}
{\bold E}_{\bold M}\circ{\bold r}_{\bold M}\circ\bar\partial
\circ{\bold E}_{\bold M}(h)
=\sum_{I,J,K,L}f_{I,J,K}(\phi(Y,W),Y,W)[d\phi(Y,W)]^I\wedge dY^J\wedge
d{\overline W}^K,\vspace{0.1in}\\
{\bold r}_{{\bold M}_0}\circ\bar\partial\circ{\bold E}_{\bold M}(h)
=\sum_{I,J,K,L}f_{I,J,K}(\phi^{(0)}(Y,W),Y,W)[d\phi^{(0)}(Y,W)]^I
\wedge dY^J\wedge d{\overline W}^K,
\end{array}
\eqno(\arabic{equation})
\newcounter{Restrictions1}
\setcounter{Restrictions1}{\value{equation}}
\addtocounter{equation}{1}$$
where $\left\{\phi_j\right\}_1^m$ and $\left\{\phi_j^{(0)}\right\}_1^m$
are defining functions of ${\bold M}$ and ${\bold M}_0$ respectively.\\
\indent
Using then estimate (\arabic{phi0-phiestimate}) from Lemma~\ref{phi-phi0}
for expressions in (\arabic{Restrictions1}) we obtain the estimate
$$\left|\bar\partial_{{\bold M}_0}{\bold E}_{\bold M}(h)
-{\bold E}_{\bold M}\left(\bar\partial_{\bold M}h\right)\right|_{1}$$
$$\leq C|h|_3|\phi - \phi^{(0)}|_2
\leq C|g|_2\left(1+ |\rho^{(0)}|_{3}\right)^{P}|h|_3.$$
\indent
Using the last estimate and the estimate
$$\sup_{i=1}^{l}|\bar\partial_{\bold M}f_i|_k <C$$
for $k \leq p-5$ we obtain the existence for $k \leq p-5$ and any
$\epsilon > 0$ of a $\delta > 0$ such that for ${\bold M}$ with
$|{\cal E}|_k < \delta$ the following estimate holds
$$\Biggl|{\bold R}^{r+1}_{\bold M}
\left(\sum_{i=1}^l \beta_i\left({\bold E}_{\bold M}(\bar\partial_{\bold M}h)
-\bar\partial_{{\bold M}_0}{\bold E}_{\bold M}(h)\right)
\bar\partial_{\bold M}f_i\right)
+\sum_{i=1}^l \beta_i\left(\bar\partial_{{\bold M}_0}{\bold E}_{\bold M}(h)
-{\bold E}_{\bold M}\left(\bar\partial_{\bold M}h\right)\right)f_i\Biggr|_k$$
$$\leq C\epsilon \left(1+|\rho^{(0)}|_{3}\right)^{P}|h|_3.$$
\indent
To estimate the last two terms of the right hand side of (\arabic{PM0-PM})
we again use the definition of $\bar\partial_{\bold M}$ and representations
$$\bar\partial_{\bold M}g_i={\bold r}_{\bold M}\circ\bar\partial g_i$$
$$=\sum_{I,J,K,L}f_{I,J,K}(\phi(Y,W),Y,W)[d\phi(Y,W)]^I\wedge dY^J\wedge
d{\overline W}^K \in C^{p-5}_{(0,r)}\left({\bold M}_0,
{\cal B}\big|_{{\bold M}_0}\right),$$
and
$${\bold E}_{\bold M}\left(\bar\partial_{{\bold M}_0}g_i\right)
={\bold E}_{\bold M}\circ{\bold r}_{{\bold M}_0}\circ\bar\partial g_i
=\sum_{I,J,K,L}f_{I,J,K}(\phi^{(0)}(Y,W),Y,W)[d\phi^{(0)}(Y,W)]^I\wedge
dY^J\wedge d{\overline W}^K.$$
\indent
Using then estimate (\arabic{phi0-phiestimate}) from Lemma~\ref{phi-phi0}
and estimates
$$\sup_{i=1}^s|g_i|_{k} \leq C(k)\hspace{0.05in}\mbox{for}\hspace{0.05in}
k \leq p-4,$$
$$|\phi|_{k} \leq C(k)\left(1+ |\rho^{(0)}|_{k+1}\right)\hspace{0.05in}
\mbox{for}\hspace{0.05in}|g|_{k}\hspace{0.05in}\mbox{small enough},$$
we obtain for $k \leq p-6$ the estimate
$$\left|{\bold E}_{\bold M}\left(\bar\partial_{{\bold M}_0}g_i\right)
-\bar\partial_{\bold M}g_i\right|_{k}
\leq C(k)|g_i|_{k+2}\left(1+ |\rho^{(0)}|_{k+2}\right)^{P(k)}|g|_{k+1}.$$
\indent
Using the last estimate we obtain the existence for $k \leq p-7$ and
for any $\epsilon > 0$ of a $\delta > 0$ such that for ${\bold M}$
with $|{\cal E}|_{k+2} < \delta$ the following estimate holds
$$\Biggl|{\bold R}^{r+1}_{\bold M}\bar\partial_{\bold M}
\left[\sum_{i=1}^s \alpha_i\left({\bold E}_{\bold M}(h)\right)
\left({\bold E}\left(\bar\partial_{{\bold M}_0}g_i\right)
-\bar\partial_{\bold M}g_i\right)\right]
-\sum_{i=1}^s \alpha_i\left({\bold E}_{\bold M}(h)\right)
\left({\bold E}\left(\bar\partial_{{\bold M}_0}g_i\right)
-\bar\partial_{\bold M}g_i\right)\Biggr|_k$$
$$\leq C(k)\epsilon\left(1+ |\rho^{(0)}|_{k+3}\right)^{P(k)}|h|_1.$$
\indent
Let now ${\bold B}_{{\bold M}_0}$ be the inverse to
${\bold F}_{{\bold M}_0}
= \bar\partial_{{\bold M}_0}{\bold P}^r_{{\bold M}_0}
+ {\bold P}^{r+1}_{{\bold M}_0}\bar\partial_{{\bold M}_0}$, which exists
according to Proposition~\ref{M0Formula}
and let $B=\sup_{k=1,\dots,p-7}\left|{\bold B}_{{\bold M}_0}\right|_k$.
Combining all the estimates above we obtain that
there exists $\delta >0$ such that for ${\bold M}$ with
$|{\cal E}|_p < \delta$ the following estimate holds for any
$h \in C^k_{(0,r)}\left({\bold M},{\cal B}\big|_{\bold M}\right),
\hspace{0.05in}(k \leq p-7)$:
$$\left|{\bold E}_{{\bold M}_0}\left[\bar\partial_{{\bold M}_0}
{\bold P}^r_{{\bold M}_0}\left({\bold E}_{\bold M}(h)\right)
+{\bold P}^{r+1}_{{\bold M}_0}
\left(\bar\partial_{{\bold M}_0} {\bold E}_{\bold M}(h)\right)\right]
-\bar\partial_{\bold M}{\bold P}^r_{\bold M}(h)
-{\bold P}^{r+1}_{\bold M}(\bar\partial_{\bold M}h)\right|_k
< \frac{B^{-1}}{4}\cdot |h|_k.$$
\indent
Applying then operator
${\bold C}_{\bold M}={\bold E}_{{\bold M}_0}\circ{\bold B}_{{\bold M}_0}
\circ{\bold E}_{\bold M}$ to the form
$${\bold E}_{{\bold M}_0}\left[\bar\partial_{{\bold M}_0}
{\bold P}^r_{{\bold M}_0}\left({\bold E}_{\bold M}(h)\right)
+{\bold P}^{r+1}_{{\bold M}_0}
\left(\bar\partial_{{\bold M}_0} {\bold E}_{\bold M}(h)\right)\right]
-\bar\partial_{\bold M}{\bold P}^r_{\bold M}(h)
-{\bold P}^{r+1}_{\bold M}(\bar\partial_{\bold M}h)$$
and using equalities
$$\begin{array}{ll}
{\bold E}_{\bold M}\circ{\bold E}_{{\bold M}_0}={\bold I},\vspace{0.1in}\\
{\bold B}_{{\bold M}_0}\circ{\bold F}_{{\bold M}_0}={\bold I},
\end{array}$$
we obtain the following estimate
$$\left|{\bold I}-{\bold C}_{\bold M}{\bold F}_{\bold M}\right|_k < \frac{1}{4}.
\eqno(\arabic{equation})
\newcounter{I-CF}
\setcounter{I-CF}{\value{equation}}
\addtocounter{equation}{1}$$
\indent
From the estimate above we obtain the existence of an inverse operator
${\bold D}_{\bold M}$ such that
$$\begin{array}{ll}
{\bold D}_{\bold M}\circ{\bold F}_{\bold M}={\bold I},\vspace{0.1in}\\
\left|{\bold D}_{\bold M}\right|_k \leq \frac{4}{3}
\left|{\bold B}_{{\bold M}_0}\right|_k.
\end{array}$$
\indent
Therefore for any
$h \in C^k_{(0,r)}\left({\bold M},{\cal B}\big|_{\bold M}\right)$
we have
$$h = {\bold D}_{\bold M}\bar\partial_{\bold M}{\bold P}^r_{\bold M}(g)
+ {\bold D}_{\bold M}{\bold P}^{r+1}_{\bold M}(\bar\partial_{\bold M}g)$$
$$= \bar\partial_{\bold M}{\bold P}^r_{\bold M}
{\bold D}_{\bold M}^2\bar\partial_{\bold M}{\bold P}^r_{\bold M}(g)
+ {\bold D}_{\bold M}{\bold P}^{r+1}_{\bold M}(\bar\partial_{\bold M}g)$$
$$= \bar\partial_{\bold M}{\bold Q}^r_{\bold M}(g)
+ {\bold Q}^{r+1}_{\bold M}(\bar\partial_{\bold M}g),$$
with ${\bold Q}^r_{\bold M} = {\bold P}^r_{\bold M}{\bold D}^2_{\bold M}
\bar\partial_{\bold M}{\bold P}^r_{\bold M}$ and
${\bold Q}^{r+1}_{\bold M} ={\bold D}_{\bold M}{\bold P}^{r+1}_{\bold M}$.\qed

\end{document}